\documentclass[11pt]{article}
\usepackage{latexsym}
\usepackage{amsfonts,amssymb,amsmath,mathtools}
\usepackage{bbold}

\newtheorem{thm}{Theorem}[section]
\newtheorem{cor}[thm]{Corollary}
\newtheorem{lem}[thm]{Lemma}
\newtheorem{pro}[thm]{Proposition}
\newtheorem{defn}[thm]{Definition}

\newcommand{\ov }{\overline }
\newcommand{\minus}{\smallsetminus}
\newcommand{\x}{\hspace{-0.025in}\times\hspace{-0.025in}}
\newcommand{\e}{\varepsilon}
\newcommand{\0}{\varnothing}

\title{Some properties of Higman-Thompson monoids \\      
    and digital circuits }

\author{J.C.\ Birget}
 
\date{\scriptsize{      
19 iii 2024              
      }           }

\begin{document}
\maketitle

\begin{abstract}
We define various monoid versions of the Thompson group $V$, and prove
connections with monoids of acyclic digital circuits.
We show that the monoid $M_{2,1}\,$ (based on partial functions) is not
embeddable into Thompson's monoid ${\sf tot}M_{2,1}$, but that 
${\sf tot}M_{2,1}$ has a submonoid that maps homomorphically onto $M_{2,1}$.
This leads to an efficient completion algorithm for partial functions and
partial circuits.
We show that the union of partial circuits with disjoint domains is an 
element of $M_{2,1}$, and conversely, every element of $M_{2,1}$ can be 
decomposed efficiently into a union of partial circuits with disjoint 
domains. 
\end{abstract}

\section{Introduction}

We consider various monoids that generalize Richard Thompson's group $V$. 
The group $V$ (which is also called $G_{2,1}$ or $V_2$) was introduced by 
Thompson \cite{Th60s, McKTh} (see also \cite{CFP}), and he proved that $V$
is an infinite finitely presented simple group.  In fact, as Matt Brin 
pointed out \cite{Brin_PC, Brin_lect, stackex}, Thompson initially derived 
$V$ as the group of units of a monoid (in the mid 1960s); 
he did not publish anything about this monoid, but talked about it in 
lectures (notes by Brin \cite{Thompson_lect2004, Thompson_lect2008}). 
Thompson's monoid, was rediscovered in \cite{BiThompsMonV3, BiThompsMon}, 
where Thompson's monoid is called ${\sf tot}M_{2,1}$. 
A partial-function version of the monoid was also introduced in 
\cite{BiThompsMonV3, BiThompsMon}, called $M_{2,1}\,$ (it had apparently not 
been considered by Thompson). In \cite{BiThompsMonV3, BiThompsMon} it was 
proved that $M_{2,1}$ is finitely generated and congruence-simple, that its 
group of units is $V$, and that its word problem over a finite generating set 
is in {\sf P}. 

In his lectures \cite{Thompson_lect2004, Thompson_lect2008} Thompson 
outlined a proof (recently reconstructed by Brin \cite{Brin_PC}) that 
${\sf tot}M_{2,1}$ is finitely presented. By using similar ideas Brin also 
proved that $M_{2,1}$ is finitely presented  (but Brin's proofs of finite 
presentation remain unpublished). More recently a proof of finite 
presentation of ${\sf tot}M_{2,1}\,$ (and more generally of the 
$n$-dimensional version $\,n\,{\sf tot}M_{2,1}$), was posted by de Witt and 
Elliott \cite{deWittElliott}.

\smallskip

Most of the time we will actually work with the Higman-Thompson monoids 
$M_{k,1}$ and ${\sf tot}M_{k,1}$, and the Higman-Thompson group $G_{k,1}$
\cite{Hig74}, although we often just say ``Thompson monoids'' and 
``Thompson group''.

\medskip

{\sf Remark about the name ``Thompson monoid'':} The Thompson group 
$F$ contains a submonoid $F^+$ such that $F$ is the group of fractions is 
$F^+$. In the literature $F^+$ is usually called the ``Thompson monoid'' 
or the ``positive submonoid'' of $F$. The present paper does not explicitly
use $F^+$, and ``Thompson monoids'' will refer to various monoids that
generalize $V$, replacing bijections by various functions. 

\bigskip

\noindent {\bf Acknowledgement:} Section 8 was motivated by discussions
with Matt Brin \cite{Brin_PC}.

\newpage 

\noindent {\bf Summary of results}

\smallskip

\noindent
We prove new connections between $M_{k,1}$, ${\sf tot}M_{k,1}$, other
Thompson monoids, and acyclic digital circuits:

\smallskip

\noindent $\bullet \ $ Sections 2, 3, and 4:  We define various 
pre-Thompson and Thompson monoids, and circuits (partial circuits and 
boolean circuits).
The Thompson monoids are obtained from pre-Thompson monoids by applying
congruences.
All these monoids (including the monoids of input-output
functions of circuits) have a generating set of the form $\Gamma \cup \tau$,
where $\Gamma$ is finite, and $\tau$ is the infinite set of bit 
transpositions; for circuits, $\Gamma$ is the set of gates.
We show that circuits are closely related to Thompson monoids, and that
important concepts from circuits (input-length, size, depth) can be 
defined in all Thompson monoids. 

\smallskip

\noindent $\bullet \ $ Section 5: The pre-Thompson monoid
${\cal RM}_2^{\sf fin}\,$ (which maps onto $M_{2,1}$) is finitely generated.

\smallskip

\noindent $\bullet \ $ Section 6: Certain pre-Thompson and Thompson monoids,
including the monoid of input-output functions of circuits, are {\sl not} 
finitely generated.

\smallskip

\noindent $\bullet \ $ Section 7: Most Thompson monoids are 
congruence-simple. An exception is the monoid of input-output functions of
boolean circuits. 

\smallskip

\noindent $\bullet \ $ Section 8: The Thompson monoid $M_{2,1}$ (based on
partial functions) is not embeddable homomorphically into Thompson's
monoid ${\sf tot}M_{2,1}$ (based on total functions). 
Similarly, the monoid of input-output functions of partial circuits is not 
embeddable homomorphically into the monoid of input-output functions of 
boolean circuits; in other words, partial circuits cannot be completed to 
boolean circuits in a homomorphic way (where the completion of the composite 
would be the composite of the completions). 
  
\smallskip

\noindent $\bullet \ $ Section 9: Although $M_{2,1}$ is not embeddable into
${\sf tot}M_{2,1}$ and not a homomorphic image of ${\sf tot}M_{2,1}$ (by 
congruence-simplicity), there is a submonoid of ${\sf tot}M_{2,1}$ that maps
homomorphically onto $M_{2,1}$. This is an inverse completion, which can be 
used to give efficient completion algorithms (for $M_{2,1}$ and for partial
circuits).

\smallskip

\noindent $\bullet \ $ Section 10: Any element of $M_{2,1}$ can be computed 
by (a finite set of) partial circuits with disjoint domains.
Conversely, the union of a finite set of partial circuits with disjoint 
domains is an element of $M_{2,1}$.  

In conclusion, the Thompson monoids and circuits are, in essence, the same
thing.

\section{Preliminary monoids }

Before defining various Thompson monoids we need some definitions, that are
fairly standard (see also \cite{BerstelPerrin, BiTh, Bi1wPerm, 
BiThompsMonV3, BiOntoTh}). 
We first introduce monoids of right-ideal morphisms; the Thompson and
Higman-Thompson monoids are then obtained from these {\em ``pre-Thompson'' 
monoids} by applying a congruence.

\smallskip 

\begin{defn} \label{DEFfunction} {\bf (function).}

A function is a triple $f = (X, r, Y)$, where $r$ $\subseteq$
$X \x Y$ is a relation that has the function property:
 \ $\big(\forall \,(x_1,y_1), (x_2,y_2) \in r \big)[\,x_1 = x_2$
$\,\Rightarrow\,$  $y_1 = y_2\,]$.

\smallskip

The {\em domain} of $f$ is $\,{\rm Dom}(f) = {\rm Dom}(r)$
$=$  $\{ x: (\exists y \in Y)[(x,y) \in r]\}$ $\,\subseteq\,$ $X$, and the 
{\em image} (a.k.a.\ {\em range}) is $\,{\rm Im}(f) = {\rm Im}(r)$
$=$  $\{y: (\exists x \in X)[(x,y) \in r]\}$ $\,\subseteq\,$ $Y$.
Thus, $r \subseteq {\rm Dom}(f) \x {\rm Im}(f) \subseteq X \x Y$.

The sets $X$ and $Y$ are called the {\em source}, respectively the
{\em target}, of $f = (X, r, Y)$.
\end{defn}
Here we usually take $\,X = Y = A^*$, or $\,X = Y = A^{\omega}$, or
$\,X = Y = A^*0^{\omega}$, where $A$ is a finite alphabet such that
$\{0,1\} \subseteq A$, and $\omega$ is the set of natural integers (as an
ordinal).

The composition of two functions $f_1 = (X_1, r_1, Y_1)$,
$f_2 = (X_2, r_2, Y_2)$, is defined by $\,(.)f_1 \circ f_2 = $
$(X_1, (.)r_1 \circ r_2, Y_2)$, where $\,(.)r_1 \circ r_2 =$
$\{(x,y) \in X_1 \x Y_2 :$  $(\exists z \in Y_1 \cap X_2)[\,(x,z) \in r_1$
$\&$  $(z,y) \in r_2\,]\,\}$.
Then $\,{\rm Dom}((.)f_1 \circ f_2) = $
${\rm Dom}(f_1) \cap f_1^{-1}({\rm Dom}(f_2))$,  and
$\,{\rm Im}((.)f_1 \circ f_2) = $
${\rm Im}(f_2) \cap f_2({\rm Im}(f_1))$.
We often write $(.)f$ for $f$ to indicate that the function is applied to
the right, and we write $f(.)$ when the function is applied to the left.

We follow the set-theoretic view of functions, not the category-theoretic
one. (In category theory, most commonly, functions are treated as total,
and composition is undefined if $Y_1 \ne X_2$.)
We will deviate from the set-theoretic view of functions in one way:
We will usually make functions act on the {\em left} of the argument, and we
write $f(.)$ to indicate this; we sometimes let functions act on the right,
and then we write $(.)f$.

A function $\,f: X \to Y$ is {\em total} iff ${\rm Dom}(f) = X$ and
$X \ne \0$.

The {\em empty function} $\,(X, \0, Y)$ is denoted by $\theta\,$ when $X$
and $Y$ are clear from the context.

\medskip

We use finite alphabets, typically of the form $A = \{0,1,\ldots,k\!-\!1\}$ 
with $k \ge 1$, and the free monoid $A^*$, consisting of all strings (i.e.,
finite sequences of elements of $A$), including the {\em empty string} $\e$.
We also consider the set $A^{\omega}$ of $\omega$-sequences over  $A$. 
For a string $x \in A^*$ the length is denoted by $|x|$. Concatenation of 
$x, y \in A^*$ is denoted by $xy$ (juxtaposition) or by $x\cdot y$. Sets 
$X, Y \subseteq A^*$ are concatenated as $XY = \{xy: x \in X, y \in Y\}$.
For $x, y \in A^*$, $\,x$ is a {\em prefix} of $y$ iff $y = x u$ for some 
$u \in A^*$; this is denoted by $x \le_{\rm pref} y$. 
We write $x <_{\rm pref} y$ iff $x \le_{\rm pref} y$
and $x \ne y$. The relation $\le_{\rm pref}$ is a partial order. 
Strings $x, y \in A^*$ are {\em prefix-comparable} iff 
$x \le_{\rm pref} y$ or $y \le_{\rm pref} x$; this is denoted by
$\,x \|_{\rm pref} \,y$. A set $P \subseteq A^*$ is a {\em prefix code} iff 
no element of $P$ is a prefix of another element of $P$, i.e., different 
elements of $P$ are prefix-incomparable \cite{BerstelPerrin}. 
A {\em maximal prefix code} of $A^*$ is a prefix code that is 
$\subseteq$-maximal among the prefix codes of $A^*$.
A set $R \subseteq A^*$ is a {\em right-ideal} of $A^*$ iff $R = R\,A^*$; 
this is equivalent to 
$\,(\forall x \in R)(\forall u \in A^*)[\, x \le_{\rm pref} u$
$\,\Rightarrow\, u \in R\,]$. 
The intersection (or union) of any set of right ideals is a right ideal.
We say that two sets $X, Y$ intersect iff $X \cap Y \ne \0$.
An {\em essential right ideal} of $A^*$ is a right ideal that intersects 
every right ideal of $A^*$.
The following is easy to prove: $\,R$ is a right-ideal of $A^*$ iff there 
exists a (unique) prefix code $P$ of $A^*$ such that $R = P\,A^*$; $\,P$ is 
called the {\em prefix code of} $R$. Moreover, $R$ is essential iff the
prefix code of $R$ is maximal. 

A {\em finitely generated right ideal} is defined to be a right ideal whose
prefix code is finite. The intersection (or union) of any two finitely 
generated right ideals is a finitely generated right ideal (see 
\cite{BiTh}).

As we saw already, in this paper, ``function'' always means 
{\em partial function}. 

A {\em right-ideal morphism} of $A^*$ is a function $f$: $A^* \to A^*$ such
that ${\rm Dom}(f)$ is a right ideal, and 
$\,(\forall x \in {\rm Dom}(f))(\forall u \in A^*)[\,f(xu) = f(x) \ u\,]$. 
The prefix code of the right ideal ${\rm Dom}(f)$ is called the {\em domain
code} of $f$, and denoted  by ${\rm domC}(f)$. One easily proves that if
$f$ is a right-ideal morphism then ${\rm Im}(f)$ is a right ideal; the prefix
code of ${\rm Im}(f)$ is called the {\em image code} of $f$, and denoted  by
${\rm imC}(f)$. The right-ideal morphism $f$ is uniquely determined by the
restriction $f|_{{\rm domC}(f)}$; this restriction is called the {\em table
of} $f$. In general, ${\rm imC}(f)$ $\subseteq f\big({\rm domC}(f)\big)$, 
but that this inclusion can be strict (as will be discussed later). 

It is easy to prove that if $R$ is a right ideal and $f$ is a right-ideal 
morphism then $f(R)$ and $f^{-1}(R)$ are right ideals. Moreover, can prove 
(see \cite{BiTh, BiThompsMonV3}) that for any two functions $f_2, f_1:$ 
 \ \ ${\rm Dom}(f_2 \circ f_1(.)) =$
$f_1^{-1}\big({\rm Im}(f_1) \cap {\rm Dom}(f_2)\big)$, \ and
 \ ${\rm Im}(f_2 \circ f_1(.)) =$
$f_2\big({\rm Im}(f_1) \cap {\rm Dom}(f_2)\big)$.  It follows that the 
composite $f_2 \circ f_1(.)$ of two right-ideal morphisms $f_2, f_1$ of 
$A^*$ is a right-ideal morphism.

\begin{defn} \label{Def_RMfin} 
 \ The monoid ${\cal RM}^{\sf fin}\,$ (also called ${\cal RM}_k^{\sf fin}$)
consists of all right-ideal morphisms of $A^*$ with {\em finite} domain 
code, with function composition as the monoid operation.

When $A = \{0,1\}$, the above monoid is called ${\cal RM}^{\sf fin}_2$.
\end{defn} 
Equivalently, $f \in {\cal RM}^{\sf fin}$ iff $f$ is a right-ideal morphism
whose table $\,f|_{{\rm domC}(f)}\,$ is finite. 

It follows from \cite[{\small Thm.\ 4.5}]{BiThompsMonV3} that if $f_2, f_1$
are right-ideal morphisms with finite domain codes then $f_2 \circ f_1(.)$ 
also has a finite domain code. Hence ${\cal RM}^{\sf fin}$ is closed under 
composition, so ${\cal RM}^{\sf fin}$ is indeed a monoid.

\bigskip 

\noindent Let is introduce a few more monoids.

\begin{defn} \label{DEFtotRM}\!\!.
 
A right-ideal morphism $f: A^* \to A^*$ is called {\em total} iff 
$\,{\rm domf}(f)$ is a {\em maximal} prefix code of $A^*$; equivalently,
${\rm Dom}(f)$ is an {\em essential} right ideal.

A right-ideal morphism $f$ is called {\em surjective} iff 
$\,{\rm imf}(f)$ is a {\em maximal} prefix code of $A^*$, i.e., 
${\rm Im}(f)$ is an {\em essential} right ideal.

We define 

\bigskip

 \ \ \  \ \ \ ${\sf tot}{\cal RM}^{\sf fin} \ = \ $
$\{f \in {\cal RM}^{\sf fin} :\,f\,$ {\rm is total} \};

\medskip

 \ \ \  \ \ \ ${\sf sur}{\cal RM}^{\sf fin} \ = \ $ 
$\{f \in {\cal RM}^{\sf fin} :\,f\,$ {\rm is surjective} \};

\medskip

 \ \ \  \ \ \ ${\sf inj}{\cal RM}^{\sf fin} \ = \ $
$\{f \in {\cal RM}^{\sf fin} :\,f\,$ {\rm is injective} \}.
\end{defn}
One proves easily that these are indeed monoids under composition. 

Since the above monoids are submonoids of ${\cal RM}^{\sf fin}$ we can 
consider intersections of any of these monoids, e.g., 
${\sf tot inj}{\cal RM}^{\sf fin}$, or 
${\sf tot sur inj}{\cal RM}^{\sf fin}$.  The latter maps homomorphically onto
the Thompson group $V$ and was studies in \cite{BiOntoTh} under the name 
${\sf riAut}(k)$.

\medskip

\noindent {\bf Remark on the words ``total'' and ``surjective'':} 
Here, total or surjective right-ideal morphisms of $A^*$ are usually not 
total, respectively surjective, as functions on $A^*$. However we will see 
later, when ${\sf tot}M_{k,1}$ and ${\sf sur}M_{k,1}$ are introduced, that 
they are total, respectively surjective, as functions on $A^{\omega}$. 

\bigskip

\noindent The following monoids have close connections with acyclic boolean 
circuits.

\begin{defn} \label{DEFlepfl} {\bf ({\sf plep}, {\sf tlep}, {\sf pfl}, and
{\sf tfl}).}

\smallskip

\noindent Let $f: A^* \to A^*$ be a right-ideal morphism. 

\smallskip

\noindent $\bullet$ \ $f$ is called {\em partial length-equality preserving}, 
or {\sf plep}, \ iff \ for all $\,x_1, x_2 \in {\rm Dom}(f):$ 

\smallskip

 \ \ \  \ \ \ $|x_1| = |x_2| \ \ \Rightarrow \ \ |f(x_1)| = |f(x_2)|$. 

\medskip

\noindent $\bullet$ \ $f$ is called {\em total length-equality preserving}, 
or {\sf tlep}, \ iff \ $f$ is total and {\sf plep}.

\smallskip

\noindent $\bullet$ \ $f$ is called {\em partial fixed-length input}, or 
{\sf pfl}, \ iff \ $f$ is {\sf plep} and
 $\,{\rm domC}(f) \subseteq \{0,1\}^m\,$ for some $m \ge 0$.

\smallskip

\noindent $\bullet$ \ $f$ is called {\em total fixed-length input}, or  
{\sf tfl}, \ iff \ $f$ is {\sf tlep} and $\,{\rm domC}(f) = \{0,1\}^m\,$ for
some $m \ge 0$.

\smallskip

\noindent $\bullet$ \ We define

\medskip

 \ \ \  \ \ \ ${\sf plep}{\cal RM}^{\sf fin} \ = \ $
$\{f \in {\cal RM}^{\sf fin} : \,f\,$ {\rm is} {\sf plep} \};

\medskip

 \ \ \  \ \ \ ${\sf tlep}{\cal RM}^{\sf fin} \ = \ $
$\{f \in {\cal RM}^{\sf fin} : \,f\,$ {\rm is} {\sf tlep} \};

\medskip

 \ \ \  \ \ \ ${\sf pfl}{\cal RM}^{\sf fin} \ = \ $
$\{f \in {\cal RM}^{\sf fin} : \,f\,$ {\rm is} {\sf pfl} \};

\medskip

 \ \ \  \ \ \ ${\sf tfl}{\cal RM}^{\sf fin} \ = \ $
$\{f \in {\cal RM}^{\sf fin} : \,f\,$ {\rm is} {\sf tfl} \}.  
\end{defn}
It is easy to show that these are monoids under composition, and that if 
$f$ is {\sf pfl} then 

\smallskip

 \ \ \  \ \ \ $f\big({\rm domC}(f)\big) \,=\, {\rm imC}(f)$ 
$\,\subseteq\, \{0,1\}^n$ \ \ for some $n \ge 0$. 

\medskip

\noindent Since the monoids ${\sf tfl}{\cal RM}^{\sf fin}$ and 
${\sf pfl}{\cal RM}^{\sf fin}$ are strongly related to acyclic circuits, we 
usually consider ${\sf tfl}{\cal RM}^{\sf fin}_2$ and
${\sf pfl}{\cal RM}^{\sf fin}_2$, i.e., we use the alphabet $A = \{0,1\}$. 

\bigskip

\bigskip

\noindent {\bf Summary of pre-Thompson monoids}

\medskip

\noindent We introduced the following monoids of right-ideal morphisms of 
$A^*$, with $|A| = k \ge 2:$

\smallskip

\noindent ${\cal RM}^{\sf fin}$ \ --- \ the monoid of all right-ideal 
morphisms with finite domain-code; also called ${\cal RM}_k^{\sf fin}$.

\smallskip

\noindent  ${\sf tot}{\cal RM}^{\sf fin}$ \ --- \ the monoid of right-ideal
morphisms whose domain-code is a maximal finite prefix code;

\noindent also called ${\sf tot}{\cal RM}_k^{\sf fin}$.

\smallskip

\noindent ${\sf sur}{\cal RM}^{\sf fin}$ \ --- \ the monoid of right-ideal
morphisms whose domain-code is finite, and whose image-code is a maximal
prefix code.

\smallskip

\noindent ${\sf inj}{\cal RM}^{\sf fin}$ \ --- \ the monoid of injective 
right-ideal morphisms whose domain-code is finite.

\medskip

\noindent Intersections of the above, in particular
${\sf totsurinj}{\cal RM}^{\sf fin}$.

\medskip

\noindent ${\sf tfl}{\cal RM}^{\sf fin}_2$ \ --- \ the monoid of right-ideal
morphisms whose domain-code is $\{0,1\}^m$ for some $m\,$ (depending on the
morphism), and whose image-code is a subset of $\{0,1\}^n$ for some $n\,$
(depending on the morphism).

\smallskip

\noindent ${\sf pfl}{\cal RM}^{\sf fin}_2$ \ --- \ the monoid of right-ideal
morphisms whose domain-code is a subset of $\{0,1\}^m$ for some $m\,$
(depending on the morphism), and whose image-code is a subset of $\{0,1\}^n$
for some $n\,$ (depending on the morphism).

\section{Circuits and monoids}

The main motivation for the monoids ${\sf tfl}{\cal RM}^{\sf fin}_2$ and
${\sf pfl}{\cal RM}^{\sf fin}_2$ is their connection to classical acyclic 
boolean circuits and partial circuits. Acyclic circuits are combinational,
so they have no synchronization issues.

Connections between Thompson monoids (and groups), and acyclic circuits 
were previously studied in \cite{BiCoNP, Bi1wPerm, BiFact, BinG}.

\subsection{Circuits }

An {\em acyclic boolean circuit} is a directed acyclic graph ({\sc dag}),
together with a vertex labeling that associates an input variable $x_i$ with 
each source vertex, an output variable $y_j$ with each sink vertex, and a 
gate with each interior vertex.
The classical {\sf not}, {\sf and}, {\sf or}, and {\sf fork} gates are used; 
{\sf not} has fan-in 1 and fan-out 1, {\sf and} and {\sf or} have fan-in 2 
and fan-out 1, {\sf fork} has fan-in 1 and fan-out 2; a source vertex has 
fan-in 0 and fan-out 1, and a sink vertex has fan-in 1 and fan-out 0. Every 
vertex has in-degree and out-degree that matches the fan-in, respectively 
fan-out, of its gate.

A {\em partial acyclic circuit} (more briefly, a {\em partial circuit}) is
defined in the same way as an acyclic boolean circuit, except that now the
set of gates is enlarged by a finite set of additional {\em partial gates} 
in ${\sf pfl}{\cal RM}^{\sf fin}_2$. 
In a partial circuit, if one or more gates produce an undefined output for
a given input, then the whole output of the circuit is undefined for that
input; so partial circuits do not use a kind of multi-valued logic. 
(Physically, an ``undefined'' input or output value of a circuit is any 
physical signal that does not correspond unambiguously to a boolean value 0 
or 1; this could e.g.\ be some ``intermediate'' signal, an undecodable 
signal, or fluctuating signals.)
Subsection 3.2 gives more information on partial circuits; the present 
subsection focuses on boolean circuits. 

\smallskip

It was proved by Hoover, Klawe, and Pippenger \cite{HKP} (see also
\cite[{\small Theorem 1.9.1}]{CloteKran}) that if {\sf and}, {\sf or},
{\sf not} are allowed to have unbounded fan-out (let $\cal G$ be this
infinite set of gates), then every circuit $C$ made from gates in $\cal G$
is equivalent (regarding its input-output function) to a circuit $C'$ with
gates from a fixed finite set of gates $\cal G'$ (equivalently, a set of
gates of bounded fan-in and bounded fan-out), such that
$\,{\rm size}(C') \le c \ {\rm size}(C)$, and
$\,{\sf depth}(C') \le c \ {\sf depth}(C)\,$ (for some constant $c \ge 1$).  

\smallskip

Two circuits $C_1$, $C_2$ are called {\em isomorphic} iff the {\sc dag}s
of $C_1$ and $C_2$ are isomorphic, and the vertices of $C_1$ and $C_2$ that
correspond to each other by the {\sc dag} isomorphism are labeled by the 
same gate, or the same input variable, or the same output variable.
This is denoted by $\,C_1 \,\widehat{=}\, C_2$.
 \ Two isomorphic circuits can be different in their geometric ``layouts''. 
E.g., wire permutations, braiding and knotting of wires, could be different,
but in such a way that the connection relation between vertices are the same.

\smallskip

As input-output functions of the {\em gates}, applied at position $j$ in an
input bit-string, are defined as follows (where $u, v \in \{0,1\}^*$, 
$ \ j  = |u| + 1$):

\medskip

 \ \ \ ${\sf not}_j: \ \ u \, x_j \,v \ \longmapsto$ 
$ \ u \, \overline{x}_j \,v$;

\medskip

 \ \ \ ${\sf and}_{j,j+1}: \ \ u \, x_j x_{j+1} \,v \ \longmapsto \ $
    $u \ (x_j \wedge x_{j+1}) \ v$;

\medskip

 \ \ \ ${\sf or}_{j, j+1}: \ \ u \, x_j x_{j+1} \,v \ \longmapsto$
     $ \ u \ (x_j \vee x_{j+1}) \ v$;

\medskip

 \ \ \ ${\sf fork}_j: \ \ u \, x_j \,v \ \longmapsto$ 
$ \ u \, x_j \, x_j \,v$.

\medskip

\noindent An acyclic circuit also uses {\em bit-position transpositions}
$\tau_{i,j}$ for $i < j\,$ (where $u,v,w \in \{0,1\}^*$, $ \ i = |u| + 1$, 
$\,j = |uv| + 2$):

\medskip

 \ \ \ $\tau_{i,j}: \ \ u\,x_i\,v\,x_j \,w \ \longmapsto$
    $ \ u\,x_j\,v\,x_i \,w$.

\medskip

\noindent The transpositions are not viewed as gates (in actual circuits 
they are usually implemented as wire-permutations or wire-crossings). 

\medskip

Since the gates can be applied at any input positions, there are
infinitely many gates of each kind (viewed as input-output functions). 

Let $ \ \kappa_{1,n}:\, x_1 x_2 x_3 \,\ldots\, x_n$ $\,\longmapsto\,$
$x_2 x_3 \,\ldots\, x_n x_1\,$ be the {\em cyclic permutation} of the first 
$n$ variables; then we have 
$ \ \kappa_{1,n}(.) \,=\,$
$\tau_{n-1,n}\,\tau_{n-2,n-1} \ \ldots \ \tau_{2,3}\,\tau_{1,2}(.)$.
 \ It is well known that the set of finite permutations of $\omega$ (the 
natural integers) is generated by the set of transpositions 
$\{\tau_{j,j+1} : j \in \omega\}$; it is also generated by the set of 
transpositions $\{\tau_{1,j} :$ $j \in \omega\}$.

Thanks to $\tau = \{\tau_{i,j} : i,j \ge 1\}$ we only need the gates
 \ ${\sf not}_1$ (= {\sf not}), \ ${\sf and}_{1,2}$ (= {\sf and}),
 \ ${\sf or}_{1,2}$ (= {\sf or}), and ${\sf fork}_1$ (= {\sf fork}). 
Indeed,

\medskip

 \ \ \  \ \ \ ${\sf not}_j \,=\,$ 
$\tau_{1,j} \circ {\sf not}_1 \circ \tau_{1,j}(.)$,

\medskip

 \ \ \  \ \ \ ${\sf and}_{j,j+1} \,=\,$  
$\kappa_{1,j} \circ {\sf and}_{1,2} \circ \tau_{2,j+1} \circ \tau_{1,j}(.)$,

\medskip

 \ \ \  \ \ \ ${\sf fork}_j \,=\,$  
$(\kappa_{1,j+1})^2 \circ {\sf fork}_1 \circ \tau_{1,j}(.)$.

\bigskip

\noindent {\bf Extension to a right-ideal morphism}

\smallskip

For every circuit (and gate) we extend the input-output function $f:$ 
$\{0,1\}^m \to \{0,1\}^n\,$ to a right-ideal morphism by defining 

\medskip

 \ \ \  \ \ \ $f(xu) = f(x) \,u\,$, \ \ for all $x \in \{0,1\}^m$,
   $ \ u \in \{0,1\}^*$.

\medskip

\noindent Thus $\,{\rm domC}(f) = \{0,1\}^m$, 
$ \ {\rm Dom}(f) = \{0,1\}^m\,\{0,1\}^*$, 
$ \ {\rm imC}(f) \subseteq \{0,1\}^n$, and $\,{\rm Im}(f) \subseteq$ 
$\{0,1\}^n\,\{0,1\}^*$; here, $m$ is the number of input bits of the circuit,
and $n$ is the number of output bits. 
The definition of the gates, given above, carries out this extension already 
with \ ${\rm domC}({\sf not}_j) = \{0,1\}^j = {\rm domC}({\sf fork}_j)$,
 \ ${\rm domC}({\sf and}_j) = \{0,1\}^{j+1} = {\rm domC}({\sf or}_j)$, 
 \ ${\rm imC}({\sf not}_j) = \{0,1\}^j$  $= {\rm imC}({\sf and}_j)$
$ = {\rm imC}({\sf or}_j)$, and 
$\,{\rm imC}({\sf fork}_j) = \{0,1\}^{j+1}$.

This extension makes $f$ an element of the pre-Thompson monoid 
${\sf tfl}{\cal RM}^{\sf fin}_2$. The extension to right-ideal morphisms 
enables us to always compose input-output functions of circuits, and the 
composite is never the empty function.

The input-output function of the circuit is obtained by composing a sequence 
of copies of the four gates {\sf not}, {\sf and}, {\sf or}, {\sf fork}, and 
elements of $\tau$.

Conversely, every total function $\,\{0,1\}^m \to \{0,1\}^n\,$ (for any
given $m, n \ge 1$) can be realized by an acyclic boolean circuit (see 
e.g.\ \cite{CloteKran, JESavage, Wegener}). 

\smallskip

Later we will also consider the extension of the input-output function $f$ of
any acyclic circuit to a function $\,\{0,1\}^{\omega} \to \{0,1\}^{\omega}$,
by defining $\,f(xu) = f(x)\,u \ $ for all $x \in {\rm domC}(f)$, $\,u \in$ 
$\{0,1\}^{\omega}$.  This makes $f$ an element of the Thompson monoid 
${\sf tlep}M_{2,1}$.

\smallskip

After extending the input-output function $f_C$ of an acyclic circuit $C$
to a right-ideal morphism of $\{0,1\}^*$, or to a function 
$\,\{0,1\}^{\omega}$ $\to$  $\{0,1\}^{\omega}$, we need to redefine the 
input- and output-length, since now
$C$ has an unbounded sequence of input variables, and an 
unbounded sequence of output variables. But $f_C$ only modifies a finite
prefix $\,x_1 \ldots x_m \in {\rm domC}(f_C)\,$ of the total input into a 
finite output prefix $\,y_1\ldots y_n \in f_C({\rm domC}(f_C))$. I.e.,
$\,f_C(x_1 \ldots x_m\,x_{m+1} \ldots x_{m+k} \,\ldots \,)$ $\,=\,$ 
$f_C(x_1 \ldots x_m) \ x_{m+1} \ldots x_{m+k} \,\ldots$ $ \ =\,$
$y_1\ldots y_n\,y_{n+1}  \ldots  y_{n+k} \,\ldots \ $, such that for all 
$k \ge 1:$ $\,y_{n+k} = x_{m+k}\,$ (i.e., $y_{n+k}$ and $x_{m+k}$ have the
same value in $\{0,1\}$). The connection from the input variable
$x_{m+k}$ to the output variable $y_{n+k}$ is just a wire (with no gate or
wire permutation operation); we call the link $(x_{m+k}, y_{n+k})$ an 
{\em identity wire}.

\begin{defn} \label{DEF_IOlenCircuits} {\bf (input- and output-length 
of a boolean circuit).}

\smallskip

\noindent Let $\,(x_i:$ $i \in \omega)\,$ be the input variables of an 
acyclic circuit $C$ with input-output function extended to a right-ideal 
morphism of $\{0,1\}^*$, and let $\,(y_j:$ $j \in \omega)\,$ be the output 
variables.  
The {\em input-length} $\,\ell_{\rm in}(C)\,$ and the {\em output-length} 
$\,\ell_{\rm out}(C)\,$ of $C$ are defined by

\medskip

$\ell_{\rm in}(C)$  $\,=\,$
$\max\{\,i > 0 :$ {\rm a gate or a wire permutation in $C$ is connected to the
input variable} $x_i \}$.

\medskip

$\ell_{\rm out}(C)$  $\,=\,$
$\max\{\,j > 0 :$ {\rm a gate or a wire permutation in $C$ uses the output 
variable} $y_j \}$.
\end{defn}

\begin{defn} \label{DEFsizedepthCircuits} {\bf (size and depth of a
circuit).}

\smallskip

\noindent $\bullet$ The {\em size} of an acyclic circuit $C$, denoted by
$\,|C|$, is the maximum of $\,\ell_{\rm in}(C)$, $\,\ell_{\rm out}(C)$, and 
the number of gates of $C\,$ (the latter is the number $\,|V_i|$ of internal 
vertices of the {\sc dag} of $C$, i.e., the non-input and non-output 
vertices). I.e., $\,|C| \,=\,$ $\max\{\ell_{\rm in}(C)$, $\ell_{\rm out}(C)$, 
$|V_i|\}$.

Wire permutations are not explicitly counted in the size, but they
influence $\ell_{\rm in}(C)$ and $\ell_{\rm out}(C)$.

\smallskip

\noindent $\bullet$ The {\em depth} of an acyclic circuit $C$, denoted by
${\sf depth}(C)$, is the depth of the underlying {\sc dag}, i.e., the maximum
length of all directed paths.
\end{defn}

\noindent
For a function in ${\sf tfl}{\cal RM}^{\sf fin}$ (or ${\cal RM}^{\sf fin}$,
etc.), given by a sequence of generators, the input- and output-length can 
also be defined. 

Note that for the composition of functions we have:
$\,u_k \circ\,\ldots\, \circ u_1(x)\,$ is defined iff
for all $\,i \in [1,k]$,
$\,u_i \circ\,\ldots\,\circ u_1(x)\,$ is defined.

\begin{defn} \label{DEF_IOlengthRMfin} {\bf (input- and output-length
in ${\sf tfl}{\cal RM}^{\sf fin}$, ${\sf tot}{\cal RM}^{\sf fin}$,
${\sf pfl}{\cal RM}^{\sf fin}$, ${\cal RM}^{\sf fin}$).}

\smallskip

\noindent {\small \rm (1)} Let $\,\Gamma$ be a finite set such that 
$\,\Gamma \cup \tau$ generates ${\sf tfl}{\cal RM}^{\sf fin}$. 
 \ For any $\,u = u_k\,\ldots\,u_1 \in$  $(\Gamma \cup \tau)^*$ with 
$u_i \in \Gamma \cup \tau$ (for $i \in [1,k]$), we define:

\medskip

$\ell_{\rm in}(u)$ $ \ = \ $
$\min \{\,m \in {\mathbb N}\,:\,$
$(\forall x \in \{0,1\}^m)[$
$u_k \circ\,\ldots\, \circ u_1(x)\,$ {\rm is defined}$\,]\,\}$;

\medskip

$\ell_{\rm out}(u)\,$ is the unique
$n \in {\mathbb N}$ such that
$ \ u_k \circ\,\ldots\, \circ u_1(\{0,1\}^m) \,\subseteq\, \{0,1\}^n$, 
 \ where $m = \ell_{\rm in}(u)$.
 
\medskip 

\noindent The same definition of $\,\ell_{\rm in}(u)\,$  (but not of 
$\,\ell_{\rm out}(u)$) can be applied for $\,{\sf tot}{\cal RM}^{\sf fin}$.

\medskip

\noindent {\small \rm (2)} Let $\,\Gamma$ be a finite set such that 
$\,\Gamma \cup \tau$ generates ${\sf pfl}{\cal RM}^{\sf fin}$. 
 \ For any $\,u = u_k\,\ldots\,u_1 \in$  $(\Gamma \cup \tau)^*$ with 
$u_i \in \Gamma \cup \tau$ (for $1 \le i \le k$), we define:

\medskip

$\ell_{\rm in}(u)$ $ \ = \ $
$\min \{\,m \in {\mathbb N}\,:\,$

\hspace{2,5cm} $(\forall x \in \{0,1\}^m)\,[$
$x\,\{0,1\}^{\omega} \,\subseteq\, {\rm Dom}(E_u)\,\{0,1\}^{\omega}$
$ \ \Rightarrow \ $
$u_k \circ\,\ldots\, \circ u_1(x)\,$ {\rm is defined}$\,]\,\}$;

\medskip

$\ell_{\rm out}(u)\,$ is the unique
$n \in {\mathbb N}$ such that
$ \ u_k \circ\,\ldots\, \circ u_1(\{0,1\}^m) \,\subseteq\, \{0,1\}^n$,

\hspace{1,25cm} where $m = \ell_{\rm in}(u)$.

\medskip

\noindent The same definition of $\,\ell_{\rm in}(u)\,$ (but not of
$\,\ell_{\rm out}(u)$) can be applied for ${\cal RM}^{\sf fin}$.

\medskip

\noindent  {\small \rm (3)} For the Thompson monoids ${\sf tot}M_{2,1}$,
${\sf tfl}M_{2,1}$, $M_{2,1}$, ${\sf pfl}M_{2,1}$, acting on 
$\{0,1\}^{\omega}$, we define the input-length (and an output-length in the
{\sf plep} case) by using words in $(\Gamma \cup \tau)^*$, where $\Gamma$ is
finite such that $\Gamma \cup \tau$ generates ${\sf tot}{\cal RM}^{\sf fin}$,
${\sf tfl}{\cal RM}^{\sf fin}$, ${\cal RM}^{\sf fin}$, 
${\sf pfl}{\cal RM}^{\sf fin}$, respectively.  I.e., we use the generating 
sets of the corresponding pre-Thompson monoids. 
\end{defn}
The following is straightforward.

\begin{lem} \label{EquivDefInputLen} {\bf (equivalent definition of 
input-length).}

\noindent Let $\Gamma$ be a finite set such that $\,\Gamma \cup \tau$ 
generates ${\cal RM}^{\sf fin}$ (or ${\sf tot}{\cal RM}^{\sf fin}$,
${\sf tfl}{\cal RM}^{\sf fin}$, ${\sf pfl}{\cal RM}^{\sf fin}$), and let 
$\,u = u_k \,\ldots\,u_1 \in $ $(\Gamma \cup \tau)^*$ with $u_i \in$ 
$\Gamma \cup \tau$ for $i \in [1,k]$.
 \ Then

\smallskip
 
 \ \ \  \ \ \ $\ell_{\rm in}(u) \ = \ \min \big\{ m \in {\mathbb N}\,:\,$
$\{0,1\}^m \,\cap\, {\rm Dom}(u_k \circ \ldots \circ u_1(.))$
$ \ \equiv_{\rm fin} \ {\rm domC}(u_k \circ \ldots \circ u_1(.))\,\big\}$.

\noindent $\Box$
\end{lem}

\begin{defn}  \label{word_length_size}  {\bf (lengths and sizes of words).}

\smallskip

\noindent Let $\Gamma$ be a finite set such that 
$\Gamma \cap \tau = \0$, and let $w \in (\Gamma \cup \tau)^*$.  Then

\smallskip

\noindent $\bullet$ \ \ ${\rm maxindex}_{\tau}(w)$ $\,=\,$
$\max\{i :\, \tau_{i-1,i} \ {\rm occurs \ in} \ w \}$; \ and

${\rm maxindex}_{\tau}(w) = 0$ \ if no element of $\tau$
occurs in $w$.

\smallskip

\noindent $\bullet$ \ \ For a finite set $S \subseteq \{0,1\}^*$,
 \ ${\rm maxlen}(S) = \max\{|s| : s \in S\}$.

For a finite set ${\cal F} \subseteq {\cal RM}^{\sf fin}$,
$ \ {\rm maxlen}({\cal F}) \,=\, \max \big\{\,|s| :\, s \in $
$\,\bigcup_{f \in {\cal F}} \big({\rm domC}(f) \cup f({\rm domC}(f))\big)
\,\big\}$.

\smallskip

\noindent $\bullet$ \ \ $|w|_{\Gamma}\,$ is the number of occurrences of 
letters of $\,\Gamma$ in $w$;

\smallskip

\noindent $\bullet$ \ \ $|w|_{\Gamma \cup \tau}\,$ is the number of 
occurrences of elements in $\Gamma \cup \tau$ in $w$; 

this is the length of $w$ where every element of 
$\Gamma \cup \tau$ is assigned length $1$;

\smallskip

\noindent $\bullet$ \ \ $|w| \,=\, $
$|w|_{\Gamma \cup \tau} \,+\, {\rm maxindex}_{\tau}(w)\,$; 
 \ this is called the \textit{\textbf{total word-length}} of $w$.
\end{defn}
The total word-length of $w$ could have been defined to be
$ \ \max\{|w|_{\Gamma \cup \tau}, \,{\rm maxindex}_{\tau}(w)\}$. 
But this does not make much difference (since for positive real numbers
$a, b$ we have $\, \max\{a,b\} < a+b < 2\,\max\{a,b\}\,$).

\begin{defn} \label{DEFevalfunction} {\bf (evaluation function).}
 \ For a generating set $\Gamma \cup \tau$ of one of the
above monoids (e.g., ${\cal RM}^{\sf fin}$), the {\em generator evaluation} 
function $\,E: w \in (\Gamma \cup \tau)^*$  $\,\longmapsto\,$ 
$E_w \in {\cal RM}^{\sf fin}\,$ is defined for any $x \in A^*$ by 

\smallskip

\hspace{2,5cm} $E_w(x) = w_k \circ \,\ldots\, \circ w_1(x)$, 

\smallskip

\noindent where $w = w_k\,\ldots\,w_1$ with $w_i \in \Gamma \cup \tau$ for 
$i \in [1,k]$. For a given $w \in (\Gamma \cup \tau)^*$, the right-ideal
morphism $\,E_w$ is called the {\em input-output function} of $w$.

For partial circuits, the {\em circuit evaluation} function $\,E: C$ 
$\,\longmapsto\,$  $E_C\,$ is such that for any partial circuit $C$, $E_C$ 
is the input-output function of $C$ (extended to a right-ideal morphism of
$\{0,1\}^*$). 
\end{defn}

\smallskip

\noindent {\large \bf Encoding of a circuit by a bitstring:}

\smallskip

\noindent
As a preliminary step, we encode the alphabet $\{0,1, \#\}$ over $\{0,1\}$.
This can be done, e.g., by the injective function $\,{\sf code}_0(.)$:
$\{0,1, \#\}^* \to \{0,1\}^*$ defined by $\,{\sf code}_0(0) = 00$,
$\,{\sf code}_0(1) = 01$, $\,{\sf code}_0(\#) = 1$, and
$\,{\sf code}_0(xy) = {\sf code}_0(x) \ {\sf code}_0(y)\,$ for all $x, y \in$
$\{0,1, \#\}^*$. The set $\{00, 01, 1\}$ is a maximal prefix code.

Let $(V,E)$ be the set of vertices, respectively edges, of an acyclic circuit 
(boolean or partial). We assume that every vertex $v$ is a bitstring of 
length $\lceil \log_2 |V| \rceil$, and different vertices are different
bitstrings. Moreover, the bitstrings used here represent the numbers 0
through $|V| - 1$ in binary representation of length 
$\lceil \log_2 |V| \rceil\,$ (with leading 0s when needed).

For every input vertex $v$ we consider the string
$\, \# v \# b \# \ell \# r \#\,$ over the alphabet $\{0,1, \#\}$. Here $b$
is the gate-type of $v$, namely 000, 001, 010, 011, 100, 101, for 
respectively Input, Output, {\sf and}, {\sf or}, {\sf not}, {\sf fork}; in
a partial circuit, we also allow the partial gate $\zeta_1$ (see Subsection
3.2), encoded by $110$.  
The bitstrings $\ell, r$ are the {\em predecessor} vertices (a.k.a.\ the 
parents) of $v$; they provide the inputs for $v$.  If $v$ is a source 
(a.k.a.\ an input vertex) then $\ell, r$ are each a copy of $v$; if $v$ has 
in-degree 1 then $r$ is a copy of $\ell$.

In addition, we require that the vertex labels $v \in V$ are chosen so that
$\,\# v \# b \# \ell \# r \#\,$ satisfies $v <_{\rm lex} \ell$ and
$v <_{\rm lex} r\,$  (unless $v$ is a source, and then  $v = \ell = r$).
So, the names of the vertices and their lexicographic order are compatible
with the {\em precedence} order (dependence order, a.k.a.\ topological       sort).

To represent the circuit by a string we concatenate all the strings that
describe the vertices; let $s \in \{0,1,\#\}^*$ be the resulting string.
We do this in {\em sorted order of the vertices} $v \in \{0,1\}^*$,          
according to the dictionary order in $\{0,1\}^*$.
This, together with the fact that $v <_{\rm lex} \ell$ and $v <_{\rm lex} r$,
implies that in an encoding of an acyclic circuit, the vertices are 
{\em topologically sorted} (i.e., the partial order defined by the acyclic 
digraph, as a Hasse diagram, is compatible with the lexicographic order of 
the vertex strings).

Finally, we consider ${\sf code}_0(s\,\#\#)$ to obtain the encoding of a 
circuit $C$, denoted by ${\sf Code}_<(C)$.  We call ${\sf Code}_<(C)$ the 
{\em precedence encoding} of the acyclic circuit $C$.
The set of all strings ${\sf Code}_<(C)$ that represent acyclic circuits is 
called the {\em precedence connection language}.
In our formulation, the precedence connection language is an infinite prefix
code.

\medskip

\noindent {\bf Decoding:} From a precedence encoding ${\sf Code}_<(C)$ of a
boolean or partial circuit $C$, a circuit isomorphic to $C$ can be 
reconstructed, up to circuit isomorphism; $C$ itself is not uniquely 
determined by ${\sf Code}_<(C)\,$ (see the remarks about isomorphic circuits
at the beginning of this subsection).

\bigskip

\noindent Let us show that the precedence connection language 
belongs to {\sc DSpace}$(\log)$.

\smallskip

\noindent 1. A finite-state machine can decode a string from $\{0,1\}^*$ to
a string in $\{0,1, \#\}^*$; i.e., apply the function ${\sf code}_0^{-1}(.)$,
and reject if this is not decodable.
A finite-state machine can also check whether the decoded string belongs to
 \ $\big(\# \,\{0,1\}^+ \,\#\, \{000, 001, 010, 011, 100, 101, 110\}$
$\#\, \{0,1\}^+ \,\#\, \{0,1\}^+ \,\# \big)^+\,$.

\smallskip

\noindent 2. We can check in log-space whether in a string

\smallskip

 \ \ \  $\# v_1 \# b_1 \# \ell_1 \# r_1 \#$
$\# v_2 \# b_2 \# \ell_2 \# r_2 \#$  $ \ \ldots \ $
$\# v_i \# b_i \# \ell_i \# r_i \#$  $ \ \ldots \ $
$\# v_{|V|} \# b_{|V|} \# \ell_{|V|} \# r_{|V|} \# \,$,

\smallskip

\noindent all the vertex strings $v_i$, $\ell_i$, and $r_i$ (over all $i$)
have the same length, whether $v_i <_{\rm lex} \ell_i$ and
$ v_i <_{\rm lex} r_i$ (unless $v_i = \ell_i = r_i$),
and whether the vertices $v_i$ occur in strict sorted lexicographic order
(with increments of 1 at each step, starting with $v_1 = $
$0^{\lceil \log_2 |V| \rceil} 1$ and ending at the binary
representation of $|V| - 1$).

\smallskip

\noindent 3. For each vertex $v$ we check its in-degree.  If $v$ is an 
{\sf and}, {\sf or}, or a $\zeta_1$ gate, we check whether $\ell$, $r$,
and $v$, are all different.
If $v$ is a {\sf not} or a {\sf fork} gate or and output vertex, we check
whether $\ell = r <_{\rm lex} v$.
If $v$ is an input vertex, we check whether $\ell = r = v$. This is easily
done in log-space.

\smallskip

\noindent 4. For each vertex $v_i$ we check its out-degree:
If $v_i$ is an output vertex, $v_i$ must not be equal to $\ell_j$ or
$r_j$ of any other vertex $v_j$.
If $v_i$ has out-degree 1 (i.e., it is of type Input, {\sf and}, {\sf or}, 
{\sf not} or $\zeta_1$) then $v_i$ must be equal to exactly one parent 
vertex $\ell_j$ or $r_j$ of another vertex $v_j$, with $i < j\,$ (because 
of topological sorting).
If $v_i$ is of type {\sf fork}, then there must be two vertices $v_j$ and
$v_k$, with $i,j,k$ all different and $i < j, k\,$ (by topological sorting),
such that $v_i$ is equal to $\ell_j$ or $r_j$, and to $\ell_k$ or $r_k$.

\bigskip

\noindent {\bf Computing properties of a circuit:} We would like to compute
the input-length, output-length, size, and depth of a boolean or partial
circuit $C$ in log-space, from its precedence encoding ${\sf Code}_<(C)$.

Given ${\sf Code}_<(C)$, a finite-state automaton can output the number of 
input vertices in unary notation (by recognizing the gate-type Input, i.e., 
$000$); similarly, the number of output vertices, or the total number of 
vertices can be computed in unary.
To compute the depth $d(C)$ of $C$ in binary notation, we give an algorithm
with space complexity $O(d(C) + \log |V|)$. The algorithm uses depth-first
search, as follows:

We use every source vertex, one after another.
We remember the currently chosen source vertex $s_i$ (in space $\log |V|$),
and do depth-first search from $s_i$, remembering the sequence of left
or right steps in the {\sc dag}, remembering left as 0, and right as 1;
the search goes left when possible; this yields a bitstring of length
$\le d(C)$; when this depth-first phase ends, the length of the search
string is computed in binary, and remembered (in space $\le \log d(C)$).
Now the search backtracks, using the search string in reverse, reading 0s
until a 1 is found, or the start vertex is found.
Now the search follows the 1 found, and then goes left again.
The new depth eventually found is compared with the stored one, and the
larger one (in lexicographic order) is kept. If the vertex is found in the 
backtrack, and a 0 was read just before the start vertex, this 0 is replaced 
by 1, and the search resumes. If a 1 was read just before the start vertex, 
the search resumes with the next source vertex (unless there is no next 
source vertex, in which case the search is over).

In a similar way, $\ell_{\rm in}(C)$ and $\ell_{\rm out}(C)\,$ (Def.\ 
\ref{DEF_IOlenCircuits}) can be computed in log-space, by doing reverse 
depth-first search (in the reverse {\sc dag}) or depth-first search.

In summary we proved the following.

\begin{lem} \label{LEMcodeDAG}
 \ There is an injective function ${\sf Code}_<(.)$ (the {\em precedence
encoding}) that maps every acyclic boolean or partial circuit to a bitstring
with the following properties:
The image set of ${\sf Code}_<(.)$ (i.e., the set of all encodings of acyclic
circuits, also called the {\em precedence connection language}) is a prefix
code that belongs to {\sc DSpace}$(\log)$.
Moreover, there is an algorithm that on input ${\sf Code}_<(C)$, for any
acyclic circuit $C$, computes the depth ${\sf depth}(C)$  with space
complexity $O({\sf depth}(C) + \log |V|)\,$ (where $|V|$ is the number of
vertices in the circuit). And there are log-space algorithms that compute
$\ell_{\rm in}(C)$, $\ell_{\rm out}(C)$, and $|C|$. 
  \ \ \   \ \ \ $\Box$
\end{lem}

\smallskip

In this paper, whenever a circuit $C$ is an input or an output of an 
algorithm, we actually use the encoding of $C$ by a bitstring as the 
input or output of the algorithm.

\medskip

Theorem \ref{circ2lepM} describes a close relation between acyclic circuits
and words over a generating set $\Gamma_{\sf\!tfl} \cup \tau$ of
${\sf tfl}{\cal RM}^{\sf fin}_2$.
This relation was already studied in \cite[{\small Prop.\ 2.4}]{Bi1wPerm},
but Theorem \ref{circ2lepM} is more detailed in its parts (2) and (3).

\begin{thm} \label{circ2lepM}  {\bf (boolean circuits versus words in
$\,{\sf tfl}{\cal RM}^{\sf fin}_2$).}

\smallskip

\noindent Let $\Gamma_{\sf\!tfl}$ be a finite set such that
$\Gamma_{\sf\!tfl} \cup \tau$ generates $\,{\sf tfl}{\cal RM}^{\sf fin}$, 
and $\,\Gamma_{\sf\!tfl} \cap \tau = \0$.
For $w \in (\Gamma_{\sf\!tfl} \cup \tau)^*$, let $\,|w| \,=\, $
$|w|_{\Gamma \cup \tau} \,+\, {\rm maxindex}_{\tau}(w)\,$.

\medskip

\noindent {\small \bf (1)} {\rm (From circuit to word):}
 \ Let $C$ be any acyclic boolean circuit, with gates of type {\sf and, 
or, not, fork}. Then a word $\,{\sf W}(C) \in$ 
$(\Gamma_{\sf\!tfl} \cup \tau)^*\,$ can be constructed such that $C$ and 
${\sf W}(C)$ have the same input-output function. Moreover,

\smallskip

 \ \ \  \ \ \ $|{\sf W}(C)| = \Theta(|C|)\,$;

\smallskip

 \ \ \  \ \ \ ${\rm maxindex}_{\tau}({\sf W}(C)) \,\le\, |C|\,$;

\smallskip

 \ \ \  \ \ \ $\ell_{\rm in}({\sf W}(C)) = \ell_{\rm in}(C)$, \ \ and 
 \ \ $\ell_{\rm out}({\sf W}(C)) = \ell_{\rm out}(C)$.

\smallskip

\noindent And $\,{\sf W}(C)$ can be computed from $C$ in log-space
in terms of $|C|$.

\medskip

\noindent {\small \bf (2)} {\rm (From word to circuit):}
 \ From any word $u \in (\Gamma_{\sf\!tfl} \cup \tau)^*$, an acyclic
boolean circuit ${\sf C}(u)$, over the gate set $\{{\sf and, or, not,
fork}\}$, can be constructed such that $\,{\sf C}(u)$ and $u$ have the
same input-output function. Moreover, 

\smallskip

 \ \ \  \ \ \ $|{\sf C}(u)| = \Theta(|u|)$;

\smallskip

 \ \ \  \ \ \ $\ell_{\rm in}({\sf C}(u)) = \ell_{\rm in}(u)$ \ \ and
 \ \ $\ell_{\rm out}({\sf C}(u)) = \ell_{\rm out}(u)$.

\smallskip

\noindent And $\,C(w))$ can be computed from $u$ in log-space in terms of 
$|u|$.

\medskip

\noindent {\small \bf (3)}  If boolean circuits use $\Gamma_{\sf\!tfl}$ as
the set of gates, then $\,{\sf W}(.)$ and $\,{\sf C}(.)$ can be chosen so 
that 

\smallskip

 \ \ \  \ \ \ ${\sf C}({\sf W}(C))$ $ \ \widehat{=} \ $  $C$ 

\smallskip

\noindent (i.e., the boolean circuits $C$ and $\,{\sf C W}(C)$ are 
isomorphic).
\end{thm}
{\sc Proof.} (1) Let us assume at first that \{{\sf and}, {\sf or},
{\sf not}, {\sf fork}\} $\subseteq$ $\Gamma_{\sf\!tfl}$. We adapt the proof
of \cite[{\small Prop.\ 2.4}]{Bi1wPerm}, but some changes are needed because
here we use ${\sf tfl}{\cal RM}^{\sf fin}_2$ instead of ${\sf plep}M_{2,1}$.

We construct ${\sf W}(C)$ by simulating the gates in the left-to-right order
in each depth-layer of $C$. The wires that leave a layer are numbered from
left to right, and similarly the wires that enter a layer are numbered from
left to right.  Each gate of the boolean circuit $C$ is simulated by a word 
over $\Gamma_{\sf\!tfl} \cup \tau$. The reasoning is the same for every
gate, so let us just focus on a gate $g = {\sf or}_{i_1,i_2,j}$ with 
in-wire numbers $i_1$ and $i_2$, and out-wire number $j$.  Here we assume 
that all gates that precede $g$ (i.e., to the left of $g$ in the layer, or at 
the right end of the previous layer if $g$ is at the left end) have already 
been applied, in order to produce a suffix of ${\sf W}(C)$; and we assume 
that the input sequence, after all the gates that precede $g$ have been 
applied, is $(x_1, x_2, \ldots, x_{i_1}, \ldots, x_{i_2}, $ $\ldots$, 
$x_{m_1})$.
By using bit-transpositions we have: \ ${\sf or}_{i_1,i_2}(.) = $
${\sf or} \circ \tau_{2,i_2} \circ \tau_{1,i_1}(.)$.
Then the complete output so far is $(y_1, x_3, \ldots, x_{i_1-1}, x_1,$
$\ldots, x_2, x_{i_2 + 1}, \ldots, x_{m_1})$.
The output of ${\sf or}_{i_1,i_2,j}(.)$ is expected to go to position $j$,
whereas ${\sf or} \circ \tau_{2,i_2} \circ \tau_{1,i_1}(.)$ sends its
output bit $y_1$ to position 1. But instead of permuting all the wires
$1$ through $j$ in order to get the output of
${\sf or} \circ \tau_{2,i_2} \circ \tau_{1,i_1}(.)$ to wire $j$, we just
leave $y_1$ on wire 1 for now.  The simulation of the next gate of $c$ will
then use the appropriate bit-transpositions $\tau_{2,\ell} \circ \tau_{1,k}$
for moving the correct input wires to the next gate.  Thus, each gate of $C$
is simulated by one generator in $\Gamma_{\sf\!tfl}$ and two wire-swaps in 
$\tau$.

Finally, when all the output variables of the boolean circuit have been 
computed, a permutation of the $n$ output wires is used in order to send the 
outputs to the correct final output wires. This will also take care of 
wire-swaps of wires that are not connected to any gate.
Any permutation of $n$ elements can be realized with $< n$ ($\le |C|$)
bit-transpositions.

\smallskip

We have $\ell_{\rm in}({\sf W}(C)) = \ell_{\rm in}(C)$ and  
$\ell_{\rm out}({\sf W}(C)) = \ell_{\rm out}(C)$ because $C$ and
${\sf W}(C)$ use the same input wires, and the same output wires.

\smallskip

The word ${\sf W}(C)$ can be computed by a log-space Turing machine, since
the only thing to be remembered in the computation is a few numbers that are
$\le |C|\,$ (the lengths of the intermediary input and output strings, and up
to three positions in the intermediary input and output strings).

If \{{\sf and}, {\sf or}, {\sf not}, {\sf fork}\} $\not\subseteq$ $\Gamma$,
any $g \in$ \{{\sf and}, {\sf or}, {\sf not}, {\sf fork}\} $\minus$ $\Gamma$
can be replaced by a word $w_g \in (\Gamma \cup \tau)^*$. The gate $g$ and
the word $w_g$ represent the same function; therefore the reasoning with
\{{\sf and}, {\sf or}, {\sf not}, {\sf fork}\} $\subseteq$ $\Gamma$
still works for the words $w_g$, except that
$|{\sf W}(C)|_{\Gamma} \le s \,|C|$, where
$s = \max\{ |w_g|_{\Gamma} : g \in \Gamma_{\sf\!tfl}\}$.

\medskip

\noindent
(2) Since $\,$\{{\sf and}, {\sf or}, {\sf not}, {\sf fork}\} $\cup$ $\tau\,$
is a generating set of ${\sf tfl}{\cal RM}^{\sf fin}_2$ we consider first the 
case where $\Gamma_{\sf\!tfl}$ $=$ \{{\sf and}, {\sf or}, {\sf not}, 
{\sf fork}\}. Then the construction of a boolean circuit ${\sf C}(u)$ from 
$u \in$ $(\Gamma_{\sf\!tfl} \cup \tau)^*$ will be an inversion of the 
construction in the proof of (1).

\smallskip

\noindent Notation: \ For $g \in \Gamma_{\sf\!tfl}\,$ (finite) we let 
 
 \ \ \ $\ell_{\rm in}(g)$ $=$ ${\rm maxlen}\big({\rm domC}(g)\big)$, \ and
 \ $\ell_{\rm out}(g)$ $=$ ${\rm maxlen}\big(g({\rm domC}(g))\big)$
  $= {\rm maxlen}\big({\rm imC}(g)\big)$.

\smallskip

\noindent Every $u \in (\Gamma_{\sf\!tfl} \cup \tau)^*$ has the form

\smallskip

 \ \ \  \ \ \ $u$ $\,=\,$
$t_r g_r \ \ldots \ t_{j+1} g_{j+1} t_j g_j \ \ldots \ t_2 g_2 t_1 g_1 t_0$,

\smallskip

\noindent where $r \ge 0$, $\,g_j$ $\in$ $\Gamma_{\sf\!tfl}$ (for
$1 \le j \le r$), and $t_j \in \tau^+$ (for $0 \le j \le r$); i.e., every
$t_j$ is a non-empty sequence of bit-transpositions, representing a finite
permutation of ${\mathbb N}_{>0}$.  The boolean circuit ${\sf C}(u)$ is 
constructed in two steps (described in detail below): 
 \ (a) For every $t_j$ we define a bipartite {\sc dag} that represents 
the wire permutation $t_j$, and then we connect $g_j$ to the {\sc dag}s of 
$t_{j-1}$ and $t_j$.
 \ (b) We eliminate all intermediary vertices and just keep the gates $g_j$,
as well as the input variables and the output variables.
Wire-swaps are automatically kept (they are implemented by the adjacency
relation of the {\sc dag} of the circuit; they do not correspond to
vertices).
In more detail:

\smallskip

\noindent (a) \ For every $j$ with $0 \le j \le r$, let $\,\ell_j$  $=$
$\max\{i \in {\mathbb N}_{>0} : t_j(i) \ne i\}$; note that $t_j(i) \ne i$ is
equivalent to $t_j^{-1}(i) \ne i$.
So for all $i > \ell_j$: $ \ t_j(i) = i$. Hence $t_j$, which is a finitary
permutation of ${\mathbb N}_{>0}$, can be viewed as a permutation of
$\,\{i : 1 \le i \le \ell_j\}$.

The bipartite {\sc dag} of $t_j$ has input vertex set
$\,\{x_i^{(j, {\rm in})} : i \in {\mathbb N}_{>0}\}$, output vertex set
$\,\{x_i^{(j, {\rm out})} : i \in {\mathbb N}_{>0}\}$, and edge set
$\,\{(x_i^{(j, {\rm in})}, x_{t_j(i)}^{(j, {\rm out})}) :$
$i \in {\mathbb N}_{>0}\}$.

For every $j$ with $1 \le j \le r$, we connect $g_j$ to the {\sc dag} of
$t_{j-1}$ and to the {\sc dag} of $t_j$ as follows:
We connect $x_i^{(j-1, {\rm out})}$ to the $i$th input variable of $g_j$,
for $1 \le i \le \ell_{\rm in}(g_j)$;
we connect the $i$th output variable of $g_j$ to $x_i^{(j, {\rm in})}$, for
$1 \le i \le \ell_{\rm out}(g_j)$;
moreover, we connect $x_i^{(j-1, {\rm out})}$ to
$x_{i - \ell_{\rm in}(g_j) + \ell_{\rm out}(g_j)}^{(j, {\rm in})}$,
for all $i > \ell_{\rm in}(g_j)$.
This connects the {\sc dag}s of all the wire permutations $t_j$
($0 \le j \le r$) and all the gates $g_j$ ($1 \le j \le r$) into one
{\sc dag}; the vertex set of this {\sc dag} is

\smallskip

 \ \ \
$\{x_i^{(j, {\rm in})} : 0 \le j \le r, \ i \in {\mathbb N}_{>0}\}$
$\,\cup\,$
$\{g_j : 1 \le j \le r\}$
$\,\cup\,$
$\{x_i^{(j, {\rm out})} : 1 \le j \le r, \ i \in {\mathbb N}_{>0}\}$.

\smallskip

\noindent The vertices $x_i^{(j, {\rm in})}$ (for $1 \le j \le r$), and
$x_i^{(j, {\rm out})}$ (for $0 \le j < r$), are called ``intermediary
vertices''.

\smallskip

\noindent The acyclic digraph constructed in (a) can be computed in 
log-space.  Indeed, the only things to remember during the computation are 
numbers $\le |u|$.

\medskip

\noindent (b) \ For all $j, k \in \{1, \ldots, r\}$ with $j \ne k$, we add
an edge $(g_j, g_k)$, from the $\alpha$th output of $g_j$
($1 \le \alpha \le \ell_{\rm out}(g_j)$) to the $\beta$th
input of $g_k$ ($1 \le \beta \le \ell_{\rm out}(g_k)$), iff
in the {\sc dag} constructed in (a) there exists a path from $g_j$ to $g_k$,
from the $\alpha$th output of $g_j$ to the $\beta$th input of $g_k$, such
that this path does not contain any other gate $g \in$ $\Gamma_{\sf\!tfl}$.
And we add an edge $(x_i^{(0, {\rm in})}, g_j)$, from the input variable
$x_i^{(0, {\rm in})}$ to the $\alpha$th input of $g_j$, iff in the {\sc dag}
constructed in (a) there exists a path from  $x_i^{(0, {\rm in})}$ to the
$\alpha$th input of $g_j$, such that this path does not contain any other
gate $g \in \Gamma_{\sf\!tfl}$.
And we add an edge $(g_j, x_i^{(r, {\rm out})})$, from the $\alpha$th output
of $g_j$ to the output variable $x_i^{(r, {\rm out})}$, iff in the {\sc dag}
constructed in (a) there exists a path from the $\alpha$th output  of $g_j$  to
the output variable $x_i^{(r, {\rm out})}$, such that this path     does not
contain any other gate $g \in \Gamma_{\sf\!tfl}$.
And we add an edge $(x_i^{(0, {\rm in})}, x_m^{(r, {\rm out})})$ iff
in the {\sc dag} constructed in (a) there exists a path from $x_i^{(0, {\rm
in})}$ to $x_m^{(r, {\rm out})}$, such that this path does not contain any
gate $g \in \Gamma_{\sf\!tfl}$.
Now we drop all intermediary vertices, and the edges incident to one or more
intermediary vertices.

The resulting {\sc dag} is a boolean circuit ${\sf C}(u)$ with input 
variables $x_i^{(0, {\rm in})}$ (for $i \in {\mathbb N}_{>0}$), output 
variables $x_i^{(r, {\rm out})}$ (for $i \in {\mathbb N}_{>0}$), and gate 
vertices labeled by $g_1, \ldots, g_r$. These input (or output) variables 
are the same as those of $u$.
Wire-swaps are automatically kept, via the adjacency relation in the
{\sc dag}.
The input-output function of the circuit ${\sf C}(u)$ is the same as the one
of the {\sc dag} constructed in (a), which is obtained by composing the 
letters of $u$; hence $u$ and ${\sf C}(u)$ have the same input-output 
function.

The number of gates of ${\sf C}(u)$ is $|u|_{\Gamma_{\sf\!tfl}}$.
 \ \ \  \ \ \ [End of (b).]

\smallskip

\noindent The acyclic digraph constructed in (b) can be computed in 
log-space.  Indeed, the only things to remember during the computation are 
numbers $\le |u|$.

The order of the gates in $u$ yields a precedence order for ${\sf C}(u)$
which is the same as the gate order in $u$.

\medskip

In case $\Gamma_{\sf\!tfl} \ne$ $\{{\sf and, or, not, fork}\}$, we use the 
above method to build an acyclic  boolean circuit over the gate set 
$\Gamma_{\sf\!tfl}$. After this we can replace each gate $\gamma$ in the 
finite set $\Gamma_{\sf\!tfl}$ by an acyclic boolean circuit
${\sf C}(\gamma)$ over $\{{\sf and, or, not, fork}\}$.  The inputs and 
outputs of ${\sf C}(\gamma)$ are the same as the ones of $\gamma$.
Since $\Gamma_{\sf\!tfl}$ is finite, this increases the size of ${\sf C}(u)$
by a constant factor.

\smallskip

Since steps (a) and (b) use log-space, the final result ${\sf C}(u)$ is 
constructed in log-space (since log-space is closed under composition).

\smallskip

We have $\ell_{\rm in}({\sf C}(u)) = \ell_{\rm in}(u)$ and 
$\ell_{\rm out}({\sf C}(u)) = \ell_{\rm out}(u)$ because $u$ and 
${\sf C}(u)$ use the same input wires and the same output wires.

\medskip

\noindent (3) \ Let us prove that ${\sf C W}(C) \ \widehat{=} \ C$.
By the construction in (1), ${\sf W}(C)$ contains the same gates as
$C$, in the same order (where in $C$ we order the gates by layers, starting
at the top, and from left to right in each layer). And ${\sf W}(C)$
describes the same connection between gates (and between input or output
variables and gates) as $C$.
By the construction in (2), ${\sf C W}(C)$ contains the same gates as
${\sf W}(C)$, in the same order. Hence, ${\sf C W}(C)$ contains the same
gates as $C$, in the same order. The connections in $C$ and in ${\sf C W}(C)$
are the same, since both ${\sf C}(.)$ and ${\sf W}(.)$ preserve the
connections between generators $g_j$. Hence, ${\sf W}(C)$ and $C$ are
isomorphic boolean circuits.

However, usually the wire permutations that implement the connections between
gates are different in ${\sf C W}(C)$ and $C$.
 \ \ \  \ \ \ $\Box$

\bigskip

\noindent {\bf Remark:} In Theorem \ref{circ2lepM}(3) we mainly care about
the fact that the boolean circuit $C$ is isomorphic to ${\sf C W}(C)$. 
Depending on the details of the encoding we could also have 
$\,{\sf WC}(u) = u$.

\bigskip

\noindent As a consequence of the Theorem we have:

\begin{pro} \label{PROPtflGen} {\bf (generators of 
$\,{\sf tfl}{\cal RM}^{\sf fin}_2$).}

\smallskip

The monoid $\,{\sf tfl}{\cal RM}^{\sf fin}_2$ is the set of input-output 
functions of acyclic boolean circuits, extended to right-ideal morphisms of 
$\{0,1\}^*$.

\smallskip

And $\,{\sf tfl}{\cal RM}^{\sf fin}_2$ has a generating set 
$\Gamma \cup \tau$, where $\Gamma$ is finite.  The finite set $\Gamma$ can 
be chosen to be \{{\sf not}, {\sf and}, {\sf or}, {\sf fork}\}.
 \ \ \  \ \ \  $\Box$
\end{pro}
We will prove later that ${\sf tfl}{\cal RM}^{\sf fin}_2$ is not finitely 
generated (Theorem \ref{PROcircNotFinGen}), but that ${\cal RM}^{\sf fin}$ 
is finitely generated (Theorem \ref{RMfinFinGen}) .

\bigskip

\noindent {\bf Remark about change of generators in
${\rm tfl}{\cal RM}^{\sf fin}_2$:} 

\smallskip

\noindent Let $\Gamma_1$ and $\Gamma_2$ be two finite sets such that 
$\Gamma_1 \cup \tau$ and $\Gamma_2 \cup \tau$ generate 
${\rm tfl}{\cal RM}^{\sf fin}_2$.
For every $\gamma \in \Gamma_2$ we can choose $k(\gamma) \in$
$(\Gamma_1 \cup \tau)^*$ be such that $\gamma$ and $k(\gamma)$ represent
the same element of ${\rm tfl}{\cal RM}^{\sf fin}_2$; and for every 
$\tau_{i,j} \in \tau$ we let $k(\tau_{i,j}) = \tau_{i,j}$. 
Since $\gamma \in {\rm tfl}{\cal RM}^{\sf fin}_2$ and $\Gamma_2 \cup \tau$ is
a generating set, $k(\gamma)$ exists.
For all $w = w_k \,\ldots\,w_1 \in$ $(\Gamma_1 \cup \tau)^*\,$ (with 
$w_i \in$ $\Gamma_1 \cup \tau$) we define $\,k(w) = k(w_k)\,\ldots\,k(w_1)$.

Then $w$ and $k(w)$ represent the same right-ideal morphism in
${\rm tfl}{\cal RM}^{\sf fin}_2$.

\subsection{Partial circuits}          

We saw already that a {\em partial acyclic circuit} (more shortly, a 
{\em partial circuit}) is defined in the same way as an acyclic boolean 
circuit, except that now elements from a finite set of additional 
{\em partial gates} in ${\sf pfl}{\cal RM}^{\sf fin}_2$ can be included 
in the circuit (see the beginning od Subsection 3.1).  

In a partial circuit, if one or more gates produce an undefined output for
a given input, then the whole output of the circuit is undefined for that 
input.

\medskip

Here only one new partial gate, called $\zeta_1$, is introduced. It turns 
out that $\,\{\zeta_1, {\sf not}, {\sf and}, {\sf or}, {\sf fork}\}\,$ is a
complete set of gates for partial circuits. The gate $\zeta_1$ is defined as
the right-ideal morphism of $\{0,1\}^*$ with 
$\,{\rm domC}(\zeta_1) = \{00, 01\}$, $ \ {\rm imC}(\zeta_1) = \{0,1\}$, and

\smallskip

 \ \ \  \ \ \ $\zeta_1(00) = 0$;

\smallskip

 \ \ \  \ \ \ $\zeta_1(01) = 1$;

\smallskip

 \ \ \  \ \ \ $\zeta_1(1u)$ \ is undefined for all $u \in \{0,1\}^*$.

\medskip  

\noindent As a gate, $\zeta_1$ has two input wires and one output wire, which
outputs the second input bit if the first input bit is 0; there is no output
if the first input bit is 1. 

\medskip

More generally, for all $m \ge 1$ we define $\,\zeta_m \in$ 
${\sf pfl}{\cal RM}^{\sf fin}_2\,$ to be the right-ideal morphism with 
$\,{\rm domC}(\zeta_m) = \{00, 01\}^m$ and 
$\,{\rm imC}(\zeta_m) = \{0,1\}^m$, such that 

\medskip

 \ \ \  \ \ \ $z_1\,\ldots\,z_m\,u \,\in\, \{00,01\}^m\,\{0,1\}^*$ 
 \ \ \ \big(with $\,z_i \in \{00,01\}\,$ for $i \in [1,m]$, and 
$u \in \{0,1\}^*$\big),

\smallskip

 \ \ \  \ \ \ $\longmapsto$
 \ \ $\zeta_m(z_1\,\ldots\,z_m \,u) \ = \ $ 
$\zeta_1(z_1) \,\ldots\, \zeta_1(z_m) \ u$.

\medskip

\noindent And $\zeta_m(x)$ is undefined if
$x \not\in \{00,01\}^m\,\{0,1\}^*$. 

\smallskip

Then $\zeta_m$ is a composition of $m$ instances of $\zeta_1$ alternating 
with bit-transpositions in $\tau$.
Inductively: 

\smallskip

 \ \ \ $\zeta_m(z_1\,\ldots\,z_{m-1} z_m u) = $
$\zeta_{m-1}(z_1\,\ldots\,z_{m-1}) \ \zeta_1(z_m) \ u$, \ so

\smallskip

 \ \ \ $\zeta_m(.) \,=\, $
$\kappa_{1,m} \circ \zeta_1 \circ \tau_{2,m+1} \circ \tau_{1,m}$ $\circ$ 
$\zeta_{m-1}(.)$, 

\smallskip

\noindent where $\kappa_{1,m}(.)$ is the cyclic permutation seen before.

\begin{pro} \label{PROplepFinGenTau} {\bf (generators for
  ${\sf pfl}{\cal RM}^{\sf fin}_2$).}

\smallskip

\noindent The monoid $\,{\sf pfl}{\cal RM}^{\sf fin}_2$ has a generating set of
the form $\,\Gamma \cup \tau$, where $\,\Gamma$ is a finite subset of 
$\,{\sf pfl}{\cal RM}^{\sf fin}_2$.

\smallskip

\noindent One can choose $\,\Gamma \,=\,$
$\{\zeta_1, {\sf not}, {\sf and}, {\sf or}, {\sf fork}\}$; so this is a 
complete set of gates for acyclic partial circuits.
\end{pro}
{\sc Proof.} Consider $f \in {\sf pfl}{\cal RM}^{\sf fin}_2$ with 
$\,{\rm domC}(f) \subseteq \{0,1\}^m$, and 
$\,f({\rm domC}(f)) = {\rm imC}(f)$ $\,\subseteq\, \{0,1\}^n$.

\smallskip

\noindent Let ${\sf code}_0(.)$ be the encoding of the 3-letter alphabet 
$\{0,1,\bot\}$ into $\{0,1\}^*$, defined as follows:

\smallskip

 \ \ \  \ \ \  \ \ \   ${\sf code}_0(0) = 00$,
 \ \ \ ${\sf code}_0(1) = 01$, \ \ \ ${\sf code}_0(\bot) = 11$;

\smallskip

\noindent and for any string $s = s_1\,\ldots\,s_k \in \{0,1\}^*$ with 
$s_i \in \{0,1\}\,$ (for $i = 1,\ldots,k$):

\smallskip

 \ \ \  \ \ \  \ \ \   ${\sf code}_0(s) =$
${\sf code}_0(s_1) \ \ldots \ {\sf code}_0(s_k)$ \hspace{2,00cm} 
(concatenation).

\smallskip

\noindent So, ${\sf code}_0(.)$ is a monoid morphism of $\{0,1\}^*$, not
a right-ideal morphism.  

\medskip

\noindent We define the following {\em completion} $F$ of $f$ (where $m$ 
and $n$ are as above): \ For all $x \in \{0,1\}^m$,

\medskip

 \ \ \  \ \ \ $F(x) = {\sf code}_0(f(x))$ \ \ \ if
 \ \ $x \in {\rm domC}(f)$ \ ($\,\subseteq \{0,1\}^m$);

\smallskip

 \ \ \  \ \ \ $F(x) = (11)^n$ ($= {\sf code}_0(\bot)^n$) \ \ \ if 
 \ \ $x \in \{0,1\}^m \minus {\rm domC}(f)$; \ and

\smallskip

 \ \ \  \ \ \ $F(xu) = F(x) \ u$, \ \ \  for all $x \in \{0,1\}^m$,
$ \ u \in \{0,1\}^*$.

\medskip

\noindent So $\,{\rm domC}(F) = \{0,1\}^m$, and
$\,{\rm imC}(F) \subseteq \{00, 01, 11\}^n$. 

\smallskip

Then $F$ is total, so $F \in {\sf tfl}{\cal RM}^{\sf fin}_2$.  Therefore by 
Prop.\ \ref{PROPtflGen}, $F$ is generated by $\Gamma_{\sf tfl} \cup \tau$ 
where $\Gamma_{\sf tfl} = \{{\sf not}, {\sf and}, {\sf or}, {\sf fork}\}$,
and $\Gamma_{\sf tfl} \cup \tau$ generates ${\sf tfl}{\cal RM}^{\sf fin}_2$.
Moreover,

\medskip

\noindent $(*)$ \hspace{2,75cm} $f(.) \,=\, \zeta_n \circ F(.)$,

\medskip

\noindent where $\zeta_n \in {\sf pfl}{\cal RM}^{\sf fin}_2$ is the 
right-ideal morphism defined above. Since $\zeta_n$ is a composition of 
instances of $\zeta_1$ and bit-transpositions in $\tau$, we conclude that 
${\sf pfl}{\cal RM}^{\sf fin}_2$ is generated by
$\,\Gamma_{\sf\!tfl} \cup \{\zeta_1\} \cup \tau$.
 \ \ \  \ \ \  $\Box$

\bigskip

\noindent {\bf Remark.}
For ${\sf plep}{\cal RM}^{\sf fin}$ and ${\sf tlep}{\cal RM}^{\sf fin}$, it 
is not known whether they have a generating sets of the form 
$\Gamma \cup \tau$ with $\Gamma$ finite.

We will prove later that ${\sf tlep}{\cal RM}^{\sf fin}$ and 
${\sf pfl}{\cal RM}^{\sf fin}_2$ are not finitely generated.

\bigskip

We now generalize Theorem \ref{circ2lepM} from 
${\sf tfl}{\cal RM}^{\sf fin}_2$ to ${\sf pfl}{\cal RM}^{\sf fin}_2$ by 
using partial circuits.  We saw in Prop.\ \ref{PROplepFinGenTau} that 
${\sf pfl}{\cal RM}^{\sf fin}_2$ is generated by a finite set 
$\Gamma_{\sf\!pfl}$ of partial gates, together with $\tau$. 
In the Theorem below the {\em total length} of $w \in$
$(\Gamma_{\sf\!pfl} \cup \tau)^*$ is defined as in Theorem \ref{circ2lepM}:

\smallskip

\hspace{2,75cm} $|w| \,=\, |w|_{\Gamma_{\sf\!pfl} \cup \tau}$  $\,+\,$
${\rm maxindex}_{\tau}(w)$.

\medskip

\noindent We for an acyclic partial circuit $C$ we define $\ell_{\rm in}(C)$,
$\,\ell_{\rm out}(C)$, $\,|C|$ and ${\sf depth}(C)$ in the same way as for 
boolean circuits (Definitions \ref{DEF_IOlenCircuits} and 
\ref{DEFsizedepthCircuits}).

\begin{thm} \label{circ2PartLepM} {\bf (partial circuits and
${\sf pfl}{\cal RM}^{\sf fin}_2$).}

\smallskip

\noindent Let $\Gamma_{\rm\!pfl} \subseteq {\sf pfl}{\cal RM}^{\sf fin}_2$
be a finite set such that $\Gamma_{\sf\!pfl} \cup \tau$ generates
${\sf pfl}{\cal RM}^{\sf fin}_2$, and $\, \Gamma_{\sf\!pfl} \cap \tau = \0$.
Then we have:

\smallskip

\noindent {\small \rm (1) \ }
For every every $w \in (\Gamma_{\sf\!pfl} \cup \tau)^*$ there exists a 
partial circuit ${\sf C}(w)$ of partial gates such that $w$ and ${\sf C}(w)$
have same input-output function. Moreover, the total length of $w$ and the
size of ${\sf C}(w)$ are linearly related, i.e., 
$ \ |{\sf C}(w)| = \Theta(|w|)$.

\medskip

\noindent {\small \rm (2) \ }
For every partial circuit $c$ there exists a word $\,{\sf W}(c)$ $\in$
$(\Gamma_{\sf\!pfl} \cup \tau)^*$ such that $c$ and ${\sf W}(c)$ have 
same input-output function. moreover,

\smallskip

 \ \ \  \ \ \  The size of $c$ and the total length of 
${\sf W}(c)$ are linearly related: \ $|{\sf W}(c)| = \Theta(|c|)$.

\smallskip

 \ \ \  \ \ \ ${\rm maxindex}_{\tau}({\sf W}(c)) \,\le\, |c|\,$.

\smallskip

 \ \ \  \ \ \ $\ell_{\rm in}({\sf W}(c)) = \ell_{\rm in}(c)$, \ \ and
 \ \ $\ell_{\rm out}({\sf W}(c)) = \ell_{\rm out}(c)$.

\medskip

\noindent {\small \rm (3) \ } The functions ${\sf W}(.)$ and ${\sf C}(.)$ 
can be chosen so that for every partial circuit $c:$
 \ ${\sf C}({\sf W}(c)) \ \widehat{=} \ c$.
\end{thm}
{\sc Proof.} This can be proved in the same way as Theorem \ref{circ2lepM}. 
Indeed, using a set $\Gamma_{\sf\!pfl}$ of partial gates and generators 
instead of a set $\Gamma_{\sf\!tfl}$ of total ones, does not change the 
proof.
 \ \ \  $\Box$

\bigskip

\noindent {\bf Remark.} Most of the time we use the alphabet $A = $
$\{0,1,\ldots,k\!-\!1\}$, for some fixed $k \ge 2$, and the Higman-Thompson
monoids and groups. However for ${\sf tfl}{\cal RM}^{\sf fin}_2$ and
${\sf pfl}{\cal RM}^{\sf fin}_2$ we only consider the alphabet $\{0,1\}$, 
i.e., $k = 2$. This is because of the connection between these two monoids 
and circuits.  

\bigskip

\noindent Theorems \ref{circ2lepM} and \ref{circ2PartLepM} are used for 
the following definition.

\begin{defn} \label{DEFdepthRM} {\bf (depth in 
${\sf tfl}{\cal RM}^{\sf fin}_2$ and ${\sf pfl}{\cal RM}^{\sf fin}_2$).}
 \ Let $\Gamma_{\sf\!tfl}$ and $\Gamma_{\sf\!pfl}$ be finite sets such that
$\Gamma_{\sf\!tfl} \cup \tau$ generates ${\sf tfl}{\cal RM}^{\sf fin}_2$ and 
$\Gamma_{\sf\!pfl} \cup \tau$ generates ${\sf pfl}{\cal RM}^{\sf fin}_2$.
Then for $w \in (\Gamma_{\sf\!tfl} \cup \tau)^*$ or $w \in$ 
$(\Gamma_{\sf\!pfl} \cup \tau)^*\,$ we define the depth of $w$ by

\smallskip 

\hspace{4,1cm}  ${\sf depth}(w) \,=\, {\sf depth}({\sf C}(w))$.  
\end{defn}

\smallskip

\section{Various definitions of various Thompson monoids}

The Thompson monoids will be defined as quotient monoids of the monoids of
right-ideal morphisms seen above. For this we define several congruence 
relations, that are themselves based on equivalence relations between prefix 
codes. For finite prefix codes, these equivalence relations turn out to be
the same.

\subsection{Equivalence relations between prefix codes}

We will use the following equivalence relations between prefix codes
$P, Q \subset A^*$. The equivalences $\equiv_{\rm end}$, $\equiv_{\rm bd}$, 
and $\equiv_{\rm poly}$ were studied in \cite{Equiv};
$\,\equiv_{\rm fin}$ was introduced in \cite[{\small Def.\ 2.12}]{BinMk1};
$\,\equiv_{0^{\omega}}$ is new.

\medskip

\bigskip

\noindent $\bullet$ \ {\boldmath $[\,\equiv_{\rm end}\,]$}
 \ \ \ {\em $P \equiv_{\rm end} Q$ \ \ iff \ \ for every right ideal
$R \subseteq A^*: \ \big[\,R \,\cap\, P A^* = \0$
 \ $\Leftrightarrow$ \ $R \,\cap\, Q A^* = \0\,$\big]. }
 \ \ This is called the {\em end-equivalence}.

\medskip

\noindent By \cite[{\small Prop.\ 2.5}]{Equiv}, $P \equiv_{\rm end} Q\,$ is 
equivalent to 

\medskip

 \ \ \  \ \ \ ${\sf cl}(P A^{\omega}) \,=\, {\sf cl}(Q A^{\omega})$, 

\medskip

\noindent where ${\sf cl}(.)$ denotes {\em closure} in the Cantor space 
topology of $A^{\omega}$.

\smallskip

The relation $\equiv_{\rm end}$ can be extended to a preorder on the set
of prefix codes of $A^*$ by defining $\,P \le_{\rm end} Q$ \ iff \ for
every right ideal
$R \subseteq A^*: \ [R \,\cap\, Q A^* = \0$
 \ $\Rightarrow$ \ $R \,\cap\, P A^* = \0]$.  \ Equivalently,
$\,{\sf cl}(P A^{\omega}) \subseteq {\sf cl}(Q A^{\omega})$.
 \ Then $\,P \equiv _{\rm end} Q$ \ iff \ $P \le_{\rm end} Q\,$ and
$\,Q \le_{\rm end} P$.

\begin{lem} 
 \ Let $P, Q \subset A^*$ be prefix codes. If $P$ is a {\em maximal} prefix 
code and $P \equiv_{\rm end} Q$, then $Q$ is maximal.
\end{lem} 
{\sc Proof.} $\,P$ is maximal iff $\,{\sf cl}(PA^{\omega}) = A^{\omega}$.
And since $P \equiv_{\rm end} Q$, $\,{\sf cl}(QA^{\omega})$ $=$
${\sf cl}(PA^{\omega})$, so $\,{\sf cl}(QA^{\omega}) = A^{\omega}$.
 \ \ \ $\Box$

\bigskip

\noindent $\bullet$ \ {\boldmath $[\,\equiv_{\rm bd}\,]$}
 \ \ \ {\em $P \equiv_{\rm bd} Q$ \ \ iff 
 \ \  $P A^{\omega} = Q A^{\omega}$.}
 \ \ This is called {\em bounded end-equivalence}.

\medskip

\noindent By \cite[{\small Prop.\ 3.2}]{Equiv}, 
$\,P A^{\omega} = Q A^{\omega}\,$ is equivalent to the following: 
(1) $\,P \equiv_{\rm end} Q$, and (2) there exists a total function 
$\,r$: ${\mathbb N} \to {\mathbb N}\,$ such that for
all $x_1 \in P$ and $x_2 \in Q$:

\smallskip

 \ \ \  \ \ \ if $\, x_1 \parallel_{\rm pref} x_2$ \ then
 \ $|x_1| \le r(|x_2|)$ and $|x_2| \le r(|x_1|)$.

\smallskip

\noindent Here, $\parallel_{\rm pref}$ denotes comparability in the prefix
order $\le_{\rm pref}$ of $A^*$.

\smallskip

It follows that $\,P \equiv_{\rm bd} Q\,$ implies $\,P \equiv_{\rm end} Q$.
The converse does not hold; e.g., $0^*1 \equiv_{\rm end} \{\e\}$, but
$0^*1 \not\equiv_{\rm bd} \{\e\}$.

$P \equiv_{\rm bd} Q$ does not imply $\,P A^* = Q A^*$.
 \ E.g., $0\,\{0,1\} \equiv_{\rm bd} \{0\}$, but
$\,0\,\{0,1\} \,\{0,1\}^* = 0\,\{0,1\}^+$  $\ne$  $0\,\{0,1\}^*$.

\smallskip

The relation $\equiv_{\rm bd}$ can be extended to a preorder on the set
of prefix codes of $A^*$ by defining $\,P \le_{\rm bd} Q$ \ iff \
$P A^{\omega} \subseteq Q A^{\omega}$.
Then $\,P \equiv_{\rm bd} Q$ \ iff \ $P \le_{\rm bd} Q\,$ and
$\,Q \le_{\rm bd} P$.

\bigskip

\noindent $\bullet$ {\boldmath $[\,\equiv_{\rm poly}\,]$}
 \ \ \ {\em $P \equiv_{\rm poly} Q$ \ \ iff
 \ \ $P \equiv_{\rm bd} Q$, and there exists a polynomial $\pi$ such that    
for all $x_1 \in P$ and $x_2 \in Q$:
 \ \ if $ \ x_1 \parallel_{\rm pref} x_2 \, $ then
$ \ |x_1| \le \pi(|x_2|)\,$ and $\,|x_2| \le \pi(|x_1|)$.
}

\medskip

The relation $\equiv_{\rm poly}$ can be extended to a preorder on the set
of prefix codes of $A^*$ by defining $\,P \le_{\rm poly} Q$ \ iff \
$P A^{\omega} \subseteq Q A^{\omega}$, and there exists a polynomial $p$     
such that for all $x_1 \in P$ and $x_2 \in Q$:

\smallskip

 \ \ \  \ \ \  if $\, x_1 \parallel_{\rm pref} x_2 \, $ then
$|x_1| \le \pi(|x_2|)$ and $|x_2| \le \pi(|x_1|)$.

\bigskip

\noindent $\bullet$ {\boldmath $[\,\equiv_{\rm fin}\,]$}
  \ \ \ {\em $P \equiv_{\sf fin} Q$ \ \ iff \ \ the symmetric difference 
$\,PA^* \vartriangle QA^*\,$ is finite. } 

\medskip

The relation $\equiv_{\rm fin}$ can be extended to a preorder on the set
of prefix codes of $A^*$ by defining $\,P \le_{\rm fin} Q$ \ iff \
$P A^{\omega} \subseteq Q A^{\omega}$, and the set difference
$\,P A^* \minus Q A^*\,$ is finite. 
 \ Then $P \equiv_{\rm fin} Q\,$ iff $\,P \le_{\rm fin} Q\,$  and
$\,Q \le_{\rm fin} P$. 

\smallskip

\noindent We give another characterization of $\equiv_{\rm fin}\,$.

\begin{defn} 
 \ Let $\,\to_1\,$ be the relation between prefix codes $P, Q \subset A^*$
defined as follows: \ 

\smallskip

 \ \ \  \ \ \ $P \to_1 Q$ \ \ iff \ \ there exists $p \in P$ such that
$ \ Q \,=\, (P \minus \{p\}) \,\cup\, pA$. 

\smallskip

\noindent This is called a {\em one-step restriction} of $P\,$ (see 
{\rm \cite[{\small Lemma 2.22}]{BinMk1}}.). 

\smallskip

\noindent Let $\,\equiv_1^*\,$ be the reflexive-symmetric-transitive closure
of $\,\to_1\,$.
\end{defn}

\begin{pro} \label{equi1setdiff}
 \ The equivalence relations $\,\equiv_{\rm fin}\,$ and $\,\equiv_1^*\,$
are the same.
\end{pro}
{\sc Proof.} This is a special case of \cite[{\small Lemma 2.25}]{BinMk1} 
for $n = 1$. 
 \ \ \ $\Box$

\bigskip

\noindent $\bullet$ {\boldmath $[\,\equiv_{0^{\omega}}\,]$}
  \ \ \ {\em $P \equiv_{0^{\omega}} Q$ \ \ iff
 \ \ $PA^*\,0^{\omega} \,=\, QA^*\,0^{\omega}$. 
  \ \ \  \ \ \  (Here we assume that $\,0 \in A$.) }

\bigskip

\noindent {\bf Remark.} Let $A = \{0,1,\ldots,k\!-\!1\}$, and consider the
base-$k$ representation of rational numbers in the semi-open interval
$\,[0,1[\,$. Then the set $\,\{0.x_1 \,\ldots\,x_n\,0^{\omega} : $ 
$x_1, \,\ldots\, , x_n \in A$, $n \ge 1\}$ \, represents all the $k$-ary 
rational numbers in $\,[0,1[\,$ in a unique way. \ (By definition, a 
$k$-ary rational number is a rational number of the form $\,m/k^n\,$ for any
$m \in {\mathbb Z}, \, n \in {\mathbb N}$.) \ So we have:

\smallskip

 \ \ \  $P \equiv_{0^{\omega}} Q$ \ \ iff 
 \ \ $\{0.pu :  p \in P,\,u \in A^*\}$  $=$  
$\{0.qv :  q \in Q,\,v \in A^*\}\,$ \ ($\,\subset {\mathbb Q}$). 

\bigskip

\noindent We now start studying the relationship between the above 
equivalence relations.

\begin{lem} \label{EquFinProperty}
 \ For {\em finite} prefix codes $P, Q \subset A^*$:
 \ \ $P \equiv_{\rm fin} Q$ \ \ iff \ \ $P \equiv_{\rm bd} Q$ 
 \ \ iff \ \ $P \equiv_{\rm end} Q$.
\end{lem}
{\sc Proof.} The first equivalence is a special case of 
\cite[{\small Lemma 2.13}]{BinMk1}, for $n = 1$. The second equivalence
holds since for every finite prefix code $P$, $ \ PA^{\omega} =$ 
${\sf cl}(PA^{\omega}) \ $ in Cantor space.
 \ \ \ $\Box$

\begin{lem} \label{LEMequivend010} {\bf \ {\boldmath ($\,\equiv_{\rm bd} \ $
$\Rightarrow$  $ \ \equiv_{0^{\omega}} \ $  $\Rightarrow$
$ \ \equiv_{\rm end}\,$).}}

\smallskip

\noindent {\small \rm (1)} For all prefix codes $\,P, Q \subseteq A^*$:
 \ \ \ $P \equiv_{\rm bd} Q$ \ \ $\Rightarrow$
 \ \ $P \equiv_{0^{\omega}} Q$ \ \ $\Rightarrow$
 \ \ $P \equiv_{\rm end} Q$.

\medskip

\noindent {\small \rm (2)} If $P$ and $Q$ are {\em finite}, these 
implications are equivalences ($\Leftrightarrow$). 

But for prefix codes in general (even for finite-state ones), the 
implications are not equivalences.
\end{lem}
{\sc Proof.} (1) [First $\Rightarrow$]
 \ Suppose $P A^{\omega} = Q A^{\omega}$.  Then
$P A^{\omega} \cap A^*0^{\omega} = Q A^{\omega} \cap A^*0^{\omega}$.
And it is easy to prove that $P A^{\omega} \cap A^*0^{\omega} =$
$P A^*\,0^{\omega}$, and similarly for $Q$.

\smallskip

\noindent [Second  $\Rightarrow]$ \ Suppose
$P A^*\,0^{\omega} = Q A^*\,0^{\omega}$, and consider any right-ideal
$R \subseteq A^*$. If $PA^* \cap R \ne \0$ then there exist $p \in P$ and 
$u \in A^*$ such that $pu \in R$. Since 
$P A^*\,0^{\omega} = Q A^*\,0^{\omega}$, there are $q \in Q$ and $v \in A^*$
such that $pu\,0^{\omega} = qv\,0^{\omega}$, hence $pu\,0^i = qv\,0^j$
for some $i,j \in {\mathbb N}$. Since $R$ is a right ideal, $pu\,0^i \in R$,
hence $qv\,0^j \in R$. So, $QA^* \cap R \ne \0$.
In a similar way one proves that if $QA^* \cap R \ne \0$ then
$PA^* \cap R$  $\ne$  $\0$.  Hence, $P \equiv_{\rm end} Q$.

\smallskip

\noindent (2) We saw in Lemma \ref{EquFinProperty} that if $P$ and $Q$ are
finite then $\equiv_{\rm bd}$ and $\equiv_{\rm end}$ are the same. Hence, 
since $\equiv_{0^{\omega}}$ is between $\equiv_{\rm bd}$ and 
$\equiv_{\rm end}$, $\,\equiv_{0^{\omega}}$ is also equal to them.

\smallskip

In general, $\equiv_{0^{\omega}}$ is different from $\equiv_{\rm end}$.
 \ Example: For $A = \{0,1\}$, Consider $0^*1$ and $\{\e\}$. 
We have $\,0^* 1 \equiv_{\rm end} \{\e\}$; but
$0^{\omega}$ $\in$  $\{0,1\}^*0^{\omega} \minus 0^*1\,\{0,1\}^*\,0^{\omega}$,
so $\,0^*1 \not\equiv_{0^{\omega}} \{\e\}$.

\smallskip

In general, $\,\equiv_{0^{\omega}}\,$ is different from $\equiv_{\rm bd}$. 
 \ Example: For $A = \{0,1\}$, Consider $\,1^*\,0\,$ and $\{\e\}$.
We have $\,1^*\,0 \not\equiv_{\rm bd} \{\e\}$, since
$\,1^{\omega} \not\in 1^*\,0\,\{0,1\}^{\omega}$.
However we can prove that 
 \ $1^*\,0\,\{0,1\}^*\,0^{\omega} = \{0,1\}^*\,0^{\omega}$, so
$\,1^*\,0 \equiv_{0^{\omega}} \{\e\}$.

\noindent Proof that 
$\,1^*\,0\,\{0,1\}^*\,0^{\omega} = \{0,1\}^*\,0^{\omega}$:

\noindent We obviously have ``$\subseteq$''. Let us prove ``$\supseteq$''. 
For any element $u\,0^{\omega} \in \{0,1\}^*\,0^{\omega}$
we have the following cases:  \\
(a) \ $u \in 0^*$: \ Then $\,u\,0^{\omega} = 0^{\omega}$  $\in$
$\,1^*\,0\,\{0,1\}^*\,0^{\omega}$. \\
(b) \ $u \in 1^*$: \ Then $\,u\,0^{\omega} = 1^n 0^{\omega} \in$
$1^*\,0\,\{0,1\}^*\,0^{\omega}$.  \\
(c.1) \ $u$ contains both 0 and 1, and starts with 1: \ Then
$\,u\,0^{\omega} = 1^n 0 v \,0^{\omega} \in 1^*\,0\,\{0,1\}^*\,0^{\omega}$. \\
(c.2) \ $u$ contains both 0 and 1, and starts with 0: \ Then
$\,u\,0^{\omega} = 0 v \,0^{\omega}$ $=$
$\e 0v\,0^{\omega} \in 1^*\,0\,\{0,1\}^*\,0^{\omega}$.
 \ \ \  \ \ \ $\Box$

\bigskip

\noindent {\bf Remark.} 
Since $\equiv_{0^{\omega}}$ implies $\equiv_{\rm end}$, we could have used
$\equiv_{0^{\omega}}$ instead of $\equiv_{\rm end}$ in the definitions of
$\equiv_{\rm bd}$ and $\equiv_{\rm poly}$ (given above).

\begin{pro} \label{EquivBDonFin} {\bf (preservation of finiteness).}
 \ Let $P, Q \subset A^*$ be prefix codes. 

If $P$ is finite and $\,P \equiv_{\rm bd} Q$, then $Q$ is finite.

This does not hold in general for $\,\equiv_{0^{\omega}}\,$ (and hence not 
for $\,\equiv_{\rm end}$).
\end{pro}
{\sc Proof.} The first claim is proved in \cite[{\small Lemma 3.15}]{Equiv}.

An example where this does not hold for $\equiv_{0^{\omega}}$ is 
$\,1^*\,0 \equiv_{0^{\omega}} \{\e\}\,$ (used in the proof of Lemma 
\ref{LEMequivend010}(2)). This implies $\,1^*\,0 \equiv_{\rm end} \{\e\}$.
  \ \ \ \ \ \  $\Box$

\begin{pro} \label{EquFinEquivalent}
 \ For {\em finite} prefix codes $P, Q \subset A^*$ the following are
equivalent:

\smallskip

 \ \ \ $P \equiv_{\rm fin} Q$,  
 \ \  $P \equiv_{\rm poly} Q$, \ \ $P \equiv_{\rm bd} Q$,
 \ \ $P \equiv_{0^{\omega}} Q$, \ \ $P \equiv_{\rm end} Q$. 
\end{pro}
{\sc Proof.} The left-to-right implications hold for all prefix codes:
The first two are obvious, and the last three hold by Lemma 
\ref{LEMequivend010}(1).  \ In the other direction:
For finite prefix codes, $\,P \equiv_{\rm end} Q\,$ implies 
$P \equiv_{\rm fin} Q\,$ by Lemma \ref{EquFinProperty}.
  \ \ \ $\Box$

\subsection{Congruence relations on ${\cal RM}^{\sf fin}$ }

In \cite{Equiv} we extended the equivalence relations
$\equiv_{\rm end}$, $\equiv_{\rm bd}$ and $\equiv_{\rm poly}$ from prefix
codes to right-ideal morphisms. We do this more generally:

\begin{defn} \label{DEFequFinMorph}
 \ Let $\equiv_{\rm X}$ is  any one of $ \ \equiv_{\rm end}$, 
$\,\equiv_{0^{\omega}}$, $\,\equiv_{\rm bd}$, $\,\equiv_{\rm poly}$, 
or $\,\equiv_{\rm fin}$
 \ For right-ideal morphisms $f, g$ of $A^*$ we define 
$\,f \equiv_{\rm X} g\,$ by

\medskip

 \ \ \  $f \equiv_{\rm X} g$ \ \ iff
 \ \ ${\rm domC}(f) \equiv_{\rm X} {\rm domC}(g)$, and
$\,f(x) = g(x)\,$ for all $\,x \in {\rm Dom}(f) \cap {\rm Dom}(g)$.
\end{defn}
In \cite{Equiv} we proved that $\equiv_{\rm end}$, $\equiv_{\rm bd}$ and
$\equiv_{\rm poly}$ are congruences on the ${\cal RM}^{\sf P}$. 
In the rest of this subsection we prove now that $\equiv_{0^{\omega}}$ and 
$\equiv_{\rm fin}$ are also congruences.

\medskip

In the next Proposition the following is used: By definition, a function 
$f: A^* \to A^*$ is \textbf{\textit{balanced}} \ iff \ there exists a 
total function $\pi$: ${\mathbb N} \to {\mathbb N}$ such that for all 
$x \in$ ${\rm Dom}(f):$ \ $|f(x)| \le \pi(|x|)$ \ and
 \ $|x| \le \pi(|f(x)|)$. 

One can show that $f$ is balanced iff for
all $y \in {\rm Im}(f): \ f^{-1}(y)$ is finite.

\begin{pro} \label{EquFinCongr}
 \ The relation $\,\equiv_{\rm fin}$ is a {\em congruence} on the monoid of
all balanced right-ideal morphisms of $A^*$.
\end{pro}
{\sc Proof.} We want to show that if $f, g, h$ are balanced right-ideal
morphisms, and $f \equiv_{\rm fin} g$, then $f h \equiv_{\rm fin} g h$
and $h f \equiv_{\rm fin} h g$.
It follows from $f \equiv_{\rm fin} g$ that $f \equiv_{\rm bd} g$ and
that $\, {\rm Dom}(f) \vartriangle {\rm Dom}(g) \,$ is finite. Since
$\equiv_{\rm bd}$ is a congruence (as was proved in \cite{Equiv}), it
follows that $fh \equiv_{\rm bd} gh$ and $hf \equiv_{\rm bd} hg$. Hence it
will suffice to show that
$\, {\rm Dom}(f h) \vartriangle {\rm Dom}(gh) \,$ and
$\, {\rm Dom}(hf) \vartriangle {\rm Dom}(hg) \,$ are finite.

For $hf$ and $hg$, consider any $x \in {\rm Dom}(hf) \minus {\rm Dom}(hg)$,
i.e., $hf(x)$ is defined but $hg(x)$ is not defined. It is not possible for
both $f(x)$ and $g(x)$ to be defined (otherwise, $f(x) = g(x)$ since
$f \equiv_{\rm fin} g$, and then both $hf(x)$ and $hg(x)$ would be defined).
Hence, $f(x)$ is defined but $g(x)$ is not defined.  Since
$\, {\rm Dom}(f) \vartriangle {\rm Dom}(g) \,$ is finite, there are only
finitely many such $x$, so $\, {\rm Dom}(hf) \minus {\rm Dom}(hg)$ is        finite.
In a similar way one proves that $\, {\rm Dom}(hg) \minus {\rm Dom}(hf) \,$
is finite.

For $fh$ and $gh$, consider any $x \in {\rm Dom}(fh) \minus {\rm Dom}(gh)$,
i.e., $fh(x)$ (hence $h(x)$) is defined but $gh(x)$ is not defined.
Hence, $h(x) \in {\rm Dom}(f) \minus {\rm Dom}(g)$, which is finite. Thus,
there are only finitely many choices for $h(x)$. Since $h$ is balanced it 
follows that there are only finitely many choices for $x \in$ $h^{-1} h(x)$,
so $\, {\rm Dom}(fh) \minus {\rm Dom}(gh) \,$ is finite. In a similar way 
one proves that $\, {\rm Dom}(gh) \minus {\rm Dom}(fh) \,$ is finite.
  \ \ \  \ \ \ $\Box$

\bigskip

\noindent In order to prove that $\equiv_{0^{\omega}}$ is a congruence we
use two lemmas.

\begin{lem} \label{LEM0omegaInters}
 \ For any prefix codes $P, Q \subseteq A^*$:
 \ \ $(P A^* \cap Q A^*)\,0^{\omega}$  $=$
$P A^*\,0^{\omega} \,\cap\, Q A^*\,0^{\omega}$.
\end{lem}
{\sc Proof.} $[\subseteq]\,$ Since $P A^* \cap Q A^* \subseteq P A^*$, we
have $(P A^* \cap Q A^*)\,0^{\omega}$  $\subseteq$
$P A^*\,0^{\omega}$; similarly $(P A^* \cap Q A^*)\,0^{\omega}$  $\subseteq$
$Q A^*\,0^{\omega}$. Hence, $(P A^* \cap Q A^*)\,0^{\omega}$  $\subseteq$ 
$P A^*\,0^{\omega} \,\cap\, Q A^*\,0^{\omega}$.

\noindent $[\supseteq]\,$ If $x \in$
$P A^*\,0^{\omega} \cap Q A^*\,0^{\omega}$ then $x = pu0^{\omega}$ $=$
$qv 0^{\omega}$, for some $p \in P$, $q \in Q$ and $u, v \in A^*$.
Hence $pu0^i = qv0^j$ for some $i, j \in {\mathbb N}$.
So $x = pu0^i 0^{\omega} = qv0^j 0^{\omega}$ for some $pu0^i = qv0^j \in$
$P A^* \cap Q A^*$; hence $x \in (P A^* \cap Q A^*)\,0^{\omega}$.
 \ \ \  \ \ \ $\Box$

\begin{lem} \label{LEM0omegahinv}
 \ For any prefix codes $P, Q \subseteq A^*$ and any right-ideal morphism
$h$ of $A^*$:

\smallskip

 \ \ \  \ \ \ $P A^*\,0^{\omega} = Q A^*\,0^{\omega}$ \ \ implies
 \ \ $h^{-1}(P A^*) \,0^{\omega} = h^{-1}(Q A^*) \ 0^{\omega}$.
\end{lem}
{\sc Proof.} By Lemma \ref{LEM0omegaInters},
$P A^*\,0^{\omega} = Q A^*\,0^{\omega}$  $=$
$(P A^* \cap Q A^*)\,0^{\omega}$.  Our result will follow if we prove

\smallskip

 \ \ \  \ \ \ $h^{-1}(P A^*) \ 0^{\omega}$  $\,=\,$
    $h^{-1}(P A^* \cap Q A^*) \ 0^{\omega}$,

\smallskip

\noindent and similarly for $Q$; then indeed, \ $h^{-1}(P A^*) \ 0^{\omega}$ 
$\,=\,$ $h^{-1}(P A^* \cap Q A^*) \ 0^{\omega}$ $\,=\,$
$h^{-1}(Q A^*) \ 0^{\omega}$.

\smallskip

It is obvious that $\,h^{-1}(P A^* \cap Q A^*) \ 0^{\omega}$  $\subseteq$
$h^{-1}(P A^*) \ 0^{\omega}$.

Let us prove that $\,h^{-1}(P A^*) \ 0^{\omega}$
$\subseteq$  $h^{-1}(P A^* \cap Q A^*) \ 0^{\omega}$.
Any $x \in h^{-1}(P A^*) \ 0^{\omega}$ is of the form $x = z 0^{\omega}$
where $h(z) = pu$ for some $p \in P$ and $u \in A^*$.
Then, since $P A^*\,0^{\omega} = Q A^*\,0^{\omega}$  we have:
$h(z)\,0^{\omega} = pu0^{\omega} = qv 0^{\omega}$, for some $q \in Q$
and $v \in A^*$. Hence, $pu0^i = qv0^j$ for some $i, j \in {\mathbb N}$.
Now, $h(z 0^i) = pu0^i = qv0^j \in P A^* \cap Q A^*$, so
$z 0^i \in h^{-1}(P A^* \cap Q A^*)$. Therefore,
$x = z 0^i \,0^{\omega} \in h^{-1}(P A^* \cap Q A^*) \ 0^{\omega}$.
 \ \ \  \ \ \ $\Box$

\begin{pro} \label{PROPequiv0omegaCongr}
 \ The relation $\equiv_{0^{\omega}}$ is a {\em congruence} on the monoid
of all right-ideal morphisms of $A^*$ (and in particular on the monoids
${\cal RM}^{\sf fin}$ and ${\cal RM}^{\sf P}$).
\end{pro}
{\sc Proof.} Let $f, g, h$ be right-ideal morphisms of $A^*$ such that
$f \equiv_{0^{\omega}} g$. We want to show that
$fh(.) \equiv_{0^{\omega}} gh(.)$ and $hf(.) \equiv_{0^{\omega}} hg(.)$.
Since $f \equiv_{0^{\omega}} g$ implies $f \equiv_{\rm end} g\,$ by Lemma
\ref{LEMequivend010}, and $\equiv_{\rm end}$ is a congruence by
\cite[{\small Prop.\ 2.17}]{Equiv}, we conclude:
$fh(.) \equiv_{\rm end} gh(.)$ and $hf(.) \equiv_{\rm end} hg(.)$.
The Proposition then follows from the Claims below.

\smallskip

\noindent {\sf Claim 1.} \ \ ${\rm Dom}(fh(.)) \ 0^{\omega}$  $\,=\,$
${\rm Dom}(gh(.)) \ 0^{\omega}$.

\smallskip

\noindent Proof of Claim 1:  By \cite[{\small Prop.\ 3.5}]{BinMk1},
$\,{\rm Dom}(fh(.)) = {\rm Dom}(h) \,\cap\, h^{-1}({\rm Dom}(f)$. By Lemma
\ref{LEM0omegaInters}, this implies that
$ \ {\rm Dom}(fh(.))\,0^{\omega}$  $=$
${\rm Dom}(h)\,0^{\omega} \,\cap\, h^{-1}({\rm Dom}(f) \ 0^{\omega}$.
By Lemma \ref{LEM0omegahinv}, and since $f \equiv_{0^{\omega}} g$, this
is equal to
$ \ {\rm Dom}(h)\,0^{\omega} \,\cap\, h^{-1}({\rm Dom}(g) \ 0^{\omega}$;
and this is ${\rm Dom}(gh(.))$.
 \ \ \  \ \ \ [End, proof of Claim 1]

\smallskip

\noindent {\sf Claim 2.} \ \ $fh(.)$ and $gh(.)\,$ agree on
$\,{\rm Dom}(fh(.)) \,\cap\, {\rm Dom}(gh(.))$.

\smallskip

\noindent Proof of Claim 2:  By Lemma \ref{LEMequivend010},
$\equiv_{0^{\omega}}$ implies $\equiv_{\rm end}$ on all prefix codes; and the
requirement of agreement on the intersection of the domains is the same for
$\equiv_{0^{\omega}}$ and $\equiv_{\rm end}$.
And $\equiv_{\rm end}$ is a congruence (by \cite[{\small Prop.\ 2.
17}]{Equiv}); now, since $fh(.) \equiv_{\rm end} gh(.)$, $fh(.)$ and $gh(.)$
agree on $\,{\rm Dom}(fh(.)) \,\cap\, {\rm Dom}(gh(.))$.
 \ \ \  \ \ \ [End, proof of Claim 2]

\smallskip

\noindent Claims 1 and 2 imply that if $\,f \equiv_{0^{\omega}} g\,$ then
$\,fh(.) \equiv_{0^{\omega}} gh(.)$.

\medskip

\noindent {\sf Claim 3.} \ \ ${\rm Dom}(hf(.)) \ 0^{\omega}$  $\,=\,$
${\rm Dom}(hg(.)) \ 0^{\omega}$.

\smallskip

\noindent Proof of Claim 3:  By \cite[{\small Prop.\ 3.5}]{BinMk1},
$\,{\rm Dom}(hf(.)) = {\rm Dom}(f) \,\cap\, f^{-1}({\rm Dom}(h))$, so
$\,{\rm Dom}(hf(.))\,0^{\omega}$  $=$
${\rm Dom}(f)\,0^{\omega} \,\cap\, f^{-1}({\rm Dom}(h)) \ 0^{\omega}$.
Hence for every $x \in {\rm Dom}(hf(.))$ we have
$x \in {\rm Dom}(f)\,0^{\omega} = {\rm Dom}(g)\,0^{\omega}$, and
$x = z 0^{\omega}$ for some $z \in A^*$ such that $f(z) \in {\rm Dom}(h)$.
Since ${\rm Dom}(h)$ is a right ideal,
$f(z 0^i) = f(z)\,0^i \in {\rm Dom}(h)$, for all $i \ge 0$.
Since ${\rm Dom}(f)\,0^{\omega} = {\rm Dom}(g)\,0^{\omega}$ and
${\rm Dom}(f)$ and ${\rm Dom}(g)$ are right ideals, we have
$z 0^i \in {\rm Dom}(f) \cap {\rm Dom}(g)$ for all $i$ large enough.  Now,
$x = z0^i 0^{\omega}$, where $z 0^i \in {\rm Dom}(f) \cap {\rm Dom}(g)$;
moreover, $f(z0^i) = g(z0^i)$, since $f$ and $g$ agree on
${\rm Dom}(f) \cap {\rm Dom}(g)$.
Hence, $f(z)\,0^{\omega} \in {\rm Dom}(h)$ is equivalent to
$g(z)\,0^{\omega} \in {\rm Dom}(h)$.
Therefore, $x \in {\rm Dom}(g)\,0^{\omega} \,\cap\,$
$g^{-1}({\rm Dom}(h)) \ 0^{\omega}$.
Thus, $\,{\rm Dom}(hf(.)) \subseteq {\rm Dom}(hg(.))$.
In the same way one proves that $\,{\rm Dom}(hg(.))$  $\subseteq$
${\rm Dom}(hf(.))$.   \ \ \  \ \ \ [End, proof of Claim 3]

\smallskip

\noindent {\sf Claim 4.} \ \ $hf(.)$ and $hg(.)\,$ agree on
$\,{\rm Dom}(hf(.)) \,\cap\, {\rm Dom}(hg(.))$.

\smallskip

\noindent Claim 4 is proved in the same way as Claim 2.

\smallskip

\noindent Claims 3 and 4 imply that if $f \equiv_{0^{\omega}} g$ then
$hf(.) \equiv_{0^{\omega}} hg(.)$.
 \ \ \  \ \ \  $\Box$

\subsection{Defining Thompson monoids}

The monoid $M_{k,1}$ is defined (see \cite{BiThompsMonV3}) as the quotient
monoid

\medskip

 \ \ \  \ \ \  $M_{k,1} \ = \ {\cal RM}^{\sf fin}/\!\!\equiv_{\rm fin}\,$. 

\medskip

\noindent
On ${\cal RM}^{\sf fin}$ the congruence $\equiv_{\rm fin}$ is equal to
various other congruences, as seen in Prop.\ \ref{EquFinEquivalent}. 
Therefore, 

\medskip

 \ \ \  \ \ \ $M_{k,1} \ = \ {\cal RM}^{\sf fin}/\!\!\equiv_{\rm fin}$
$ \ = \ {\cal RM}^{\sf fin}/\!\!\equiv_{\rm poly}$
$ \ = \ {\cal RM}^{\sf fin}/\!\!\equiv_{\rm bd}$
$ \ = \ {\cal RM}^{\sf fin}/\!\!\equiv_{0^{\omega}}$
$ \ = \ {\cal RM}^{\sf fin}/\!\!\equiv_{\rm end}\,$.

\medskip

\noindent Thompson's original monoid can be defined by

\medskip

 \ \ \  \ \ \ ${\sf tot}M_{k,1} \,=\,$
${\sf tot}{\cal RM}^{\sf fin}/\!\!\equiv_{\rm fin}\,$,

\smallskip

\noindent where the same monoid is obtained if $\equiv_{\rm fin}$ is 
replaced by the other congruences. (Thompson originally defined this monoid
as the monoid generated by a certain finite set of generators acting on the
Cantor space $\{0,1\}^{\omega}$, and on other sets.)

\medskip

\noindent The Thompson group $V$ can be defined by 

\medskip

 \ \ \  \ \ \ $V$
$ \,=\, {\sf totsurinj}{\cal RM}^{\sf fin}/\!\!\equiv_{\rm fin}\,$,

\smallskip

\noindent where the same monoid is obtained if $\equiv_{\rm fin}$ is 
replaced by the other congruences.
 \ The Thompson group $V$ is the {\em group of units} of each one of the
monoids $M_{k,1}$, ${\sf tot}M_{k,1}$, ${\sf sur}M_{k,1}$, 
${\sf inj}M_{k,1}$, and their intersections. 

\bigskip

\noindent We also define 

\medskip

 \ \ \  \ \ \ ${\sf tfl}M_{2,1}$  
$ \,=\, {\sf tfl}{\cal RM}^{\sf fin}_2/\!\!\equiv_{\rm fin}\,$, \ \ and

\medskip

 \ \ \  \ \ \ ${\sf pfl}M_{2,1}$  
$ \,=\, {\sf pfl}{\cal RM}^{\sf fin}_2/\!\!\equiv_{\rm fin}\,$,

\medskip

\noindent where the same monoid is obtained if $\equiv_{\rm fin}$ is
replaced by the other congruences. These are called the Thompson monoids of 
total, or partial, \textbf{\textit{fixed-length input}} functions.

These Thompson monoids have another characterization:

\begin{lem} \label{lepVSfl} 
 \ The monoids   
 \ ${\sf plep}{\cal RM}^{\sf fin}/\!\!\equiv_{\rm fin}$ \ and
 \ ${\sf pfl}{\cal RM}^{\sf fin}_2/\!\!\equiv_{\rm fin}\,$ \ are isomorphic.

The monoids \ ${\sf tlep}{\cal RM}^{\sf fin}/\!\!\equiv_{\rm fin}$ \ and
 \ ${\sf tfl}{\cal RM}^{\sf fin}_2/\!\!\equiv_{\rm fin}\,$ are isomorphic.
\end{lem}
{\sc Proof.} For every $f \in {\sf plep}{\cal RM}^{\sf fin}\,$ let
$m = {\rm maxlen}({\rm domC}(f))$. Then $\,f|_{\{0,1\}^m}$  $\in$ 
${\sf pfl}{\cal RM}^{\sf fin}_2$, and 
$\,f \equiv_{\rm fin} f|_{\{0,1\}^m}$\,.
 \ So, ${\sf plep}{\cal RM}^{\sf fin}/\!\!\equiv_{\rm fin}\,$ is isomorphic 
to $\,{\sf pfl}{\cal RM}^{\sf fin}_2/\!\!\equiv_{\rm fin}\,$ by the 
isomorphism $ \ [f]_{\rm fin} \,\longmapsto\, [f|_{\{0,1\}^m}]_{\rm fin}\,$, 
with $f$ ranging over ${\sf plep}{\cal RM}^{\sf fin}_2$; here 
$[f]_{\rm fin}$ denotes the $\,\equiv_{\rm fin}$-class of $f$. 

A similar reasoning applies for {\sf tlep} and {\sf tfl}.
  \ \ \  \ \ \ $\Box$

\medskip

\noindent As a consequence of Lemma \ref{lepVSfl}, we use the notations
${\sf tfl}M_{2,1}$ and ${\sf tlep}M_{2,1}$ interchangeably, and similarly
for ${\sf pfl}M_{2,1}$ and ${\sf plep}M_{2,1}$. 

\bigskip

\noindent We also consider the subgroup

\medskip

 \ \ \  \ \ \ ${\sf lp}M_{k,1}$
$ \ = \ {\sf lp\,totsurinj}{\cal RM}^{\sf fin}/\!\!\equiv_{\rm fin}\,$,

\smallskip

\noindent where 

\medskip

 \ \ \  \ \ \ ${\sf lp}{\cal RM}^{\sf fin}$ $\,=\,$
$\{f \in {\cal RM}^{\sf fin} \,: \ |f(x)| = |x| \ $ for all 
$x \in {\rm Dom}(f)\,\}$.

\medskip

\noindent The latter is called the monoid of {\em length-preserving} 
right-ideal morphisms.

The group ${\sf lp}V$ is the subgroup of length-preserving elements of $V$
(see \cite{BiFact}), and it is the {\em group of units} of the monoids 
${\sf plep}M_{k,1}$ and ${\sf tlep}M_{k,1}$. 

\bigskip

We next study faithful actions of $M_{k,1}$ and its submonoids and subgroups.
It is well known that $G_{k,1}$ acts faithfully on $A^{\omega}$, and 
$M_{k,1}$ acts similarly on $A^{\omega}\,$ (see 
\cite[{\small Subections 1.2-1.3}]{BiLR}).

Here we will show that $M_{k,1}$ and $G_{k,1}$ also act faithfully on
$A^*\,0^{\omega}\,$. As opposed to $A^{\omega}$, the set $A^*\,0^{\omega}\,$
is countable.

\smallskip

We extend right-ideal morphisms of $A^*$ to functions on
$A^{\omega}$, as follows.  For any right-ideal morphism $f$, any 
$x \in {\rm Dom}(f)$, and any $z \in A^{\omega}$, we define

\bigskip

\noindent {\small \sf (a)} \hspace{4.5cm}
$f(x\,) \,=\, f(x) \ z$. 

\medskip

\noindent In particular, we can let $z = 0^{\omega}$, and in that case $f$
is just extended to $A^*\,0{\omega}$.
 \ The extension of $f$ to $A^{\omega}$ or to $A^*\,0^{\omega}$ will also be
called $f$; and usually we will view right-ideal morphisms of $A^*$ as 
functions on $A^* \cup A^{\omega}$.

\begin{lem} \label{y1y2z}
 \ If $y_1, y_2 \in A^*$ are such that
$\,(\forall z \in A^*)[\,y_1 z\,0^{\omega} = y_2 z\,0^{\omega}\,]$, \ then 
$\,y_1 = y_2$.
\end{lem}
{\sc Proof.} Since $y_1 z0^{\omega} = y_2 z0^{\omega}$, the strings $y_1$ 
and $y_2$ are prefix-comparable. Let us assume $y_1 \le_{\rm pref} y_2$ 
(the case $y_2 \le_{\rm pref} y_1$ is similar). Then $y_2 = y_1 v$ for some
$v \in A^*$, so $y_1 z0^{\omega}$  $=$ $y_1 v z0^{\omega}$, which implies 
$z0^{\omega} = v z0^{\omega}$, for all $z \in A^*$. Suppose (for a 
contradiction) that $v \ne \e$. Let $a \in A$ denote the left-most letter in
$v$. Since $z$ is arbitrary, we can pick $z \in A \minus \{a\}$. Then
$z0^{\omega} \ne v z0^{\omega}$ (since they start with different letters), 
so we have a contradiction. Hence $v = \e$, so $y_1 = y_2$.
 \ \ \  $\Box$

\begin{pro} \label{PROPactOnA0omega}\!\!.

\smallskip

\noindent {\small \rm (1)} \ For any two right-ideal morphisms $f, g$ of 
$A^*$ (and their extensions to $A^{\omega}$ and to $A^*\,0^{\omega}$): 

\smallskip

 \ \ \  \ \ \ $f \equiv_{0^{\omega}} g$ \ \ iff \ \ $f|_{A^*\,0^{\omega}}$ 
$\,=\,$ $g|_{A^*\,0^{\omega}}$ \ \ (i.e., $f$ and $g$ have the same action on 
$A^*\,0^{\omega}$).

\medskip

\noindent {\small \rm (2)} \ The monoid $M_{k,1}$ and the Higman-Thompson
group $G_{k,1}$ act faithfully on $A^*0^{\omega}$ by the action
{\small \sf (a)}.
\end{pro}
{\sc Proof.} \ (1) follows from the next two Claims.

\smallskip

\noindent {\sf Claim 1.} If two right-ideal morphisms $f, g$ act the same on
$A^*0^{\omega}$ then $\,f \equiv_{0^{\omega}} g$.

\smallskip

\noindent Proof of Claim 1: \ Since $f$ and $g$ act the same on
$A^*0^{\omega}$, they have the same domain in $A^*0^{\omega}$; i.e.,
${\rm domC}(f)\,A^*\,0^{\omega}$ $=$ ${\rm Dom}(f)\,0^{\omega}$ $=$
${\rm Dom}(g)\,0^{\omega}$ $=$ ${\rm domC}(g)\,A^*\,0^{\omega}$. Hence,
${\rm domC}(f)$  $\equiv_{0^{\omega}}$  ${\rm domC}(g)$.

Moreover, since $f$ and $g$ act the same on $A^*0^{\omega}$ they agree on
${\rm Dom}(f)\,0^{\omega}$ $=$ ${\rm Dom}(g)\,0^{\omega}$ $=$
${\rm Dom}(f)\,0^{\omega}$ $\,\cap\,$  ${\rm Dom}(g)\,0^{\omega}$  $=$
$({\rm Dom}(f)$  $\cap$  ${\rm Dom}(g))\,0^{\omega}$; the latter equality
holds by Lemma \ref{LEM0omegaInters}.  Hence, for all $x \in$
${\rm Dom}(f)$ $\cap$  ${\rm Dom}(g)$ and all $z \in A^*$:
$\,f(x)\,z0^{\omega}$  $=$  $g(x)\,z0^{\omega}$.  Since this holds for all
$z$, we conclude (by Lemma \ref{y1y2z}) that for all $x \in$
${\rm Dom}(f)$ $\cap$  ${\rm Dom}(g):$ \ $f(x) = g(x)$.
So $f$ and $g$ agree on ${\rm Dom}(f)$  $\cap$  ${\rm Dom}(g)$. Thus, 
$f \equiv_{0^{\omega}} g$.
 \ \ \ [This proves Claim 1.]

\medskip

\noindent {\sf Claim 2.} If $f, g \in {\cal RM}^{\sf P}$ satisfy
$f \equiv_{0^{\omega}} g\,$ then $f$ and $g$ act the same on
$A^*0^{\omega}$.

\smallskip

\noindent Proof of Claim 2: Since $f \equiv_{0^{\omega}} g$ we have
${\rm domC}(f)\,A^*\,0^{\omega}$ $=$  ${\rm domC}(g)\,A^*\,0^{\omega}$,
hence $f$ and $g$ have the same domain in $A^*0^{\omega}$.
Moreover, $f \equiv_{0^{\omega}} g$ implies that $f$ and $g$ agree on
${\rm Dom}(f) \cap {\rm Dom}(g)$, by the definition of
$f \equiv_{0^{\omega}} g$.  This implies that $f$ acts the same as $g$ on
${\rm Dom}(f)\,0^{\omega}$ $=$ ${\rm Dom}(g)\,0^{\omega}$.
 \ \ \ [This proves Claim 2.]

\medskip

\noindent (2) Since for finite prefix codes, $\,\equiv_{0^{\omega}}\,$
coincides with $\,\equiv_{\rm fin}\,$, it follows that $M_{k,1}$ acts
faithfully on $A^*0^{\omega}$ by the action {\small \sf (a)}. Hence, any
submonoid of $M_{k,1}$, in particular $G_{k,1}$, also acts faithfully on
$A^*0^{\omega}$ by the action {\small \sf (a)}.
 \ \ \  \ \ \ $\Box$

\bigskip

\noindent {\bf Remark.} We saw already (in the Remark after the definition of 
$\equiv_{A^*0^{\omega}}$) that  $A^*0^{\omega}$ is in 1-to-1 correspondence 
with the set ${\mathbb Q}_k$ of {\em $k$-ary rationals} in $[0,1[\,$, via the 
map $\,x0^{\omega} \in A^*0^{\omega}$ $\longmapsto$ $0.x \in$ 
${\mathbb Q}_k \cap [0,1[\,$.  So the action of $M_{2,1}$ on $A^*0^{\omega}$ 
is also an action on the set ${\mathbb Q}_k \cap [0,1[\,$.

\bigskip

\medskip

\noindent {\bf Summary of Thompson monoids}

\smallskip

\noindent We have defined the following Thompson monoids:

\medskip

\noindent $\bullet$ \ \ $M_{k,1}$, the action monoid of ${\cal RM}^{\sf fin}$
on $A^{\omega}\,$ (and on $A^*0^{\omega}$), with $|A| = k$.

\medskip

\noindent $\bullet$ \ \ ${\sf tot}M_{k,1}$  $\,=\,$
$\{f \in M_{k,1} :\, {\rm Dom}(f)$
is an {\em essential} right ideal of $A^*\, \}$.

\smallskip

 \ \ \ Equivalently, $f$ acting on $A^{\omega}$ (or on $A^*0^{\omega}$)
is a total function. 

 \ \ \ (For $k=2$, this is Thompson's original monoid.)

\medskip

\noindent $\bullet$ \ \ ${\sf sur}M_{k,1}$   $\,=\,$
$\{f \in M_{k,1} :\, {\rm Im}(f)$
is an {\em essential} right ideal of $A^*\,\}$.

\smallskip

 \ \ \ Equivalently, $f$ acting on $A^{\omega}$ (or on $A^*0^{\omega}$)
is a surjective function.

\medskip

\noindent $\bullet$ \ \ ${\sf inj}M_{k,1}$   $\,=\,$
$\{f \in M_{k,1} :\, f$ is injective on $A^{\omega}$\}.

\medskip

\noindent $\bullet$ \ \ ${\sf plep}M_{k,1}$   $\,=\,$
$\{f \in M_{k,1} :\,$
$(\forall x_1, x_2 \in {\rm Dom}(f) \subseteq A^*)[\, |x_1| = |x_2|$
$\Rightarrow$  $|f(x_1)| = |f(x_2)|\,] \,\}$.

\smallskip

 \ \ \ These functions are called {\em (partial) length-equality preserving}.

\medskip

\noindent $\bullet$ \ \ ${\sf tlep}M_{k,1}$  $\,=\,$
$\{f \in {\sf tot}M_{k,1} :\,$
$(\forall x_1, x_2 \in {\rm Dom}(f) \subseteq A^*)[\, |x_1| = |x_2|$
$\Rightarrow$ $|f(x_1)| = |f(x_2)|\,] \,\}$.

\smallskip

 \ \ \ These functions are called {\em total and length-equality preserving}.

\medskip

\section{Finite generation of {\boldmath ${\cal RM}^{\sf fin}_2$} }

The proof that the monoid  ${\cal RM}^{\sf fin}_2$ is finitely generated is 
carried out in several steps, and is based in part on earlier results. In 
\cite{BiThompsMon} it was shown that $M_{2,1}$ is finitely generated. 
Recall that ${\cal RM}^{\sf fin}$ uses the alphabet $\{0,1\}$.

In \cite{BiOntoTh} it was shown that the following submonoid of 
${\cal RM}_k^{\sf fin}$ is finitely generated: 

\bigskip

\hspace{0.75cm} ${\sf tot sur inj}RM^{\sf fin} \ = \ $
$\{f \in {\cal RM}_k^{\sf fin} \, : \, f$ is {\em injective}, \ and both 
${\rm domC}(f)$ and ${\rm imC}(f)$ are 

\hspace{6,4cm} finite {\em maximal} prefix codes\}.

\bigskip

\noindent In \cite{BiOntoTh} this monoid was called ${\sf riAut}(k)$ (for 
``right-ideal automorphisms'').

\subsection{Normal functions}

One of the difficulties in the study of right-ideal morphisms of $A^*$ comes 
from the existence of non-normal right-ideal morphisms, as defined next.

\begin{defn} \label{DEFnormal}
 \ A right-ideal morphism $f$ of $A^*$ is called {\bf normal} \ \ iff
 \ \ $f({\rm domC}(f)) = {\rm imC}(f)$.
\end{defn}
It is easy to prove \cite[{\small Lemma 5.7}]{Equiv} that this is 
equivalent to \ $f^{-1}({\rm imC}(f)) = {\rm domC}(f)$.

\medskip

\noindent For any right-ideal morphism $f$ of $A^*:$ 

\smallskip

 \ \ \  \ \ \ ${\rm imC}(f) \subseteq f({\rm domC}(f))$, \ \ and 
 \ \ $f^{-1}({\rm imC}(f)$ $\subseteq$ ${\rm domC}(f)$,   

\smallskip

\noindent by \cite[{\small Lemmas 5.1 and 5.2}]{Equiv}.
These inclusions could be strict, as illustrated by the following 
example. \ Let $f$ be the right-ideal morphism defined by 
${\rm domC}(f) = \{00, 01\}$, \ ${\rm imC}(f) = \{0\}$, and for all 
$x \in \{0,1\}^*$: $\, f(00 x) = 0 x$, $f(01 x) = 00 x$.  
So, $f(\{00,01\}) = \{0, 00\}$, hence $f$ is non-normal. 

All injective right-ideal morphisms are normal 
\cite[{\small Lemma 5.2}]{Equiv}). So non-normal right-ideal morphisms are 
not encountered in the study of the Thompson group $G_{2,1}$ and its 
subgroups.
 
In \cite[{\small Theorem 4.5B}]{BiThompsMonV3} it was shown that every 
right-ideal morphism in ${\cal RM}^{\sf fin}$ is 
$\, \equiv_{\rm fin}$-equivalent to a normal right-ideal morphism in 
${\cal RM}^{\sf fin}$. So in the study of $M_{2,1}$, non-normal morphisms
can ultimately avoided  (although they show up as intermediate products in
calculations).   

Every function in ${\sf pfl}{\cal RM}^{\sf fin}$ is normal (as is easily 
proved). However, there exist non-normal functions in 
${\sf tlep}{\cal RM}^{\sf fin}$. E.g., for $f$ given by the table 
$\{(0,0), (10,00)\}$ we have $f({\rm domC}(f)) = \{0,00\}$  $\ne$ 
$\{0\} = {\rm imC}(f)$.

\bigskip

For the finite generation of ${\cal RM}^{\sf fin}$ the following is crucial.

\begin{thm} {\bf (injective-normal factorization).} \label{NonNorm2gen}

\smallskip 
 
\noindent For every {\em non-normal} right-ideal morphism $f \in$ 
${\cal RM}^{\sf fin}$ there is an {\em injective} right-ideal morphism 
$j \in {\cal RM}^{\sf fin}$ and a {\em normal} right-ideal morphism 
$\nu \in {\cal RM}^{\sf fin}$ such that

\smallskip

 \ \ \  \ \ \ $f(.) \,=\, \nu \circ j(.)$. 

\smallskip

\noindent Moreover, if $f$ is total then $j$ and $\nu$ are total. 
\end{thm}
{\sc Proof.}  Let $P = {\rm domC}(f)$ and $Q = {\rm imC}(f)$. 
For every $p \in P$ there exists a unique $q_p \in Q$ and a unique 
$u_p \in A^*$ such that $\, f(p) = q_p \, u_p$. 
We define the right-ideal morphisms $j$ and $\nu$ as follows: 

\smallskip

 \ \ \  \ \ \ ${\rm domC}(j) = {\rm domC}(\nu) = P$; 

\smallskip

\noindent and for every $p \in P$:

\smallskip

 \ \ \  \ \ \ $j(p) = p \, u_p$, \ \ and 

\smallskip

 \ \ \  \ \ \ $\nu(p) = q_p$.

\smallskip

\noindent Then $\, \nu \circ j = f$. Indeed, for any $p \in P$ we have
$\, \nu(j(p)) = \nu(p \, u_p) = \nu(p) \ u_p = q_p \, u_p = f(p)$. 

The morphism $j$ is injective, since $p$ is the unique prefix of $p\,u_p$ 
that belongs to the prefix code $P$. We have
$ \ {\rm imC}(j) = \{p \, u_p : p \in P\}$.
This is a prefix code, since if $p_1 \, u_{p_1}$ is a prefix of 
$p_2 \, u_{p_2}$ for $p_1, p_2 \in P$, then $p_1$ and $p_2$ are 
prefix-comparable; hence $p_1 = p_2$ (since $P$ is a prefix code); and then 
$u_{p_1} = u_{p_2}$ (since $p$ uniquely determines $u_p$). 

For the right-ideal morphism $\nu$ we have $ \ \nu(P) = Q = {\rm imC}(\nu)$,
hence $\nu$ is normal.

\smallskip

We saw that ${\rm domC}(j) = {\rm domC}(\nu) = P$. So, if $P$
is a maximal prefix code then $j$ and $\nu$, as constructed above, are total.
 \ \ \  \ \ \ $\Box$

\medskip

\noindent Since injective right-ideal morphisms are normal (by Lemma 5.2 in 
\cite{Equiv}), Theorem \ref{NonNorm2gen} implies:

\begin{pro} \label{NonNormgenbyNorm} 
 \ Every non-normal right-ideal morphism in ${\cal RM}^{\sf fin}$, or in 
${\rm tot}{\cal RM}^{\sf fin}$, is the composite of two {\em normal} 
right-ideal morphisms in ${\cal RM}^{\sf fin}$, respectively in
${\rm tot}{\cal RM}^{\sf fin}$.
  \ \ \  \ \ \ $\Box$
\end{pro}

\subsection{Proof that ${\cal RM}^{\sf fin}_2$ is finitely generated}

The proof has the following outline (and has some overlap with the proof of 
finite generation of $M_{k,1}$ in \cite{BiThompsMonV3}):
 \ (1) We saw that ${\sf totsurinj}{\cal RM}^{\sf fin}$ (a.k.a.\ 
${\sf riAut}(k)$) is finitely generated (see \cite{BiOntoTh}). 
 \ (2) We prove that the submonoid of injective elements of 
${\cal RM}^{\sf fin}_2$ is finitely generated.
 \ (3) We show that the set of all normal elements of ${\cal RM}^{\sf fin}_2$
is finitely generated. 
 \ (4) By Prop.\ \ref{NonNormgenbyNorm} it follows that 
${\cal RM}^{\sf fin}_2$ is finitely generated.

\begin{pro} \label{injRMfinFinGen}
 \ The submonoid $ \ {\sf inj}{\cal RM}^{\sf fin}_2 = \, \{f \in$
${\cal RM}^{\sf fin}_2 : f$ {\rm is injective}\} \ is finitely generated.
\end{pro}
{\sc Proof.} Consider any $f \in {\sf inj}{\cal RM}^{\sf fin}_2$, and let
$P = {\rm domC}(f)$ and $Q = {\rm imC}(f)$. Since injective right-ideal 
morphisms are normal, $f(P) = Q$ and $f^{-1}(Q) = P$. By injectiveness,
$|P| = |Q|$.
 
\medskip

\noindent $\bullet$ Case 1: \ $P$ and $Q$ are both maximal prefix codes.

Then $f \in {\sf totsurinj}{\cal RM}^{\sf fin}_2$, so $f$ is generated by 
any finite generating set of ${\sf totsurinj}{\cal RM}^{\sf fin}_2$. 

\smallskip

\noindent $\bullet$ Case 2: \ $P$ and $Q$ are both non-maximal prefix codes.

Then by \cite[{\small Lemma 3.1}]{BiThompsMon} there exist maximal prefix 
codes $P_M, Q_M \subset \{0,1\}^*$ such that $P \subset P_M$, 
$\,Q \subset Q_M$, and $|P_M| = |Q_M|$. 
Consider the following ``standard'' maximal prefix code
$C \subset 0^*1$ of cardinality $|P_M|$:

\medskip

 \ \ \   $C \ = \ $
$0^{|P_M| - 2} \{0,1\} \ \cup \ \{0^r 1 : \, 0 \le r \le |P_M| - 3\}$.

\smallskip

\noindent Since $P$ is not empty and is a non-maximal prefix code, we 
have $|P_M| \ge 2$. The second set in the above union is empty when 
$|P_M| = 2$. 
It is easy to see that $C$ is a maximal prefix code, and $|C| = |P_M|$.

Since $|C| = |P_M| = |Q_M|$, and these are maximal prefix codes, there exist 
right-ideal morphisms $\,g_1: C \{0,1\}^* \to P_M \{0,1\}^*$, 
$ \ g_2: Q_M \{0,1\}^* \to C \{0,1\}^*$, with $g_2, g_1$ chosen in 
$\,{\sf totsurinj}{\cal RM}^{\sf fin}_2$ in such a way that 
$\, g_2\,f\,g_1 = {\sf id}_D$. Here, $D \subset C$ is chosen to consist
of the first $|P|$ elements of $C$ in dictionary order. More precisely, 

\medskip

 \ \ \   
$D \ = \ 0^{|P_M| - 2} \{0,1\}$ \ $\cup$ 
 \ $\{0^r 1 : \, |P_M| - |P| \le r \le |P_M| - 3 \}$. 

\smallskip

\noindent Since $g_1$ and $g_2$ are injective, it follows from 
$g_2 f g_1(.) = {\sf id}_D \,$ that $\,f\, =$ 
$\,g_2^{-1} \ {\sf id}_D \ g_1^{-1}(.)$.
 
It is straightforward to prove that 
$\, D \equiv_{\rm fin} 0^{|P_M| - |P|} \{0,1\}$,
so $\, {\sf id}_D \equiv_{\rm fin} {\sf id}_{0^{|P_M| - |P|} \{0,1\}}\,$ 
(the latter is the identity restricted to the right-ideal
$\,0^{|P_M| - |P|} \{0,1\} \, \{0,1\}^*$).
It was proved in \cite[{\small Lemma 3.3}]{BiThompsMon} that for all 
$j > 0$: $ \ {\sf id}_{0^j \{0,1\}}$ is generated by 
$\, \{(0 \to 00), \ (00 \to 0), \ {\sf id}_{0\{0,1\}}\}$. (In 
\cite[{\small Lemma 3.3}]{BiThompsMon}, this was proved in $M_{2,1}$, but
the proof does not use $\equiv_{\rm fin}$, just composition, and the
restrictions implied by composition.)

Finally, $\, 0^{|P_M| - |P|} \{0,1\} \ \equiv_{\rm fin} \ D$ 
 \ $\subset$ \ $C$ \ $=$ \ ${\rm imC}(g_1^{-1}) = {\rm domC}(g_2^{-1})$. 
Since we have $\,D\,\{0,1\}^*$ $\,\subseteq\,$ 
$0^{|P_M| - |P|} \{0,1\}\,\{0,1\}^*$,
$\,{\sf id}_D$ is obtained from ${\sf id}_{0^{|P_M| - |P|} \{0,1\}}$ by
restriction. So can just use composition (with $g_1^{-1}$ on the right and
$g_2^{-1}$ on the left) to implement this $\equiv_{\rm fin}$-equivalence; 
hence

\medskip

 \ \ \ $g_2^{-1} \ {\sf id}_{0^{|P_M| - |P|} \{0,1\}} \ g_1^{-1}$
 \ $=$ \ $g_2^{-1} \ {\sf id}_D \ g_1^{-1}$ \ ($= f$).

\smallskip

\noindent Since $g_2, g_1 \in {\sf totsurinj}{\cal RM}^{\sf fin}_2$, it 
follows that $f$ is generated by a finite generating set of 
${\sf totsurinj}{\cal RM}^{\sf fin}_2$ together with the finite set
$\, \{(0 \to 00), \ (00 \to 0), \ {\sf id}_{0\{0,1\}}, \ $
${\sf id}_{\{0,1\}}\}$ $\subset$ ${\sf inj}{\cal RM}^{\sf fin}$.

\medskip

\noindent $\bullet$ Case 3: \ $P$ is a maximal prefix code, and $Q$ is 
non-maximal.

Then by \cite[{\small Lemma 3.1}]{BiThompsMon} there exists a maximal prefix
code $Q_M \subset \{0,1\}^*$ such that $Q \subset Q_M$ and 
$|P| = |Q| < |Q_M|$. Consider the following ``standard'' maximal prefix code
$C' \subset 0^*1$ of cardinality $|Q_M|$:

\medskip

 \ \ \   $C'$ $ \ = \ $
$0^{|Q_M| - 2} \{0,1\} \ \cup \ \{0^r 1 : \, 0 \le r \le |Q_M| - 3 \}$.

\medskip

\noindent As in Case 2, we have $|Q_M| \ge 2$, and $C'$ is a maximal 
prefix code with $|C| = |Q_M|$, where $C$ is the ``standard'' maximal prefix
code of cardinality $|P|$:

\medskip

 \ \ \   $C$ $ \ = \ $
$0^{|P| - 2} \{0,1\} \ \cup \ \{0^r 1 : \, |P| - 3 \ge r \ge 0\}$.

\smallskip

\noindent If $|P| = 2$ then $P = \{0,1\}$ and this formula yields $C =P$.
If $|P| = 1$ then $P = \{\varepsilon\}$, and we choose $C = P$.

As in Case 2, there exist $g_2, g_1 \in {\sf totsurinj}{\cal RM}^{\sf fin}_2$
with $g_1: C \{0,1\}^* \to P \{0,1\}^*$, and 
$g_2: Q_M \{0,1\}^* \to C' \{0,1\}^*$, where $g_2, g_1$ are chosen so that 
$g_2\,f\,g_1: C \{0,1\}^* \to C' \{0,1\}^*$ preserves the 
dictionary order, and maps $C$ bijectively onto the first $|P|$ elements of 
$C'$ in the dictionary order. 
Let $S$ be this set of the first $|P|$ elements of $C'$, so

\medskip

 \ \ \   
$S \ = \ \ 0^{|Q_M| - 2} \{0,1\} \ \ \cup \ \ $
$\{0^r 1 : \, |Q_M| - |P| \le r \le |Q_M| - 3 \}$.
 
\medskip

\noindent Then $\,P \equiv_{\rm fin} \{\varepsilon\}$, 
 \ $S  \,\equiv_{\rm fin}\, 0^{|Q_M| - |P| -1} \{0,1\}$ $ \ \cup \ $
  $\{0^r 1 : \, |Q_M| - |P| - 2 \ge r \ge 0\}$, 
 \ both $P$ and $S$ are final segments of  $0^* 1$ in the dictionary order, 
and $\, g_2fg_1(.) \,$ preserves the dictionary order.  Hence,  
$\, g_2\,f\,g_1(.) \,$ $\equiv_{\rm fin}$ $\,(\e \to 0^{|Q_M| - |P|})$. 
Moreover, 

\smallskip

 \ \ \ $S\,\{0,1\}^*$  $ \ \subseteq \ $  
$\big(0^{|Q_M| - |P| -1} \{0,1\}$ $ \ \cup \ $
  $\{0^r 1 : \, |Q_M| - |P| - 2 \ge r \ge 0\}\big)\,\{0,1\}^*$,

\smallskip

\noindent so $\equiv_{\rm fin}$ can be implemented by composition (with 
$g_1^{-1}$ on the right and $g_2^{-1}$ on the left).  Hence
 \ ($f =$) $g_2^{-1}\,g_2\,f\,g_1\,g_1^{-1}$ $ \ = \ $ 
$g_2^{-1} \,(\e \to 0^{|Q_M| - |P|}) \, g_1^{-1}$.

By \cite[{\small Lemma 3.3}]{BiThompsMon}, for any $j \ge 1$, 
$ \ (\e \to 0^j)\,$ is generated by 
$\,\{(\e \to 0), \ (0 \to 00), \ (00 \to 0)\}$. 
So $f$ is generated by a finite generating set of 
$\,{\sf totsurinj}{\cal RM}^{\sf fin}_2$ and the finite set 
$\,\{(\e \to 0), \ (0 \to 00), \ (00 \to 0)\}$ $\subset$ 
${\sf inj}{\cal RM}^{\sf fin}_2$.

\medskip

\noindent $\bullet$ Case 4: \ $P$ is a non-maximal prefix code, and $Q$ 
is maximal.

This is proved by applying Case 3 to $f^{-1}$, and by using the inverses of
the generators of Case 3.
 \ \ \  \ \ \ $\Box$

\medskip

\noindent Prop.\ \ref{injRMfinFinGen} was proved for the binary alphabet 
$A = \{0,1\}$, but the same proof could be adapted to $A =$
$\{0,1,\ldots,k\!-\!1\}$ for any $k \ge 2$.

\begin{lem} \label{NormFinGen}
 \ The normal elements of ${\cal RM}^{\sf fin}_2$ are generated by a finite 
set of normal elements.
\end{lem}
{\sc Proof.} Let $f \in {\cal RM}^{\sf fin}_2$ be a normal element with
$P = {\rm domC}(f)$, $Q = {\rm imC}(f)$; by normality, $Q = f(P)$ and 
$P = f^{-1}(Q)$. The proof is similar in outline to the proof of finite 
generation of $M_{k,1}\,$ \cite[{\small Theorem 3.5}]{BiThompsMon}. 

If $\,|P| = |Q|\,$ then $f$ is injective; and for 
${\sf inj}{\cal RM}^{\sf fin}_2$ we already have a finite generating set, by 
Prop.\ \ref{injRMfinFinGen}).
For the remainder of the proof we have $\,|P| > |Q|$. 

\medskip

\noindent {\sf Claim.}  {\it Let $f \in {\cal RM}^{\sf fin}$ be a 
{\em normal} right-ideal morphism, with $P = {\rm domC}(f)$, 
$\,Q = {\rm imC}(f)$  $=$  $f(Q)$, such that $\,|P| > |Q|$. Then $f$ is 
equal to a composite of $\,\le |P| - |Q|\,$ normal right-ideal morphisms 
$\,\varphi_i \in$ ${\cal RM}^{\sf fin}_2\,$ with 
$\,P_i = {\rm domC}(\varphi_i)\,$ and 
$\,Q_i = {\rm imC}(\varphi_i) = \varphi_i(P_i)$, such that 
$ \ |P_i| - |Q_i| \le 1$.
} 

\medskip

\noindent Proof of the Claim: We use induction on $|P| - |Q|$. In the base
case, $|P| - |Q| = 1$, there is nothing to prove.  Assume now 
$|P| - |Q| \ge 2$; then one, or both, of the next two cases apply.

\smallskip

Case 1: There exist $x_1, x_2, x_3 \in P$ with $f(x_1) = f(x_2) = f(x_3)$
($= y_1 \in Q$), with $x_1, x_2, x_3$ all different.

Then $f(.) = \psi_2 \psi_1(.)$, where $\psi_2, \psi_1$ are defined as 
follows: ${\rm domC}(\psi_1) = P$, $\,{\rm imC}(\psi_1)$ $=$ 
$P \minus \{x_1\}$, and $\,\psi_1(x_1) = \psi_1(x_2) = x_2$;
$\psi_1(x) = x$ for all $x \in P \minus \{x_1\}$. And 
${\rm domC}(\psi_2) = P \minus \{x_1\}$ and ${\rm imC}(\psi_2) = Q$, and 
$\psi_2(x) = f(x)$ for all $x \in P \minus \{x_1\}$; in particular,
$\psi_2(x_2) = \psi_2(x_3) = y_1$. Then $\psi_2, \psi_1$ are normal. For 
the sizes of their domain code and image code we have: 
$|P| - |P \minus \{x_1\}|$ $=1$ for $\psi_1$. And $\,|P \minus \{x_1\}| - |Q|$
 $=$ $|P| - |Q| -1 < |P| - |Q|$, so the inductive step applies to $\psi_2$.

\smallskip

Case 2: There exist $x_1, x_2 \in P$ with $f(x_1) = f(x_2)$ ($= y_1 \in Q$), 
and there exist $x_3, x_4 \in P$ with $f(x_3) = f(x_4)$ ($= y_2 \in Q$),
with $y_1 \ne y_2$, and $x_1, x_2, x_3, x_4$ all different.

Then $f(.) = \psi_2 \psi_1(.)$, where $\psi_2, \psi_1$ are the same as in
Case 1.  

\smallskip

In either case, $f$ is factored into one normal function with difference 1,
and one with difference $|P| - |Q| -1 < |P| - |Q|$. So, by induction, $f$ is
factored into $\,\le |P| - |Q|\,$ normal functions with difference 1. 
 \ \ \ [This proves the Claim.]

\medskip

By the Claim it suffices to show that every normal $f \in$ 
${\cal RM}^{\sf fin}_2$ with $\,|P| - |Q| = 1\,$ is generated by a finite 
set of normal elements. Let us denote $|P| = n$, 
$\,P = \{p_1, \ldots, p_n\}$, $\,Q = \{q_1, \ldots, q_{n-1}\}$, such that 
$f(p_i) = q_i$ for $i = 1, \ldots, n\!-\!1$, and $f(p_n) = q_{n-1}$.
Except for $\,f(p_{n-1}) = f(p_n) = q_{n-1}$, the indexing is arbitrary. 
We define a ``standard'' prefix code $C \subset 0^*1$ of size $|C| = |P|:$

\smallskip

 \ \ \  \ \ \    
$C \ = \ 0^{|P| -2} \{0,1\} \ \cup \ \{0^r 1 : \, 0 \le r \le |P| -3 \}$.

\smallskip

\noindent This is the same code $C$ as in Prop.\ \ref{injRMfinFinGen}. 
When $|P| = 2$ the formula yields $C = \{0,1\}$. Since $|P| = |Q| + 1$, we
have $|P| \ge 2$. We enumerate $C$ in increasing dictionary order as 
$(c_1, \ldots, c_n)$; so $c_1 = 0^{n - 1}$, and 
$c_i = 0^{n - i} 1$ for $2 \le i \le n$; in particular, $c_n = 1$.

We factor $f$ as $f(.) = \psi_3 \psi_3 \psi_1(.)$, where 
$\psi_3, \psi_3, \psi_1$ are defined as follows on their domain and image
codes: 

\smallskip

 \ \ \ $\psi_1: P \to C$ \ is the bijection $p_i \mapsto c_i$ for 
$1 \le i \le n$.

\smallskip

 \ \ \ $\psi_2: C \to C \minus \{c_n\}$ \ is given by $c_i \mapsto c_i$ 
for $1 \le i \le n\!-\!1$, and $c_n \mapsto c_{n-1}$.

\smallskip

 \ \ \ $\psi_3:  C \minus \{c_n\} \to Q$ \ is the bijection 
$c_i \mapsto q_i$ for $1 \le i \le n\!-\!1$.

\smallskip

\noindent Clearly, $f(.) = \psi_3 \psi_3 \psi_1(.)$. Since $\psi_1$ and
$\psi_3$ are injective they are finitely generated, by Prop. 
\ref{injRMfinFinGen}.

We still want to factor $\,\psi_2: C \to C \minus \{c_n\}\,$ over a fixed
finite generating set. 

If $n=2$ then $c_n = 1$, and $c_1 = c_{n-1} = 0$.
Then $ \ \psi_2 \ = \ $
$\left( \hspace{-.08in} \begin{array}{l|l}
0 & 1  \\
0 & 0
\end{array} \hspace{-.08in} \right)$.

\smallskip

If $n>2$ then $\psi_2(c_i) = c_i$ for $1 \le i \le n\!-\!1$, and
$\psi_2(c_n) = c_{n-1}$. And $c_1 = 0^{n-1}$, and $c_i$  $=$ 
$0^{n-i} 1\,$ for $2 \le i \le n$.  Therefore,

\medskip

$\psi_2 \ = \ $
$\left( \hspace{-.08in} \begin{array}{l|l|l|l|l|l|l|l|l}
0^{n-1} & 0^{n-2}1 & 0^{n-3}1 & \ \ \ldots \ \ & 0^{n-i}1 & \ \ldots \ & 001 & 01 & 1\\  
0^{n-1} & 0^{n-2}1 & 0^{n-3}1 & \ \ \ldots \ \ & 0^{n-i}1 & \ \ldots \ & 001 & 01 & 01
\end{array} \hspace{-.08in} \right)$.

\smallskip

\noindent Since $\psi_2$ is the identity function for its first $n\!-\!1$ 
elements in the dictionary order, we have 

\smallskip

$\, \psi_2 \ \equiv_{\rm fin} \ $
$\left( \hspace{-.08in} \begin{array}{l|l|l}
00 & 01 & 1  \\
00 & 01 & 01
\end{array} \hspace{-.08in} \right)$
$ \ \equiv_{\rm fin} \ $
$\left( \hspace{-.08in} \begin{array}{l|l}
0 & 1  \\
0 & 01
\end{array}  \hspace{-.08in} \right)$ $ \ = \ $ $\{(0,0),\,(1,01)\}$. 

\smallskip

\noindent  Hence, since $\,\psi_2 \subseteq \{(0,0),\,(1,01)\}$, 

\smallskip

 \ ($f =$) 
$\psi_3 \circ \psi_2 \circ \psi_1$ $=$  
$\psi_3$ $\circ $
$\left( \hspace{-.08in} \begin{array}{l|l}
0 & 1  \\
0 & 01
\end{array}  \hspace{-.08in} \right)$ $\circ$ $\psi_1$.

\medskip

\noindent The middle function is not normal, but it can be factored into 
two normal functions: 

\medskip

$\left( \hspace{-.08in} \begin{array}{l|l}
0 & 1  \\
0 & 01
\end{array}  \hspace{-.08in} \right)$ $=$
$\left( \hspace{-.08in} \begin{array}{l|l}
0 & 1  \\
0 & 0
\end{array} \hspace{-.08in} \right)$ $\circ$
$\left( \hspace{-.08in} \begin{array}{l|l}
0 & 1  \\
0 & 11
\end{array} \hspace{-.08in} \right)$.

\smallskip

\noindent Hence, $f$ is generated by a finite set of generator of 
$\, {\sf inj}{\cal RM}^{\sf fin} \,$, together with the two new generators

\medskip

\noindent
$\left( \hspace{-.08in} \begin{array}{l|l}
0 & 1  \\
0 & 0
\end{array} \hspace{-.08in} \right)$  \ and \     
$\left( \hspace{-.08in} \begin{array}{l|l}
0 & 1  \\
0 & 11
\end{array}  \hspace{-.08in} \right)$.  
 \ \ \  \ \ \  \ \ \ $\Box$

\bigskip

\noindent From Lemma \ref{NormFinGen} and Prop.\ \ref{NonNormgenbyNorm} we 
now conclude:

\begin{thm} \label{RMfinFinGen}\!\!{\bf .} 

\smallskip

\noindent  The monoid ${\cal RM}^{\sf fin}_2$ is finitely generated.  

\noindent Moreover, ${\cal RM}^{\sf fin}_2$ has a finite generating set 
consisting of normal elements.  \ \ \  \ \ \  \ \ \ $\Box$
\end{thm}

\smallskip

\section{{\boldmath ${\sf tlep}{\cal RM}^{\sf fin}_2$, 
$\,{\sf tlep}M_{2,1}$, $\,{\sf plep}{\cal RM}^{\sf fin}_2$, and
$\,{\sf plep}M_{2,1}$, \\  
 \ are not finitely generated} }

The question whether ${\sf tlep}M_{2,1}$ is finitely generated was raised
in \cite[{\small Section 2}]{Bi1wPerm}, and here we answer it negatively. 
First, we need a special version of a well-known fact.

\begin{lem} \label{injComposite}
 \ For all $g, f \in {\cal RM}^{\sf fin}$: \ If $g \circ f(.)$ is total 
and injective, then $f$ is total and injective, and $g$ is defined
on all of $ \ {\rm imC}(f) \,A^{\omega}$.  
\end{lem}
{\sc Proof.} 
1. $f$ is total: By contraposition, suppose $f$ is not total, i.e., there 
exists $x \in A^*$ such that $\,xA^* \cap {\rm Dom}(f) = \0$. Then 
$g(f(xA^*))$ $=$ $\0$, so $\,xA^* \cap {\rm Dom}(gf(.)) = \0$, hence $gf(.)$
is not total.

\noindent
2. $g$ is defined on ${\rm imC}(f) \,A^{\omega}$: If $g \circ f(.)$ is total,
then $g$ is defined on every $z \in {\rm imC}(f) \,A^{\omega}$. Indeed, if 
that were not the case then there would exist $f(x) \in {\rm imC}(f)$ (for 
some $x \in A^*$) such that $g(f(x)\,A^*) = \0$. Then $gf(xA^*) = \0$, hence
$gf(.)$ would not be total. 

\noindent
3. $f$ is injective: By contraposition, suppose $f$ is not injective, i.e., 
there exist $a, b \in {\rm Dom}(f)$ with $a \neq b$ and $f(a) = f(b)$. Then 
$az \neq bz$ and $f(a) \, z = f(b) \, z$ for all $z \in A^*$.  Hence 
$g \circ f(a z) = g \circ f(b z)$, and this is defined for long enough $z$, 
since $g$ is defined on all of $\,{\rm imC}(f) \,A^{\omega}$.  Hence 
$g \circ f(.)$ is not injective.
 \ \ \ $\Box$

\begin{thm} \label{lepMnotFinGen}\!\!\!.

\smallskip

\noindent The monoids $\,{\sf tlep}{\cal RM}^{\sf fin}_2$, 
$\,{\sf tfl}{\cal RM}^{\sf fin}_2$, and $\,{\sf tlep}M_{2,1}$ 
$(= {\sf tfl}M_{2,1})$, as well as 
$\,{\sf plep}{\cal RM}^{\sf fin}_2$,
$\,{\sf pfl}{\cal RM}^{\sf fin}_2$, and $\,{\sf plep}M_{2,1}$
$(= {\sf pfl}M_{2,1})$, are \textbf{\textit{not finitely 
generated}}.
\end{thm}
{\sc Proof.} (1) Case of {\sf tfl} monoids: 

\noindent Since ${\sf tlep}M_{2,1}$ ($= {\sf tfl}M_{2,1}$) is a 
homomorphic image of ${\sf tfl}{\cal RM}^{\sf fin}_2$ and of 
${\sf tlep}{\cal RM}^{\sf fin}_2$, it is enough to prove that 
${\sf tlep}M_{2,1}$ is not finitely generated. Recall that 
${\sf tlep}M_{2,1}$ acts faithfully on $\{0,1\}^{\omega}$. 
 \ For any $\,N > 0\,$ let 

\smallskip

 \ \ \  \ \ \ $\Gamma_{\! N} \,=\,$ 
$\{f \in {\cal RM}^{\sf fin} \,:\, (\exists n,m \le N)[$
${\rm domC}(f) = \{0,1\}^m$ \ and
 \ ${\rm imC}(f) \subseteq \{0,1\}^n\,]\,\}$.

\smallskip

\noindent If ${\sf tlep}M_{2,1}$ were finitely generated, then 
${\sf tlep}M_{2,1}$ would be generated by $\Gamma_{\! N}$ for some $N$. Let 
us show that ${\sf not}_{N+1}\,$ is not generated by $\Gamma_{\! N}$. 
For all $\,x_1\ldots x_N x_{N+1} \in \{0,1\}^{N+1}\,$ and $u \in $
$\{0,1\}^{\omega}$,

\smallskip

 \ \ \  \ \ \  ${\sf not}_{N+1}(x_1\ldots x_N \, x_{N+1}\,u)$ 
$ \ = \ $  $x_1\ldots x_N\, \ov{x}_{N+1}\, u$,

\smallskip

\noindent where $\ov{x}_{N+1}$ is the negation of the bit $x_{N+1}$.

\smallskip

Suppose, for a contradiction, that ${\sf not}_{N+1}(.) =$ 
$\gamma_T\,\ldots\,\gamma_1(.)\,$ for some $T \ge 1$, and 
$\gamma_i \in \Gamma_{\! N}$ for $i \in [1, T]$. For every $i \in [1, T]$
let ${\rm domC}(\gamma_i) = \{0,1\}^{m_i}$ and 
${\rm imC}(\gamma_i) \subseteq \{0,1\}^{n_i}$; since $\gamma_i \in$
$\Gamma_N$ we have $m_i, n_i \le N$.

Since ${\sf not}_{N+1}$ is total and injective it follows from Lemma 
\ref{injComposite} that $\,\gamma_t\,\ldots\,\gamma_1(.)\,$ is total and 
injective for all $t \in [1,T]$. Totalness and injectiveness imply that for
every $t$ there exists $s \in \{0,1\}^{N+1}$ such that 
$\,|\gamma_t\,\ldots\,\gamma_1(s)| \ge N\!+\!1$.
And since $\,\gamma_t\,\ldots\,\gamma_1(.)\,$ is {\sf tfl}, we have 
$\,|\gamma_t\,\ldots\,\gamma_1(s)| = |\gamma_t\,\ldots\,\gamma_1(x)|\,$ for
all $x$ such that $|x| = |s|$. 
 \ Hence for all $t \in [1, T]$ and all $x \in \{0,1\}^{N+1}:$
 
\smallskip

\noindent {\small $(\star)$}
 \ \ \  \ \ \  \ \ \ $|\gamma_t\,\ldots\,\gamma_1(x)| \ \ge \ N\!+\!1\,$. 

\smallskip

\noindent Since $N\!+\!1 > N \ge m_i\,$ for all $i \in [1,T]$, every 
$\gamma_t\,\ldots\,\gamma_1(x)$ is defined (for all $t \in [1, T]$ and 
$x \in \{0,1\}^{N+1}$).
And the application of $\gamma_{t+1}$ to $\,\gamma_t\,\ldots\,\gamma_1(x)\,$
changes only a prefix of length $\, \le N\,$ (since $m_{t+1}, n_{t+1}$ 
$\le N$). 

Consider now any $z = p^{(0)} z_{N+1} v \in \{0,1\}^{\omega}$, with 
$p^{(0)} \in \{0,1\}^N$ and $v \in \{0,1\}^{\omega}$. Then for all 
$t \in [1,T]:$ \ $\gamma_t\,\ldots\,\gamma_1(p^{(0)} z_{N+1} v)\,$ is 
defined, and  
$ \ \gamma_T\,\ldots\,\gamma_t\,\ldots\,\gamma_1(p^{(0)} z_{N+1} v)$
$=$ $p^{(0)}\,\ov{z}_{N+1}\,v$.  \ Moreover there are strings $\,p^{(1)}$,
$\,\ldots\,$, $p^{(t)}$, $\,\ldots\,$, $p^{(T)}$ $\in \{0,1\}^*$ such that

\medskip

 \ \ \  \ \ \ $\gamma_1(p^{(0)} z_{N+1} v) \,=\, p^{(1)} z_{N+1} v$ \ with 
      \ $|p^{(1)}| \ge N$

\medskip

 \ \ \  \ \ \ $\gamma_2 \gamma_1(p^{(1)} z_{N+1} v) \,=\, p^{(2)} z_{N+1} v$ 
 \ with \ $|p^{(2)}| \ge N$,

\medskip
 
 \ \ \  \ \ \  \ \ $\vdots$

\medskip

 \ \ \  \ \ \ $\gamma_t \gamma_{t-1}\,\ldots\,$
                $\gamma_2 \gamma_1(p^{(t-1)} z_{N+1} v)$
        $\,=\,$ $p^{(t)} z_{N+1} v$ \ with \ $|p^{(t)}| \ge N$.

\medskip

 \ \ \  \ \ \  \ \ $\vdots$

\medskip

 \ \ \  \ \ \ $\gamma_T\,\ldots\,\gamma_2 \gamma_1(p^{(T-1)} z_{N+1} v)$  
        $\,=\,$ $p^{(T)} z_{N+1} v$ \ with \ $|p^{(T)}| \ge N$.

\medskip

\noindent The relations $|p^{(t)}| \ge N$ hold by relation 
{\small $(\star)$}. And the suffix $\,z_{N+1} v\,$ is not modified because 
every $\gamma_t$ only changes a prefix of length $\,\le N$ of $p^{(t-1)}$, 
while $p^{(t-1)}$ itself has length $\,\ge N$.  

\smallskip

\noindent So $\,\gamma_T\,\ldots\,\gamma_1(.)$ $=$ ${\sf not}_{N+1}\,$ does
not change $z_{N+1}$, which contradicts the definition of ${\sf not}_{N+1}$.

\medskip

\noindent (2) \ Case of {\sf pfl} monoids:

\noindent As in (1) it is sufficient to prove the result for 
${\sf plep}M_{2,1}$. \ For any $N > 0$ let

\smallskip

 \ \ \  \ \ \ $\Gamma_{\! N} \,=\,$
$\{f \in {\cal RM}^{\sf fin} \,:\, (\exists n,m \le N)[$
${\rm domC}(f) \subseteq \{0,1\}^m$ \ and 
 \ ${\rm imC}(f) \subseteq \{0,1\}^n\,]\,\}$.

\smallskip

\noindent If ${\sf plep}M_{2,1}$ were finitely generated, then it would be
generated by $\Gamma_{\! N}$ for some $N$. Let us show that 
${\sf not}_{N+1}\,$ is not generated by $\Gamma_{\! N}$.
 \ (The difference with part (1) is that ``${\rm domC}(f) =$''  is replaced
by ``${\rm domC}(f) \subseteq$''.)

Suppose, for a contradiction, that ${\sf not}_{N+1}(.) =$
$\gamma_T\,\ldots\,\gamma_1(.)\,$ for some $T \ge 1$, and
$\gamma_i \in \Gamma_{\! N}$ for $i \in [1, T]$. For every $i \in [1, T]$,
let ${\rm domC}(\gamma_i) \subseteq \{0,1\}^{m_i}$ and
${\rm imC}(\gamma_i) \subseteq \{0,1\}^{n_i}$; since $\gamma_i \in$
$\Gamma_N$ we have $m_i, n_i \le N$.

Since ${\sf not}_{N+1}$ is total and injective it follows from Lemma
\ref{injComposite} that $\,\gamma_t\,\ldots\,\gamma_1(.)\,$ is total and
injective for all $t \in [1,T]$, and that $\gamma_{t+1}$ is defined on all
of $\,{\rm imC}(\gamma_t\,\ldots\,\gamma_1) \ \{0,1\}^{\omega}$. 
Totalness and injectiveness imply that for every $t$ there exists 
$s \in \{0,1\}^{N+1}$ such that
$\,|\gamma_t\,\ldots\,\gamma_1(s)| \ge N\!+\!1$.
And since $\,\gamma_t\,\ldots\,\gamma_1(.)\,$ is {\sf pfl}, we have
$\,|\gamma_t\,\ldots\,\gamma_1(s)| = |\gamma_t\,\ldots\,\gamma_1(x)|\,$ for
all $x \in {\rm Dom}(\gamma_t\,\ldots\,\gamma_1)$ such that $|x| = |s|$.
 \ Hence, since $\gamma_t\,\ldots\,\gamma_1$ is total we have for all 
$t \in [1, T]$ and all $x \in \{0,1\}^{N+1}:$

\smallskip

\noindent {\small $(\star)$}
 \ \ \  \ \ \  \ \ \ $|\gamma_t\,\ldots\,\gamma_1(x)| \ \ge \ N\!+\!1\,$.

\smallskip

\noindent Since $N\!+\!1 > N \ge m_i\,$ for all $i \in [1,T]$, every
$\gamma_t\,\ldots\,\gamma_1(x)$ is defined (for all $t \in [1, T]$ and
$x \in \{0,1\}^{N+1}$).
And the application of $\gamma_{t+1}$ to $\,\gamma_t\,\ldots\,\gamma_1(x)\,$
changes only a prefix of length $\, \le N\,$ (since $m_{t+1}, n_{t+1}$
$\le N$).

Consider now any $z = p^{(0)} z_{N+1} v \in \{0,1\}^{\omega}$, with
$p^{(0)} \in \{0,1\}^N$ and $v \in \{0,1\}^{\omega}$. Then all $t \in [1,T]:$
$ \ \gamma_t\,\ldots\,\gamma_1(p^{(0)} z_{N+1} v)\,$ is defined, and
$\,\gamma_T\,\ldots\,\gamma_t\,\ldots\,\gamma_1(p^{(0)} z_{N+1} v)$
$=$ $p^{(0)}\,\ov{z}_{N+1}\,v$.  \ Moreover there are strings $\,p^{(1)}$,
$\,\ldots\,$, $p^{(t)}$, $\,\ldots\,$, $p^{(T)}$ $\in \{0,1\}^*$ such that
the relations of part (1) hold. \ The rest of the proof is exactly as in 
part (1). 
 \ \ \  \ \ \ $\Box$ 

\bigskip

\noindent {\large \bf Consequences for circuits}

\smallskip

We saw in Theorem \ref{circ2lepM} that the set of input-output functions of 
acyclic boolean circuits, extended to right-ideal morphisms, is precisely the
monoid ${\sf tfl}{\cal RM}^{\sf fin}$, and in Theorem \ref{circ2PartLepM} 
that the set of input-output functions of acyclic partial circuits, extended
to right-ideal morphisms, is precisely the monoid 
${\sf pfl}{\cal RM}^{\sf fin}$.
Then Theorem \ref{lepMnotFinGen} implies:

\begin{thm} \label{PROcircNotFinGen}
 \ The monoids of input-output functions of acyclic boolean circuits, or 
partial circuits, extended to right-ideal morphisms of $\{0,1\}^*$, are 
{\em not finitely generated} under function-composition.
 \ \ \  \ \ \ $\Box$
\end{thm}
However, ${\sf tfl}{\cal RM}^{\sf fin}$ and ${\sf pfl}{\cal RM}^{\sf fin}$ 
are generated under composition by a set of the form $\Gamma \cup \tau$, 
where $\Gamma$ is finite, and $\tau$ is the set of bit-position 
transpositions (or just the set of transpositions $\tau_{i,i+1}$ of 
adjacent positions).

\section{Some algebraic properties of the monoids \\   
{\boldmath $M_{2,1}$, 
$\,{\sf tot}M_{2,1}$, $\,{\sf sur}M_{2,1}$, and $\,{\sf plep}M_{k,1}$} }

Here we give a few algebraic properties of Higman-Thompson monoids.

\subsection{Regularity}

\begin{pro} \label{RMfinRegular}
 \ The monoid ${\cal RM}^{\sf fin}$ is regular, i.e., for every
$f \in {\cal RM}^{\sf fin}$ there exists $f' \in {\cal RM}^{\sf fin}$ such
that $\,f f' f = f$. Moreover, $f'$ can be chosen to be injective with, in
addition, ${\rm domC}(f') = {\rm imC}(f)$.
\end{pro}
{\sc Proof.} For any $f \in {\cal RM}^{\sf fin}$ and any
$y \in {\rm imC}(f)$, let us choose an element $x'_y \in f^{-1}(y)$.
It is a fact that $f^{-1}({\rm imC}(f)) \subseteq {\rm domC}(f)\, $ (see
\cite[{\small Lemma 5.5}]{Equiv}), hence $x'_y \in {\rm domC}(f)$.
We define $f' \in {\cal RM}^{\sf fin}$ by $\, f'(y u) = x'_y u \,$ for every
$y \in {\rm imC}(f)$ and $u \in \{0,1\}^*$; so,
${\rm domC}(f') = {\rm imC}(f)$. Hence $f' \in {\cal RM}^{\sf fin}$, since
it is a right-ideal morphism, and ${\rm imC}(f)$ is finite. Since for
different elements $y$ the sets $f^{-1}(y)$ are disjoint, $f'$ is injective.
(Note also that $f^{-1}(y)$ is also finite for every $y \in {\rm imC}(f)$.)

Then for any $x v \in {\rm Dom}(f)$ with $x \in {\rm domC}(f)$:
$f f' f(x v) = f f'(yu v)$, where $f(x) = yu$ for some
$y \in {\rm imC}(f)$, $u \in \{0,1\}^*$.
Hence, $f f'(yu v) = f(x'_y u v)$, where $x'_y = f'(y) \in {\rm domC}(f)$;
then, $f(x'_y u v) = y uv = f(x) \, v$. Thus, $f f' f(x v) = f(xv)$, so
$f'$ is an inverse of $f$.
 \ \ \ $\Box$

\bigskip

A similar argument show that all the pre-Thompson and Thompson monoids 
are regular.

\medskip

For a generating set $\Gamma \cup \tau$ of ${\cal RM}^{\sf fin}$, with
$\Gamma$ finite, the problem of finding an inverse is {\sf coNP}-hard.
Moreover, on input $w, w' \in (\Gamma \cup \tau)^*$, deciding whether $w'$
represents an inverse of $w$ (i.e., $E_w \,E_{w'}\,E_w$  $=$  $E_w$) is
{\sf coNP}-complete (see \cite{BiCoNP, BinG}).

\subsection{Congruence-simplicity}

We know that $M_{k,1}$ and ${\sf plep}M_{k,1}$ are congruence-simple
\cite{BiThompsMonV3, BiThompsMon}. We show here that Thompson's original
monoid ${\sf tot}M_{2,1}$ is also congruence-simple, but that 
${\sf tlep}M_{k,1}$ is not congruence-simple.

\medskip

The proofs below use elements of $M_{k,1}$ denoted by $(u \to v)$, where
$u, v \in A^*$. The function $(u \to v)$ has a singleton table $\{(u,v)\}$;
so ${\rm Dom}((u \to v)) = u A^{\omega}$, $\,{\rm Im}((u \to v)) = $
$v A^{\omega}$; and for all $uz \in uA^{\omega}$: $\,(u \to v)(uz) = vz$.
Since $\{\e\}$ is a maximal prefix code, $(\e \to u) \in {\sf tot}M_{k,1}$.

More generally, for any $v \in A^*$ and any finite prefix code 
$C \subseteq A^*$ we define the function $(C \to v)$ by
${\rm domC}((C \to v)) = C$, $\,{\rm imC}((C \to v)) = \{v\}$, and for
all $c \in C$ and $z \in A^{\omega}:$ $ \ (C \to v)(cz) = vz$. 
If $C$ is a maximal prefix code then $\,(C \to v) \in {\sf tot}M_{k,1}$.

\begin{thm} \label{THMcongsimplityMplepM}
 \ For all $k \ge 2$, the monoids $\,M_{k,1}\,$ and $\,{\sf plep}M_{k,1}\,$ 
are {\em $0$-${\cal J}$-simple} and {\em congruence-simple}.
\end{thm}
{\sc Proof.}  For $M_{k,1}$ this was proved in 
\cite[{\small \rm Prop.\ 2.2 and Thm.\ 2.3}]{BiThompsMonV3}, and in 
\cite{BiThompsMon}. 

For ${\sf plep}M_{k,1}$, the proofs of 
\cite[{\small \rm Prop.\ 2.2, Thm.\ 2.3}]{BiThompsMonV3}
work without any change, since  ${\sf plep}M_{k,1}$ has a {\bf 0}, and all
the multipliers used in the proofs are of the form $(u \to v)$ for various
$u, v \in A^*$.  And $(u \to v)$ obviously belongs to ${\sf plep}M_{k,1}$.
 \ \ \  \ \ \ $\Box$

\begin{thm} \label{THMcongsimplityTotM}
 \ For all $k \ge 2$, the monoid $\,{\sf tot}M_{k,1}\,$ is
{\em congruence-simple}.
\end{thm}
{\sc Proof.} Let $\equiv$ be any congruence on ${\sf tot}M_{k,1}$ that is 
not the equality relation. We will prove that $\equiv$ has just one 
congruence class.

\medskip

\noindent {\sf Claim 1.} There exist strings $y_0, y_1 \in A^*$ that are 
not prefix-comparable, such that $\,(\e \to y_0) \equiv (\e \to y_1)$.

\smallskip

\noindent Proof of Claim 1: Since $\equiv$ is not equality, there exist
$\psi, \varphi \in {\sf tot}M_{k,1}$ such that $\psi \equiv \varphi$ and
$\psi \ne \varphi$. By restriction we can represent $\psi$ and $\varphi$
by right-ideal morphisms that are normal (Def.\ \ref{DEFnormal}), with 
${\rm imC}(\psi)$ $\cup$ ${\rm imC}(\varphi) \subseteq A^n\,$ for some
$n > 0$.  Hence there are $x_0 \in A^*$ and $y_0, y_1 \in A^n$ such that
$\,\psi(x_0) = y_0 \ne y_1 = \varphi(x_0)$. Since $|y_0| = |y_1|$,
$\,y_0 \ne y_1$ implies that $y_0, y_1$ are not prefix-comparable.
Then $\psi \equiv \varphi$ implies

\smallskip

 \ \ \ $(\e \to  y_0) \ = \ \psi \circ (\e \to x_0)(.)$  
$ \ \equiv \ $  $\varphi \circ (\e \to x_0)(.) \ = \ (\e \to  y_1)$.

\smallskip

\noindent [This proves Claim 1.]

\medskip

\noindent {\sf Claim 2.} For all $x, y \in A^*:$ 
 \ $(\e \to x) \,\equiv\, (\e \to y)$.

\smallskip

\noindent Proof of Claim 2: By Claim 1 there are $u, v \in A^*$ that are 
not prefix-comparable, such that $\,(\e \to v) \equiv (\e \to v)$.
Since $u$ and $v$ are not prefix-comparable, there exists a maximal finite 
prefix code $C \subseteq A^*$ such that $\{u,v\} \subseteq C$; \ in fact,
in the proof of Claim 1 we have $u, v \subseteq A^n$.  For any
$x,y \in A^*$, let $f \in$ ${\sf tot}M_{k,1}$ be the function defined by
${\rm domC}(f) = C$, $\,f(C) = \{x,y\}$, with $\,f(u) = x$, $\,f(v) = y$, 
and $\,f(z) = x$ for all $z \in C \minus \{u,v\}$. (The set $\{x,y\}$ need
not be a prefix code.)  
Then $\,(\e \to u) \equiv (\e \to v)\,$ implies 
 
\smallskip

 \ \ \ $(\e \to x) = f \circ (\e \to u)(.)$  $\,\equiv\,$  
$f \circ (\e \to v)(.) = (\e \to y)$. 

\smallskip

\noindent [This proves Claim 2.]

\medskip

\noindent {\sf Claim 3.} For all $x \in A^*:$ 
 \ $(\e \to x) \,\equiv\, {\mathbb 1}$.

\smallskip

\noindent Proof of Claim 3: This follows from Claim 2 by letting $y = \e$.

\noindent [This proves Claim 3.]

\medskip

\noindent {\sf Claim 4.} For every finite maximal prefix code 
$P \subset A^*:$  \ $(P \to \e) \,\equiv\, {\mathbb 1}$.

\smallskip

\noindent Proof of Claim 4: Here ${\mathbb 1}$ denotes the identity function
on $A^* \cup A^{\omega}$. By Claim 3, $\,(\e \to p) \equiv {\mathbb 1}\,$
for every $p \in P$. This implies (for any particular $p \in P$): 

\smallskip

 \ \ \ ${\mathbb 1} = (\e \to \e) = (P \to \e) \circ (\e \to p)(.)$  
$\,\equiv\,$  $(P \to \e) \circ {\mathbb 1}(.) = (P \to \e)$.

\smallskip

\noindent [This proves Claim 4.]

\medskip

\noindent {\sf Claim 5.} For all $f \in{\sf tot}M_{2,1}:$ 
$ \ {\mathbb 1} \equiv f$.

\smallskip

\noindent Proof of Claim 5:  By restriction, any $f \in {\sf tot}M_{2,1}$ 
can be represented by a normal right-ideal morphism in
${\sf tot}{\cal RM}^{\sf fin}\,$ (that we also call $f$).  
Let $P = {\rm domC}(f)$, and let $Q$ be a finite maximal prefix code such 
that $f(P) = {\rm imC}(f) \subseteq Q$.
Since $P$ and $Q$ are finite maximal prefix codes, Claim 4 implies that 
 \ ${\mathbb 1} \equiv (P \to \e) \equiv (Q \to \e)$.

Then $\,(Q \to \e) \equiv {\mathbb 1}\,$ implies
 \ $(P \to \e) = (Q \to \e) \circ f(.) \equiv {\mathbb 1} \circ f(.) = f$.
Now $\,{\mathbb 1} \equiv (P \to \e)\,$ implies $\,{\mathbb 1} \equiv f$.

\noindent [This proves Claim 5 and the Theorem].
 \ \ \  \ \ \  \ \ \ $\Box$

\bigskip

\noindent Combined with the results of de Witt and Elliott
\cite{deWittElliott} and Brin \cite{Brin_PC}, Theorem 
\ref{THMcongsimplityTotM} yields (for all $k \ge 2$):

\begin{cor} \label{CORtotM}\!\!.

\smallskip

\noindent
The monoid $\,{\sf tot}M_{k,1}$ is an infinite, finitely presented,
congruence-simple monoid that has the 
Higman-Thompson group $G_{k,1}$ as its group of units.
 \ \ \  \ \ \ $\Box$
\end{cor}

\smallskip

\noindent 
We will now prove that ${\sf tlep}M_{k,1}$ is not congruence-simple.

\begin{defn} \label{DEFdelta} {\bf (the input-output length difference 
function $\delta$).}

\smallskip

\noindent The input-output length difference function $\delta$ is defined by

\medskip

 \ \ \  $\delta: \ f \in {\sf tlep}{\cal RM}^{\sf fin}$  
$ \ \longmapsto \ \delta(f) = |f(x)| - |x| \,\in\, {\mathbb Z}$, 

\medskip

\noindent for any $x \in {\rm domC}(f)$ .
\end{defn}

\begin{pro} \label{PROPiolDiffHomLEP}\!\!\!.

\smallskip

\noindent {\small \rm (0)} \ For all $f \in {\rm tlep}{\cal RM}^{\sf fin}$,
$\,\delta(f)$ is well-defined (i.e., it just depends on $f$, not on $x$).

\smallskip

 \ Moreover, $\,\delta(f) = |f(x)| - |x| \ $ for all $x \in {\rm Dom}(f)\,$ 
(not just for $x \in {\rm domC}(f)$).

\medskip

\noindent {\small \rm (1)} \ For all $f_1, f_2 \in$
${\rm tlep}{\cal RM}^{\sf fin}:$ \ If $\,f_1 \equiv_{\rm fin} f_2\,$ then
 \ $\delta(f_1) = \delta(f_2)$.

\smallskip

 \ So $\delta(f)$ is well-defined for all $\,f \in {\rm tlep}M_{2,1}$.

\medskip

 \ If $\,{\rm domC}(f) = \{0,1\}^m\,$ for some $m \ge 0$, then
$\,f({\rm domC}(f)) = {\rm imC}(f)$ $\subseteq$ $\{0,1\}^n$, where 

 \ $\,n = m + \delta(f)$. \ Hence $\delta$ is surjective. 

\medskip

\noindent {\small \rm (2)} \ For all $f_1, f_2 \in$
${\rm tlep}{\cal RM}^{\sf fin}:$
 \ \ $\delta(f_2 \circ f_1) = \delta(f_2) + \delta(f_1)$.

\smallskip

 \ Hence, $\,\delta: {\rm tlep}M_{2,1} \twoheadrightarrow {\mathbb Z} \ $
and $ \ \delta: {\rm tlep}{\cal RM}^{\sf fin} \twoheadrightarrow$
${\mathbb Z} \ $ are surjective monoid homomorphisms.

\medskip

\noindent {\small \rm (3)} \ Let $\Gamma$ be a finite set such
that $\,\Gamma \cup \tau$ generates $\,{\sf tfl}{\cal RM}^{\sf fin}\,$ or
$\,{\rm tlep}M_{2,1}$. Let 
$\,E: (\Gamma \cup \tau)^* \to {\sf tfl}{\cal RM}^{\sf fin}\,$ be the 
generator evaluation function, i.e., for every $w = w_k \ldots w_1 \in$ 
$(\Gamma \cup \tau)^*$ such that $w_i \in \Gamma \cup \tau$ we have 
$\,E_w(.) = w_k \circ \,\ldots\, \circ w_1(.)$.
 \ Then

\smallskip

 \ \ \  \ \ \ $\delta(E_w) \ = \ \sum\,(\delta(w_i) : 1 \le i \le k, \ $
$w_i \in \Gamma)$.
\end{pro}
{\sc Proof.} (0) This follows from the fact that $f$ is {\sf lep} 
(length-equality preserving).
 \ (1) If $f$ and $g$ differ by a one-step restriction (i.e., 
${\rm domC}(g) = ({\rm domC}(f) - \{x\}) \cup xA$, for some $x \in$
${\rm domC}(f)$), then $\delta(f) = \delta(g)$. Then we use the fact that
$\equiv_{fin}$ is the same as $\equiv_1^*\,$ (by Prop.\ \ref{equi1setdiff}).
 \ (2) is straightforward, and (3) follows from (2), using the fact that
$\delta(\tau_{i,j}) = 0$. 
 \ \ \  \ \ \ $\Box$

\bigskip

For any $w \in (\Gamma \cup \tau)^*$, the integer $\delta(E_w)\,$ (in
sign-magnitude binary representation) is easily computed (by using (3)
above)

\bigskip

\noindent The existence of the homomorphism $\delta(.)$ (by (2) above)
implies the following:

\begin{cor} \label{CORtlepCngrSim}\!\!.

\smallskip

\noindent  For all $k \ge 2$, the monoids $\,{\rm tlep}M_{k,1}$, 
$\,{\rm tlep}{\cal RM}^{\sf fin}$, and $\,{\sf tfl}{\cal RM}^{\sf fin}\,$ 
are {\em not congruence-simple}.
\end{cor}

\noindent The ``kernel'' of $\delta: {\rm tlep}M_{2,1} \to {\mathbb Z}$ is 
the monoid 

\smallskip

 \ \ \  \ \ \  $\delta^{-1}(0)$ $\,=\,$ 
${\sf tlp}M_{k,1}$ $\,=\,$  $\{f \in {\sf tlep}M_{k,1} :$ 
$(\forall x \in {\rm Dom}(f))\,[\,|f(x)| = |x|\,]\,\}$,

\smallskip

\noindent where {\sf lp} stands for {\em length-preserving}.

\bigskip

As a (partial) function, $\delta$ can also be defined for 
${\sf plep}{\cal RM}^{\sf fin}$, $\,{\sf pfl}{\cal RM}^{\sf fin}$, and 
${\sf plep}M_{k,1}$.
However, $\delta(\theta)$ is undefined (where $\theta$ is the empty 
function), so on these monoids $\delta$ is not a total function. The
following is proved in the same way as Prop.\ \ref{PROPiolDiffHomLEP}.

\begin{pro} \label{PROPiolDiffpartial}\!\!\!.

\smallskip

\noindent {\small \rm (0)} \ For all $f \in$ 
${\rm plep}{\cal RM}^{\sf fin} \minus \{\theta\}$, $\,\delta(f)$ is 
well-defined (i.e., it just depends on $f$, not on $x$).

\medskip

\noindent {\small \rm (1)} \ For all $f_1, f_2 \in$
${\rm plep}{\cal RM}^{\sf fin} \minus \{\theta\}:$ \ If 
$\,f_1 \equiv_{\rm fin} f_2\,$ then \ $\delta(f_1) = \delta(f_2)$.

\smallskip

 \ So $\delta(f)$ is well-defined for all $\,f \in$
${\rm plep}M_{2,1} \minus \{\theta\}$.

\medskip

\noindent {\small \rm (2)} \ For all $f_1, f_2 \in$
${\rm plep}{\cal RM}^{\sf fin} \minus \{\theta\}\,$ 
such that $f_2 \circ f_1(.) \ne \theta$,
  \ \ $\delta(f_2 \circ f_1) = \delta(f_2) + \delta(f_1)$.

\medskip

\noindent {\small \rm (3)} \ Let $\Gamma$ be a finite set such
that $\,\Gamma \cup \tau$ generates $\,{\sf pfl}{\cal RM}^{\sf fin}\,$ or
$\,{\rm plep}M_{2,1}$. Let
$\,E: (\Gamma \cup \tau)^* \to {\sf pfl}{\cal RM}^{\sf fin}\,$ be the
generator evaluation function, i.e., for every $w = w_k \ldots w_1 \in$
$(\Gamma \cup \tau)^*$ such that $w_i \in \Gamma \cup \tau$ we have
$\,E_w(.) = w_k \circ \,\ldots\, \circ w_1(.)$. Then we have, if
 \ $E_w \ne \theta:$ 

\smallskip

 \ \ \  \ \ \ $\delta(E_w) \ = \ \sum\,(\delta(w_i) : 1 \le i \le k, \ $
$w_i \in \Gamma)$.

\noindent $\Box$
\end{pro}
But $\,\delta$:$\,{\sf plep}M_{k,1} \to {\mathbb Z}\,$ cannot be extended to
a semigroup homomorphism, since ${\sf plep}M_{k,1}$ is congruence-simple 
(by \ref{THMcongsimplityMplepM}).

\subsection{The monoid {\boldmath $\,{\sf sur}M_{k,1}$} }

\smallskip 

\begin{lem} \label{LEMsurMidempot}
 \ The monoid $\,{\sf sur}M_{k,1}$ has only one idempotent, namely the 
identity element.
\end{lem}
{\sc Proof.} If $f = f^2$ then for all $x \in {\rm Dom}(f)$  $\subseteq$ 
$A^{\omega}$:
$\,y = f(x) = f\,f(x) = f(y)$.  Since $f$ is surjective, $y$ ranges over all of
$A^{\omega}$. Hence, $y = f(y)$ for all $y \in A^{\omega}$.  
 \ \ \ $\Box$

\begin{pro} \label{PROPtotMNOTINsur}
 \ For all $h, k \ge 2$, the monoids $\,{\sf tot}M_{h,1}$ and $M_{h,1}$ 
are {\em not} embeddable in ${\sf sur}M_{k,1}$. 
\end{pro}
{\sc Proof.} By Lemma \ref{LEMsurMidempot}, the only idempotent is the 
identity function.  But ${\sf tot}M_{h,1}$ and $M_{h,1}$ contain infinitely 
many idempotents, and an embedding preserves idempotents and distinctness.      
 \ \ \ $\Box$ 

\bigskip

\noindent It follows that $\,{\sf sur}M_{k,1}$ is not isomorphic to
$\,{\sf tot}M_{h,1}$ nor to $M_{h,1}$.

\smallskip

\section{Non-embeddability of {\boldmath $M_{h,1}$ into 
$\,{\sf tot}M_{k,1}$, \\       
and of $\,{\sf plep}M_{h,1}$ into $\,{\sf tlep}M_{k,1}$} } 

\smallskip

\begin{thm} \label{THMMnotintotM} 
 \ For all $h, k \ge 2$, the monoids $M_{h,1}$ and ${\sf plep}M_{h,1}$ are 
{\em not} embeddable into $\,{\sf tot}M_{k,1}$.
\end{thm}
I.e., there is no injective monoid homomorphism
$\,M_{h,1} \hookrightarrow {\sf tot}M_{k,1}$.

\smallskip

\noindent As a consequence, the monoid ${\sf plep}M_{h,1}$ is {\em not} 
embeddable into ${\sf tlep}M_{k,1}\,$ ($\, \subseteq {\sf tot}M_{k,1}$).

\bigskip

The Theorem follows from a series of Lemmas that distinguish $M_{h,1}$ 
from the submonoids of ${\sf tot}M_{k,1}$.

\begin{lem} \label{Lemma1} \!\!\!.

\smallskip

\noindent {\small \rm (1)} For every $f \in {\sf tot}M_{k,1}$ and for 
every infinite set $X \subseteq \{0,1\}^{\omega}$, \ $f(X)$ is infinite.

\medskip

\noindent {\small \rm (2)} For every $f \in {\sf tot}M_{k,1}$ and for 
every {\em finite} set $H \subseteq A^{\omega}$ such that 
$\,|H| > |{\rm domC}(f)|:$

\medskip

\hspace{2cm} $|f(H)| \,>\, |H|/|{\rm domC}(f)|$.
\end{lem}
{\sc Proof.} (1) Since $P = {\rm domC}(f)$ is a finite maximal prefix code,
we have: $\,PA^{\omega} = A^{\omega}$.  Hence, $X \cap PA^{\omega} = X$  
is infinite.
Therefore (since $P$ is finite), there exists $p \in P$ such that
$X\,\cap\, p\,A^{\omega}$  is infinite, hence
$\,X\,\cap\,p\,A^{\omega} \,=\, p\,S$, where $S$ is an infinite subset of
$A^{\omega}$.
Then $f(X)$ contains the infinite set $\, f(p S) = f(p) \ S$.      

\smallskip

\noindent (2) Let $P = {\rm domC}(f)$; this is a finite maximal prefix code,
so $\,P\,A^{\omega} = A^{\omega}$.
There exists $p \in P$ such that  $\,|H \cap p\,A^{\omega}| \ge |H|/|P|$.
Indeed, otherwise we would have $\,|H \cap p\,A^{\omega}| < |H|/|P|\,$ for 
all $p \in P$; this would imply $\,(|H| =) \ |H \cap P\,A^{\omega}|$  
$\,\le\,$ $\sum_{p \in P} |H \cap p\,A^{\omega}|$  $\,<\,$  
$\sum_{p \in P} |H|/|P|$ $=$  $|P| \, (|H|/|P|) = |H|$; this would imply 
$\,|H| < |H|$. 

For $p \in P$ such that  $\,|H \cap p\,A^{\omega}| \ge |H|/|P|$, \ let 

\smallskip

 \ \ \   \ \ \ $H \cap p\,A^{\omega} \,=\,$ 
$p \ \{s_i \ : \ i = 1, \ \ldots \ , \ |H|/|P|, \ |H|/|P|+1, \ \ldots \ \}$,

\smallskip

\noindent where all $s_i\,$ $(\in A^{\omega})$  are different.
Then $f(H)$ has the subset  

\smallskip

 \ \ \   \ \ \ $f(p) \ \{s_i \ : \ $ 
$i = 1, \ \ldots \ , \ |H|/|P|, \ |H|/|P| + 1, \ \ldots \ \}$, 

\smallskip

\noindent so $ \ |f(H)| > |H|/|P|$.         
 \ \ \  \ \ \ $\Box$

\begin{lem} \label{Lemma2}
 \ Let $g$ be any {\em non-torsion} element of $\,{\sf tot}M_{k,1}\,$
(i.e., $\{g^i : i \in \omega\}$ is infinite). Then for all $N > 0$ there 
exists $u \in A^{\omega}$ such that \ $|\{g^i(u) : i \ge 0\}| > N$.
\end{lem}
{\sc Proof.}  By contraposition we assume that there exists $N > 0$ such 
that for all $u \in A^{\omega}:$ $ \ |\{g^i(u) : i \ge 0\}| \le N$.
Since $\,|\{g^i(u) : i \ge 0\}|\,$ is finite, there exist $n, m\,$ 
(depending on $u$) such that  $g^n(u) = g^{n+m}(u)$, and 
$\,0 \le n < n\!+\!m \le N$.
Hence, letting $M = N!$ we have for all $u$: $ \ g^M(u) = g^{2M}(u)$.
So, $g$ is a torsion element.      
 \ \ \ $\Box$

\begin{lem} \label{Lemma3}
 \ If $S$ is a subsemigroup of $\,{\sf tot}M_{k,1}$ that contains a 
non-torsion element, then $S$ has no left-zero.  
\end{lem}
{\sc Proof.} Let $g \in S$ be a non-torsion element. Let us assume, for a
contradiction, that $S$ has a left-zero $z\,$ (i.e., 
$z\,s(.) = z(.)\,$ for all $s \in S$). 
Then $\,z \, g^i(.) = z(.)\,$ for all $i \ge 0$.
Let $P = {\rm domC}(z)$, and let $\,N > 2\,|P|$.
By Lemma \ref{Lemma2}, there exists $u \in A^{\omega}$, such that 
$\,|\{g^i(u) : i \ge 0\}| > N$.
By Lemma \ref{Lemma1}, $\,|z(\{g^i(u) : i \ge 0\})|  >  N/|P|$; and $N$ can
be chosen so that $N/|P| \ge 2$.
On the other hand, $\,z g^i(.) = z(.)$, so
$ \ |z(\{g^i(u) : i \ge 0\})| = |\{z(u)\}| = 1$.
This leads to the false statement $1 \ge 2$.       
 \ \ \ $\Box$

\medskip

\noindent As a consequence, ${\sf tot}M_{k,1}$ itself has no left-zero 
(and hence no zero). However, some subsemigroups have left-zeros.

\bigskip

\noindent {\bf Proof of Theorem \ref{THMMnotintotM}:} 

\smallskip

The monoids ${\sf plep}M_{h,1}$ and $M_{h,1}$ contain a zero (namely the
empty map), and non-torsion elements (e.g., the element represented by 
the right-ideal morphism $(\e \to 0)$).

Assume, for a contradiction, that $M_{h,1}$ (or ${\sf plep}M_{h,1}$) is 
embeddable in ${\sf tot}M_{k,1}$.  An embedding preserves non-torsion 
elements, and the zero in the image submonoid. So the embedded copy of 
$M_{h,1}$ (or ${\sf plep}M_{h,1}$) in ${\sf tot}M_{k,1}$ is a subsemigroup 
with a left zero and a non-torsion element. But this contradicts Lemma 
\ref{Lemma3}.   
\hspace{1,5cm} {\boldmath $\Box$}

\bigskip

\noindent Theorem \ref{THMMnotintotM} for $k \ge 3$ and $h = 2$ can also be
derived in an interesting way from Theorem \ref{THMMnotintotM} for $k = 2$. 
First, a Lemma and a Proposition that are of independent interest are proved:

\begin{lem} \label{MaxPrefCode}\!\!.

\noindent For any finite maximal prefix code $C \subseteq \{0,1\}^*$ such 
that $|C| \ge 2$ we have: \ $\{0,1\}^{\omega} \,=\, C^{\omega}$.
\end{lem}
{\sc Proof.} Obviously, $C^{\omega} \subseteq \{0,1\}^{\omega}$.
 \ Conversely, let us show that any $z = (z_i : i \in \omega)\,$ (with 
$z_i \in$ $\{0,1\}$) is factored in a unique way over $C$. 
Let $\ell = {\rm maxlen}(C)$.
Since $C$ is maximal, $z$ has a unique prefix $c_1 \in C$, of length 
$|c_1| = l_1$ with $1 \le l_1 \le \ell$; so 
$\,z = c_1 z_{l_1} \, \ldots \, $. 
Suppose now, by induction, that a prefix of length $\,\ge n$ of $z$ has been
factored over $C$; i.e., $\,z = c_1\,\ldots\,c_n\,z_{l_n} \, \ldots \, $,
where $c_1,\,\ldots\,,c_n \in C$ with $|c_1\,\ldots\,c_n| = l_n \ge n$. 
Then $\,z_{l_n} \, \ldots \, z_{l_n+\ell}$ has a prefix in $C$; let us
denote this element of $C$ by $c_{n+1}$, with $|c_1\,\ldots\,c_n c_{n+1}| =$
$l_{n+1}$.
Then $\,z = c_1\,\ldots\,c_n c_{n+1}\,z_{l_{n+1}} \, \ldots \, $.
 
Moreover, every finite prefix of $z$ can be factored in at most one way 
over $C$, since $C$ is a prefix code.
 \ \ \  \ \ \ $\Box$

\medskip

\noindent However, $C^* \ne \{0,1\}^*\,$ unless $C = \{0,1\}$.

\begin{pro} \label{PROPtotMktotM2}
 \ For every $k \ge 3$, $ \ M_{k,1}$ is embeddable into $M_{2,1}$,
and $\,{\sf tot}M_{k,1}$ is embeddable into $\,{\sf tot}M_{2,1}$.
\end{pro}
{\sc Proof.} This is proved in the same way as $G_{k,1}$ is embedded into 
$G_{2,1}$, using Higman's coding method. Let $A =$
$\{0,1,\,\ldots\,,k\!-\!1\}$, let $C \subseteq \{0,1\}^*$ be any maximal
prefix code of cardinality $|C| = k$, and let $\,\eta:\, i \in A$
$\longmapsto$  $\eta(i) \in C\,$ be any bijection. 
For any $x = x_1 \,\ldots\, x_n \in A^*$ with $x_i \in$ $A$, 
consider the coding $\eta(x) = $ $\eta(x_1) \,\ldots\, \eta(x_n)$  
$\in C^* \subseteq \{0,1\}^*$; then $\eta:A^* \to \{0,1\}^*$ is an injective 
monoid morphism.
For $x = x_1 \,\ldots\, x_n \,\ldots \ \in A^{\,\omega}$ with $x_i \in A$
we define $\eta(x) = $ $\eta(x_1) \,\ldots\, \eta(x_n)$
$\,\ldots \ $  $\in C^{\omega}$; since $C$ is a maximal prefix code of 
$\{0,1\}^*$ we have $ \, C^{\omega} = \{0,1\}^{\omega}\,$ (by Lemma 
\ref{MaxPrefCode}). 

For any $f \in {\sf tot}M_{k,1}$ let $P = {\rm domC}(f)$. This is a finite
prefix code of $A^*$, which is maximal if $f \in {\sf tot}M_{k,1}$.  
Then $\eta(P)$ is a finite prefix code of $\{0,1\}^*$, which is maximal if
$P$ is maximal.  We encode $f$ by $f^{\eta} \in$ ${\sf tot}M_{2,1}$, 
defined by 

\smallskip

 \ \ \  \ \ \ $f^{\eta}(\eta(x)) = \eta(f(x))$ 

\smallskip

\noindent for every $x \in A^* \cup A^{\,\omega}$. 
 \ Equivalently, for every $z \in C^* \cup C^{\omega}$,
$ \ f^{\eta}(z) = \eta \circ f \circ \eta^{-1}(z)$;  recall that 
$\,\eta: A^* \cup A^{\omega} \to C^* \cup C^{\omega}\,$ is total and 
injective. 

Then $f^{\eta} \in M_{2,1}$ (and $f^{\eta} \in {\sf tot}M_{2,1}$ if 
$f \in {\sf tot}M_{k,1}$), with ${\sf domC}(f^{\eta}) = \eta(P)$. 
And $f \mapsto f^{\eta}$ is a monoid morphism $\,M_{k,1} \to M_{2,1}\,$ 
and $\,{\sf tot}M_{k,1} \to {\sf tot}M_{2,1}$, since 
$\,(g f)^{\eta} = \eta g f \eta^{-1}$  $=$ 
$\eta g \eta^{-1} \eta f \eta^{-1} = g^{\eta} f^{\eta}$.
And $f \mapsto f^{\eta}$ is injective since $\eta$ and $\eta^{-1}$ are 
injective
 \ \ \ $\Box$

\bigskip

\noindent {\sc Proof} of Theorem \ref{THMMnotintotM} for $k \ge 3$ and 
$h = 2$, derived from Theorem \ref{THMMnotintotM} for $k = 2$:

Assume for a contradiction that $M_{2,1}$ is embeddable in 
${\sf tot}M_{k,1}$. By Prop.\ \ref{PROPtotMktotM2}, ${\sf tot}M_{k,1}$ is 
embeddable into ${\sf tot}M_{2,1}$. Now,
$M_{2,1} \hookrightarrow {\sf tot}M_{k,1}$
$\hookrightarrow {\sf tot}M_{2,1}$, so $M_{2,1} \hookrightarrow$ 
${\sf tot}M_{2,1}$. This contradicts the Theorem for $k = 2$.
 \ \ \  \ \ \ $\Box$

\smallskip

\section{Completion operations for $M_{2,1}$, $\,{\sf pfl}M_{2,1}$,
 $\,{\cal RM}^{\sf fin}$, and $\,{\sf pfl}{\cal RM}^{\sf fin}$ }

The goal of a completion operation of a partial function $f$ is to find a 
{\em total} function $\ov{f}$, such that $f$ can be recovered from $\ov{f}$
and ${\rm Dom}(f)$. 
A completion operation $\,\ov{(.)}: f \to \ov{f}\,$ may be a function 
(i.e., only one $\ov{f}$ is found for each $f$); we call this a {\em 
deterministic} completion.  
A completion operation could also be a relation (if $f$ is completed in 
several ways); we call this a {\em nondeterministic} completion.

Theorem \ref{THMMnotintotM} is a non-completability result. 
This section gives several forms of completability.

\smallskip 

We work with the monoids of right-ideal morphisms ${\cal RM}^{\sf fin}$ and
${\sf pfl}{\cal RM}^{\sf fin}$; the corresponding results for $M_{2,1}$
and ${\sf pfl}M_{2,1}$ then hold too. 
Below, for any $h \ge 2$, $\,{\cal RM}_h^{\sf fin}$ denotes the monoid 
${\cal RM}^{\sf fin}$ of right-ideal morphisms $A^* \to A^*$ for an alphabet
$A$ of size $|A| = h$.

The following concept is used in all of our completion constructions.

\begin{defn} \label{DEFcomplprefcode} {\bf (complementary prefix code).}
For a prefix code $P \subseteq A^*$, a {\em complementary prefix code} of
$P$ is a prefix code $Q \subseteq A^*$ such that 
{\small \rm (1)} $\,P \cup Q\,$ is a maximal prefix code, and 
{\small \rm (2)} $\, P A^* \cap Q A^* = \0$.
\end{defn}
Every finite prefix code has finite complementary prefix codes. See 
\cite[{\small Def.\ 5.2}]{BiCoNP}, 
\cite[{\small Def.\ 3.29, Lemma 3.30}]{BiLR}, and also
\cite[{\small Prop.\ 5.2}]{BiEvalProb}.

\subsection{Deterministic completion}

\begin{defn} \label{DEFcompletion}
 \ A {\em deterministic completion} operation of $\,{\cal RM}_h^{\sf fin}\,$
in $\,{\sf tot}{\cal RM}_k^{\sf fin}\,$ (and similarly of $M_{h,1}$ in 
${\sf tot} M_{k,1}$) is a total function 

\smallskip

 \ \ \  \ \ \ $\ov{(.)}: \, f \in {\cal RM}_h^{\sf fin}$  $\,\longmapsto\,$
$\ov{f} \in {\sf tot}{\cal RM}_k^{\sf fin}$, 

\smallskip

\noindent together with a function 

\smallskip

 \ \ \  \ \ \ $\varrho:$
$\{\,\ov{f}: f \in {\cal RM}_h^{\sf fin}\} \,\subseteq$
${\sf tot}{\cal RM}_k^{\sf fin}$  
$ \ \longmapsto \ {\cal RM}_h^{\sf fin}$, 

\smallskip

\noindent such that for all $f \in {\cal RM}_h^{\sf fin}$: 

\smallskip

 \ \ \  \ \ \ $f \,\subseteq\, \varrho(\,\ov{f}\,)$.

\medskip

\noindent In other words, for all $x \in$ ${\rm Dom}(f)$: 
 \ \ $f(x) = \big[\varrho(\,\ov{f}\,)\big](x)$.
\end{defn}

A completion is called {\em homomorphic} iff for all 
$\,f_2, f_1 \in {\cal RM}_h^{\sf fin}$:
 \ $\ov{f_2 \cdot f_1}(.) \,=\, \ov{f}_2 \cdot \ov{f}_1(.)$.
 \ By Theorem \ref{THMMnotintotM}, there exists no injective homomorphic 
completion of $M_{h,1}$.

\medskip

\noindent We give two examples of deterministic completions.

\bigskip

\noindent {\bf 1. A classical completion operation for 
${\sf pfl}{\cal RM}_2^{\sf fin}$ and ${\cal RM}_2^{\sf fin}$ }

\smallskip

\noindent For $f \in {\sf pfl}{\cal RM}_2^{\sf fin}$ with 
${\rm domC}(f)$ $\subseteq$ $\{0,1\}^m\,$ and $\,{\rm imC}(f) \subseteq$
$\{0,1\}^n$, we can define a completion $\tilde{f} \in$ 
${\sf tfl}{\cal RM}_2^{\sf fin}$ as follows:

\smallskip

 \ \ \ ${\rm domC}(\tilde{f}) \,=\, \{0,1\}^m$, 
 \ \ ${\rm imC}(\tilde{f}) \,=\, {\rm imC}(f) \cup \{0^n\}$ 
  $\,\subseteq\, \{0,1\}^n$, 

\smallskip

\noindent and for all $x \in \{0,1\}^m:$

\smallskip

 \ \ \ $\tilde{f}(x) = f(x)$ \ \ if \ \ $x \in {\rm domC}(f)$, 
 
\smallskip

 \ \ \ $\tilde{f}(x) = 0^n$ \ \  \ \ if \ \ $x \in$ 
        $\{0,1\}^m \minus {\rm domC}(f)$.

\bigskip

\noindent More generally, for $f \in {\cal RM}^{\sf fin}$ we can define a
completion $\tilde{f} \in$ ${\sf tot}{\cal RM}^{\sf fin}$ as follows:
 
\smallskip

 \ \ \ ${\rm domC}(\tilde{f}) \,=\, {\rm domC}(f) \,\cup\, Q$, 

\smallskip

\noindent where $Q$ is a finite complementary prefix code of 
${\rm domC}(f)$ in $A^*$.  \ For all $x \in P \cup Q$,

\smallskip

 \ \ \ $\tilde{f}(x) = f(x)$ \ \ if \ \ $x \in {\rm domC}(f) = P$, 

\smallskip

 \ \ \ $\tilde{f}(x) = 0$ \ \  \ \ if \ \ $x \in Q$.

\smallskip

\noindent Hence, $f \subseteq \tilde{f}$; \ so for the function $\varrho$ we
can just pick the identity function, i.e., $\varrho(\tilde{f}) = \tilde{f}$.
Then we have: $\,f = \tilde{f} \ $ iff $ \ f$ is total (i.e., ${\rm domC}(f)$
is a maximal prefix code); hence $\tilde{\tilde{f}} = \tilde{f}$. 

The completion operation $\widetilde{(.)}$ is {\em not injective}. In
particular, $\,\tilde{f}_1$ $=$ $\tilde{f}_2\,$ if the following two 
conditions hold: (1) $\,f_1 \subseteq f_2$; (2) if $f_1(x)$ is undefined then 
$f_2(x)$ is undefined or $f_2(x) \in 0A^*$.

The operation $\widetilde{(.)}$ is {\em not homomorphic}. For example, if 
$f = (1 \to 0)$ and $g = (0 \to 1)$ with $A = \{0,1\}$, then 
$gf(.) = (1 \to 1)$, so $\,\widetilde{gf}(.)$ $=$ ${\mathbb 1}_{\{0,1\}}$.
And $\tilde{f}$ has the table $\{(1,0), (0,0)\}$ and $\tilde{g}$ has the 
table $\{(0,1), (1,0)\}$, so $\tilde{g}\,\tilde{f}(.)$ has the table
$\{(0,1), (1,1)\}$.
Hence $\,\widetilde{gf}(.)$ $\ne$ $\tilde{g}\,\tilde{f}(.)$.

\bigskip

Before describing another completion we give a preliminary construction.
For any prefix code $P \subseteq \{0,1\}^*$, $ \ {\sf spref}(P)\,$ denotes 
the set of {\em strict prefixes} of the elements of $P$ (i.e., the prefixes
not in $P$).

\begin{lem} \label{prefcodeLarger} {\bf (expanding a maximal prefix code 
to a larger alphabet).}

\smallskip

\noindent Let $\{0,1,2\}$ be a 3-letter alphabet. If $P$ is a finite prefix 
code in $\{0,1\}^*$ then

\smallskip

 \ \ \  \ \ \ $P \,\cup\, ({\sf spref}(P)) \,2\,$

\smallskip

\noindent is a finite prefix code in $\{0,1,2\}^*$, which is maximal in 
$\{0,1,2\}^*\,$ iff $\,P$ is maximal in $\{0,1\}^*$. 
\end{lem}
{\sc Proof.} See \cite[{\small Lemma 9.1}]{BiCoNP} and 
\cite[{\small Lemma 1.5}]{BiNewEmb}.)
 \ \ \  $\Box$

\bigskip 

\noindent {\bf 2. An injective non-homomorphic completion of
${\cal RM}_2^{\sf fin}$ or of $M_{2,1}$ }

\smallskip

\noindent Every element $f \in {\cal RM}_2^{\sf fin}$ (or of $M_{2,1}$) is
completed to an element $\ov{f} \in$ ${\sf tot}{\cal RM}_2^{\sf fin}$
(or in ${\sf tot}M_{2,1}$) as follows. 
We first choose a finite complementary prefix code 
$\,Q \subseteq \{0,1\}^*$ of $\,{\rm domC}(f)$. Then $\,P =$ 
${\rm domC}(f) \,\cup\,Q$, ($\subseteq$ $\{0,1\}^*$) is a finite maximal 
prefix code of $\{0,1\}^*$, and $ \ P \,\cup\, ({\sf spref}(P)) \,\bot\,$ 
is a finite maximal prefix code of $\{0,1,\bot \}^*$ \ (by Lemma 
\ref{prefcodeLarger}). \ Now $\ov{f}$ is the right-ideal morphism of 
$\{0,1,\bot \}^*$ defined as follows: 

\medskip

 \ \ \  \ \ \ ${\rm domC}(\ov{f}) \,=\, {\rm domC}(f) \,\cup\, Q$   
    $\,\cup\,$
    $\big( {\sf spref}({\rm domC}(f) \,\cup\, Q)\big) \cdot \bot \,$; 

\smallskip

 \ \ \  \ \ \ $\ov{f}(x) \,=\, f(x)$ \ \ \ if $x \in {\rm domC}(f)$;

\smallskip

 \ \ \  \ \ \ $\ov{f}(x) \,=\, \bot$ \ \ \ if $x \in Q$   $\,\cup\,$
       $\big( {\sf spref}({\rm domC}(f) \,\cup\, Q)\big) \cdot \bot \,$. 
 
\medskip 

\noindent The function $f \in {\cal RM}_2^{\sf fin}$  $\longmapsto$ $\ov{f}$
$\in$ ${\sf tot}{\cal RM}_3^{\sf fin}\,$ is total and injective; it is 
not a monoid morphism, since $M_{2,1}$ is not embeddable into 
${\sf tot}M_{3,1}\,$ (by Theorem \ref{THMMnotintotM}). Taking $\varrho$ to 
be the identity function, we obtain a completion $\,\ov{(.)}:$ 
${\cal RM}_2^{\sf fin} \to {\sf tot}{\cal RM}_3^{\sf fin}\,$ (and
$M_{2,1} \to {\sf tot}M_{3,1}$).

\smallskip

By encoding the elements of ${\sf tot}{\cal RM}_3^{\sf fin}$ as in Prop.\ 
\ref{PROPtotMktotM2}, $\,{\sf tot}{\cal RM}_3^{\sf fin}$ is embedded into 
${\sf tot}{\cal RM}_2^{\sf fin}$.  This yields a (non-homomorphic) completion 
of ${\cal RM}_2^{\sf fin}$ in ${\sf tot}{\cal RM}_2^{\sf fin}$. Taking 
$\varrho$ to be the inverse of the embedding function 
$\,{\sf tot}{\cal RM}_3^{\sf fin} \to {\sf tot}{\cal RM}_2^{\sf fin}$, we 
obtain a completion $\,\ov{(.)}: {\cal RM}_2^{\sf fin} \to$ 
${\sf tot}{\cal RM}_2^{\sf fin}$.

\bigskip

Although deterministic completions are simple and fairly natural, they have 
some drawbacks: They cannot be made homomorphic 
(if they are injective). It is difficult to compute $\ov{f}$ or $\tilde{f}\,$
if $f$ is given by a sequence of generators (in $\Gamma$ or in 
$\Gamma \cup \tau$); even when $f$ is given by a table, it is difficult to 
find $\ov{f}$ or $\tilde{f}\,$ because we do not have an efficient method for
finding a complementary prefix code.

\subsection{Inverse homomorphic completion}

We construct a general completion operation that has (inverse) homomorphic 
properties, and that provides an efficient computation of $\ov{f}\,$ as a 
word of generators, if $f$ is given by a word of generators.  
Here we will work with $M_{h,1}$, not with ${\cal RM}_h^{\sf fin}$.

\begin{defn} \label{DEFgeneralcompl} {\bf (completion defined by an 
inverse homomorphism).}

\noindent {\small \rm (1)}  An {\em inverse homomorphic completion} of 
$M_{h,1}$ in ${\sf tot} M_{k,1}$ is a surjective semigroup homomorphism 

\smallskip

 \ \ \  \ \ \  \ \ \ $\varrho: \ {\rm Dom}(\varrho) \,\subseteq\,$
${\sf tot}M_{k,1} \ \twoheadrightarrow \ M_{h,1} = {\rm Im}(\varrho)$.

\medskip

\noindent The inverse homomorphic completion of ${\sf plep}M_{2,1}$ in 
${\sf tlep}M_{2,1}$ is defined is the same way as for $M_{2,1}$ in 
${\sf tot}M_{2,1}$. 

\medskip

\noindent {\small \rm (2)}  A {\em completion} of $f \in M_{2,1}$ (or 
${\sf plep}M_{2,1}$) is the choice of any element in $\varrho^{-1}(f)$. 
\end{defn}
The homomorphism $\varrho$ is an {\em inverse completion}; one could also
call it a ``completion removal operation'', or an ``un-completion''.
The completion of $f$ is nondeterministic (being a choice within 
$\varrho^{-1}(f)$, whereas the inverse completion $\varrho(f)$ is
deterministic.

\medskip

\noindent {\bf Remark:} A more general inverse homomorphic completion of 
$M_{h,1}$ in ${\sf tot} M_{k,1}\,$ can be defined by a semigroup 
homomorphism (not necessarily surjective) \ $\varrho:\,$ 
${\rm Dom}(\varrho) \ (\,\subseteq\, {\sf tot}M_{k,1}) \ \to \ M_{h,1}$, 
 \ such that for every $f \in M_{h,1}$ there exists
$s \in {\rm Dom}(\varrho)$ satisfying
 \ $f \,\subseteq\, \varrho(s)$.

\smallskip

\begin{thm} \label{THMnondetCompletion} {\bf (inverse homomorphic
completion).}

\smallskip

\noindent {\small \rm (1)} \ The monoid $M_{2,1}$ is a homomorphic image of 
a submonoid of $\,{\sf tot}M_{2,1}$.  \ More precisely, 

\medskip

\hspace{1,1cm} $f \,\in\, F_0 \, =\,$
${\sf Fix}_{{\sf tot}M_{2,1}}\!(0\!\cdot\!\{0,1\}^{\omega})$
$ \ \ (\,\subseteq\, {\sf tot}M_{2,1})$ 

\medskip

\hspace{1,5cm} {\Large  $\stackrel{\rho}{\twoheadrightarrow}$ \ }
$f \ \cap \ 1\,\{0,1\}^{\omega}\!\times\!1\,\{0,1\}^{\omega}$ 

\smallskip 

\hspace{2,6cm} $\in \ $ 
$\big\{\,(f$  $\,\cap\,$
$1\,\{0,1\}^{\omega}\!\times\!1\,\{0,1\}^{\omega})$
  $\,\in M_{2,1} \,:\,  f \in F_0 \big\}$

\medskip

\hspace{1,5cm} {\Large  $\stackrel{\iota}{\to}$ \ }
$M_{2,1}\,$,

\medskip

\noindent where  

\medskip

\hspace{1,2cm} 
${\sf Fix}_{{\sf tot}M_{2,1}}\!(0\!\cdot\!\{0,1\}^{\omega})$
$\,=\,$
$\{f \in {\sf tot}M_{2,1} :\,$
$(\forall z \in \{0,1\}^{\omega})\,[\,f(0z) = 0z \,]\,\}$ \ \ (fixator);

\medskip

\hspace{1,2cm} $\rho$ is a surjective monoid homomorphism; 

\medskip

\hspace{1,2cm} $\iota$ is a monoid isomorphism;

\medskip

\hspace{1,2cm} 
$\,\varrho = \iota \circ \rho: \, F_0 \twoheadrightarrow M_{2,1}\,$ is a 
surjective monoid homomorphism.

\medskip

\noindent The restriction of $\varrho$ to 
$\,\varrho^{-1}({\sf tot}M_{2,1})\,$ is injective, i.e., elements of 
${\sf tot}M_{2,1}$ have a unique completion. 

\bigskip

\noindent {\small \rm (2)} \ The monoid ${\sf plep}M_{2,1}$ is a 
homomorphic image of a submonoid of $\,{\sf tlep}M_{2,1}$. \ More precisely,

\medskip

\hspace{1,1cm} $f \,\in\, S_0 \, =\,$
${\sf Stab}_{{\sf tlep}M_{2,1}}\!(0\!\cdot\!\{0,1\}^{\omega})$
$ \ \ (\,\subseteq\, {\sf tlep}M_{2,1})$

\medskip

\hspace{1,5cm} {\Large  $\stackrel{\rho}{\twoheadrightarrow}$ \ }
 $\big\{\,(f$  $\,\cap\,$
$1\,\{0,1\}^{\omega}\!\times\!1\,\{0,1\}^{\omega})$
  $\,\in {\sf plep}M_{2,1} \,:\,  f \in S_0 \big\}$

\medskip

\hspace{1,5cm} {\Large  $\stackrel{\iota}{\to}$ \ }
${\sf plep}M_{2,1}\,$,

\medskip

\noindent where

\medskip

\hspace{1,2cm}
${\sf Stab}_{{\sf tlep}M_{2,1}}\!(0\!\cdot\!\{0,1\}^{\omega})$
$\,=\,$
$\{f \in {\sf tlep}M_{2,1} :\,$
$(\forall z \in \{0,1\}^{\omega})\,[\,f(0z) \in 0\,\{0,1\}^{\omega}\,]\,\}$
 \ \ (stabilizer);

\medskip

\hspace{1,2cm} $\iota$ and $\varrho\,$ are as in {\rm (1)}

\bigskip

\noindent {\small \rm (3)} \ The same results hold for ${\cal RM}^{\sf fin}$
(instead of $M_{2,1}$), and ${\sf pfl}{\cal RM}^{\sf fin}$ (instead of 
${\sf plep}M_{2,1}$).  
\end{thm} 
\noindent {\sc Proof.} In (1) the proof works in the same way for 
${\cal RM}^{\sf fin}$ and ${\sf tot}{\cal RM}^{\sf fin}$ as for $M_{2,1}$
and ${\sf tot}M_{2,1}$. And in (2) the proof works in the same way for
${\sf pfl}{\cal RM}^{\sf fin}$ and ${\sf tfl}{\cal RM}^{\sf fin}$ as for
${\sf plep}M_{2,1}$ and ${\sf tlep}M_{2,1}$. We only give the proof for 
$M_{2,1}$ and ${\sf plep}M_{2,1}$. For ${\cal RM}^{\sf fin}$, 
$\{0,1\}^{\omega}$ is replaced by $\{0,1\}^*$.

\medskip

\noindent (1) A function $f \in {\sf tot}M_{2,1}$ belongs to 
the fixator $F_0$\, iff $\,f$ has a table of the form 
\[\left[ \begin{array}{ccccccc}
0 & 1x_1 & \ldots & 1x_k & 1x'_1 &\ \ldots \ & 1x'_h \\ 
0 & 1y_1 & \ldots & 1y_k & 0y'_1 &\ \ldots \ & 0y'_h
\end{array}        \right]\, ,
\]
\noindent where $ \ {\rm domC}(f) = \{0\} \,\cup\, 1P \,\cup\, 1P'$, and 
where \ $P = \{x_1, \,\ldots\,, x_k\}$,
  \ $P' = \{x'_1, \,\ldots\,, x'_h\}\,$ are any finite prefix codes such 
 that $P\,\{0,1\}^* \,\cap\, P'\,\{0,1\}^* = \0$.
Moreover, $\,\{y_1, \,\ldots\, , y_k\}\,$ is any set of size $\,\le |P|$, 
and $\,\{y'_1, \,\ldots\, , y'_h\}\,$ is any set of size $\,\le |P'|$. 
Since $f$ is total, ${\rm domC}(f)$ is a maximal finite prefix code; hence,
$P \cup P'$ is a maximal finite prefix code, so $P'$ is a complementary 
finite prefix code of $P$. 

\smallskip

The function $\,f^{\rho} \,=\, (f$   $\,\cap\,$ 
$1\,\{0,1\}^{\omega}\!\times\!1\,\{0,1\}^{\omega}) \,\in\, M_{2,1}$
 \ is such that $f^{\rho}$ is undefined if the input, or the output, is in 
$\,0\,\{0,1\}^{\omega}$. Thus the table of $f^{\rho}\,$ is
\[ f^{\rho} \ = \ \left[ \begin{array}{ccc}
1x_1 & \ldots & 1x_k \\               
1y_1 & \ldots & 1y_k 
\end{array}        \right] .
\] 

\noindent The image set of $\rho$ is the set of all elements of $M_{2,1}$ 
that have a table of the above form for $f^{\rho} \ $ (since $P$ can be any
finite prefix code),  

\smallskip

\noindent The function 
\[\iota: \ \left[ \begin{array}{ccc}
1x_1 & \ldots & 1x_k \\
1y_1 & \ldots & 1y_k
\end{array}        \right]
 \ \ \longmapsto 
 \ \ \left[ \begin{array}{ccc}
x_1 & \ldots & x_k \\
y_1 & \ldots & y_k
\end{array}        \right]
\]
\noindent is an isomorphism from $\rho(F_0)$ onto $M_{2,1}$.  It is 
surjective because $\{(x_1,y_1),$ $\,\ldots\, ,$  $(x_k,y_k)\}$ can be the 
table of any element of $M_{2,1}$. 
 \ (Remark: The function $\rho$ could also be defined on all of
${\sf tot}M_{2,1}$, but then it would not be a homomorphism,
since ${\sf tot}M_{2,1}$ is congruence-simple.)

\medskip

Let us prove that the restriction of $\varrho$ to 
$\,\varrho^{-1}({\sf tot}M_{2,1})\,$ is injective: 
 \ If $\varphi = \{(x_1,y_1),\,\ldots\, , (x_k,y_k)\}$ is the table of an
element of ${\sf tot}M_{2,1}$ then $P = \{x_1,\,\ldots\, , y_k\}$ is a
maximal prefix code, and $P' = \0$. Hence $\varrho^{-1}(\varphi)$ consists
of just the element given by the table
$\, \{(0,0), (1x_1, 1y_1),\,\ldots\, , (1x_k, 1y_k)\}$. 

\medskip

\noindent We still need to prove the following. 

\medskip

\noindent {\sf Claim:} \ \ $\rho: \,f \in F_0 \longmapsto f^{\rho}$ 
$\in \rho(F_0)$ \ is a semigroup homomorphism.

\smallskip

\noindent Proof of Claim: 
For $g,f \in F_0$ and $x \in \{0,1\}^{\omega}$ there are several cases. 
We use parentheses $(.)$ for application of functions, and we use brackets 
$[.]$ for grouping in the composition of functions. 

\smallskip

\noindent $\circ$ \ \ $x,\,f(x) \in 1\,\{0,1\}^{\omega}:$ 

Then $\,g^{\rho}(f^{\rho}(x)) = g^{\rho}(f(x)) = g(f(x))$. 
Since $gf$ is total, $g(f(x))$ belongs either to $1\,\{0,1\}^{\omega}$
or to $0\,\{0,1\}^{\omega}$.

If $\,g(f(x)) \in 1\,\{0,1\}^{\omega}$ then $g(f(x)) = [gf]^{\rho}(x)$.
Hence $g^{\rho}(f^{\rho}(x)) = [gf]^{\rho}(x)$.
 
If $\,g(f(x)) \in 0\,\{0,1\}^{\omega}$ then $[gf]^{\rho}(x)$ is undefined.
Moreover, $g^{\rho}(f^{\rho}(x))$  $=$ $g(f(x)) \in 0\,\{0,1\}^{\omega}$,
i.e., the image of $f(x)$ under $g$ is in $0\,\{0,1\}^{\omega}$; therefore,
$g^{\rho}(f^{\rho}(x))$ is undefined. So 
$g^{\rho}(f^{\rho}(x)) = [gf]^{\rho}(x)$, both being undefined.

\smallskip

\noindent $\circ$ \ \ $x \in 1\,\{0,1\}^{\omega}$ and $f(x) = 0 z \in$ 
$0\,\{0,1\}^{\omega}:$

Then $f^{\rho}(x)$ is undefined, hence $\,g^{\rho}(f^{\rho}(x))$ is 
undefined.
And $[gf]^{\rho}(x) \in \{gf(x)\} \cap 1\,\{0,1\}^{\omega}$, where 
$gf(x) = g(0z) \in 0\,\{0,1\}^{\omega}\,$ (since $g \in F_0$); hence 
$[gf]^{\rho}(x) \in \{gf(x)\} \cap 1\,\{0,1\}^{\omega}$ $=$  $\0$, 
so $[gf]^{\rho}(x)$ is also undefined. 

\smallskip

\noindent $\circ$ \ \ $x = 0 z \in 0\,\{0,1\}^{\omega}:$ 

Then $f^{\rho}(x)$ is undefined, hence $\,g^{\rho}(f^{\rho}(x))$ is 
undefined.
Since $g, f \in F_0$, $gf(0z) = 0z$.
Now $[gf]^{\rho}(x) \in \{gf(0z)\} \cap 1\,\{0,1\}^{\omega}$ $=$
$\{0z\} \cap 1\,\{0,1\}^{\omega}$ $=$  $\0$, so
$[gf]^{\rho}(x)$ is also undefined.

\smallskip

\noindent This proves the Claim and part (1) of the Theorem.

\medskip

\noindent The same proof works for ${\cal RM}^{\sf fin}$ and 
${\sf tot}{\cal RM}^{\sf fin}$, by replacing $\{0,1\}^{\omega}$ by 
$\{0,1\}^*$.

\medskip

\noindent (2) The proof is similar to the proof of (1), but now we have to
be careful about lengths.

A function $f \in {\sf tlep}M_{2,1}$ belongs to the stabilizer 
$S_0\,$ iff $\,f$ has a table of the form
\[\left[ \begin{array}{ccccccccc}
0u_1 & \ldots & 0u_r & 1x_1 & \ldots & 1x_k & 1x'_1 &\ \ldots \ & 1x'_h \\
0v_1 & \ldots & 0v_r & 1y_1 & \ldots & 1y_k & 0y'_1 &\ \ldots \ & 0y'_h
\end{array}        \right]\, ,
\]
\noindent where $\,{\rm domC}(f) = 0U \,\cup\, 1P \,\cup\, 1P'$, and where
 \ $U = \{u_1, \,\ldots\,, u_r\}$,
 \ $P = \{x_1, \,\ldots\,, x_k\}$,
 \ $P' = \{x'_1, \,\ldots\,, x'_k\}\,$ are finite prefix codes such that
$P\,\{0,1\}^* \,\cap\, P'\,\{0,1\}^* = \0$.
Moreover, $\{y_1, \,\ldots\, , y_k\}$ is a set of size $\,\le |P|$, 
 \ $\{y'_1, \,\ldots\,, y'_h\}$ is a set of size $\,\le |P'|$,
and $\{v_1, \,\ldots\,, v_r\}$ is a set of size $\,\le |U|$.
Since $f$ is total, ${\rm domC}(f)$ is a maximal finite prefix code; hence,
$U$ and $P \cup P'$ are each a maximal finite prefix code, so $P'$ is a
complementary finite prefix code of $P$.  
 \ Since $f$ is {\sf tlep}, we have the additional properties that 
for all $i$ (in the appropriate ranges), 

\smallskip

 \ \ \  \ \ \ $|v_i| - |u_i|$  $\,=\,$  $|y_i| - |x_i|$  $\,=\,$  
$|y_i'| - |x_i'|$  $\,=\,$  $\delta(f)$.

\smallskip

\noindent The remainder of the proof is the same as for (1).  

\medskip

\noindent The same proof works for ${\sf pfl}{\cal RM}^{\sf fin}$ and 
${\sf tfl}{\cal RM}^{\sf fin}$, by replacing $\{0,1\}^{\omega}$ by
$\{0,1\}^*$.  
 \ \ \  \ \ \ $\Box$

\bigskip

\noindent {\bf Remark.} The reason why the fixator is not used in (2) is 
that we want $\ov{f}$ to be length-equality preserving ({\sf lep}).
As a result, in (2), the restriction of $\varrho$ to
$\,\varrho^{-1}({\sf tlep}M_{2,1})\,$ is not injective.

\bigskip

\medskip

\noindent {\bf Inverse-homomorphic nondeterministic completion algorithm:} 

\smallskip

For a finite prefix code $P \subseteq A^*$ we define the 
{\em size of} $P\,$ by $ \ \|P\| = \,\sum_{x \in P} |x|$, \ where $|x|$ is 
the length of $x$.

For a right-ideal morphism $g \in {\cal RM}^{\sf fin}$ given by a table 
$\,\{(x_i, y_i) : 1 \le i \le k\}$, the {\em size of the table} is 
defined by \ $\|g\| = \,\sum_{i=1}^k |x_i y_i|$.

By \cite[{\small Def.\ 3.29, Lemma 3.30}]{BiLR} and 
\cite[{\small Prop.\ 5.2}]{BiEvalProb}, every finite prefix code
$P \subseteq A^*$ has a complementary prefix code $Q \subseteq A^*$ such  
that $\,{\rm maxlen}(P) = {\rm maxlen}(Q)$. 
If $P \subseteq A^m$ then $\,A^m \minus P\,$ is a complementary prefix code
of $P$.

In Prop.\ \ref{NondetComplAlgo}, $\varrho$ is as in Theorem 
\ref{THMnondetCompletion}.

\begin{pro} \label{NondetComplAlgo}\!\!.

\smallskip

\noindent {\small \rm (1)} ${\cal RM}^{\sf fin}$ and 
${\sf tot}{\cal RM}^{\sf fin}:$

The following nondeterministic algorithm, on input $g \in M_{2,1}$
(given by the table of a right-ideal morphism in ${\cal RM}^{\sf fin}$ that
represents $g$), outputs an element $\,\ov{g} \in$  ${\sf tot}M_{2,1}\,$ 
(given by the table of a right-ideal morphism in  
${\sf tot}{\cal RM}^{\sf fin}$), such that $\,\varrho(\,\ov{g}\,) = g$.

\medskip

\noindent {\small \rm (2)} ${\sf pfl}{\cal RM}^{\sf fin}$ and 
${\sf tfl}{\cal RM}^{\sf fin}:$

A variant of this nondeterministic algorithm, on input $g \in$ 
${\sf plep}M_{2,1}$ (given by the table of a right-ideal morphism in
${\sf pfl}{\cal RM}^{\sf fin}$), outputs an element $\,\ov{g} \in$ 
${\sf tlep}M_{2,1}$ (given by the table of a right-ideal morphism in 
${\sf tfl}{\cal RM}^{\sf fin}$), such that $\,\varrho(\,\ov{g}\,) = g$.
The algorithm can find such a $\,\ov{g}\,$ in time $O(\|g\|)$.
\end{pro}
Algorithm:

\smallskip

\noindent For any function $g \in {\cal RM}^{\sf fin}$ with table
\[ g \,=\, \left[ \begin{array}{ccc}
x_1 & \ldots & x_k \\
y_1 & \ldots & y_k
\end{array}  \right] 
\] 
\noindent where $\,{\rm domC}(g) = P = \{x_1, \,\ldots\,, x_k\}$ is a 
finite prefix code, a completion can be obtained by choosing 
\[ \ov{g} \ = \  \left[ \begin{array}{ccccccc}
0 & 1x_1 & \ldots & 1x_k & 1x'_1 &\ \ldots \ & 1x'_h \\
0 & 1y_1 & \ldots & 1y_k & 0y'_1 &\ \ldots \ & 0y'_h
\end{array}        \right]\, ,
\]
\noindent where $\,P' = \{x'_1, \,\ldots\,, x'_h\}\,$ is an arbitrary 
complementary finite prefix code of $P$, and 
$\,g' = \{(x'_1,y'_1),$ $\,\ldots\,,$  $(x'_h,y'_h)\}\,$ is
an arbitrary element of ${\cal RM}^{\sf fin}$ with domain code $P'$.
(If $g$ is total, i.e.\ $P$ is maximal, then $P' = \0$, so $g'$ is the 
empty function $\theta$.)

\medskip

If $g \in {\sf plep}M_{2,1}$, represented by a right-ideal morphism in
${\sf pfl}{\cal RM}^{\sf fin}$, then $P = \{x_1, \,\ldots\,, x_k\}$
$\subseteq \{0,1\}^m$, and $\{y_1, \,\ldots\,, y_k\} \in \{0,1\}^n$, with
$n = m + \delta(g)$. We chose 
\[ \ov{g} \ = \  \left[ \begin{array}{ccccccccc}
0u_1 & \ldots & 0u_r & 1x_1 & \ldots & 1x_k & 1x'_1 &\ \ldots \ & 1x'_h \\
0v_1 & \ldots & 0v_r & 1y_1 & \ldots & 1y_k & 0y'_1 &\ \ldots \ & 0y'_h
\end{array}        \right]\, ,
\]
\noindent where $U = \{0,1\}^m$ (no choice). We arbitrarily choose a set
$\,\{v_1, \,\ldots\,, v_r\}$ $\subseteq$  $\{0,1\}^n$, and we choose the
pairing arbitrarily so that $\,\{(u_1,v_1),\,\ldots\,, (u_r,v_r)\}$
$\in {\sf tfl}{\cal RM}^{\sf fin}$.
We let $P' = \{0,1\}^m \minus P\,$ (no choice), we arbitrarily choose a set
$\{y'_1, \,\ldots\,, y'_h\}$ $\subseteq$ $\{0,1\}^n$, and we choose the 
pairing $\,g' = \{(x'_i,y'_i) : i = 1,\ldots,h\}$ in 
${\sf pfl}{\cal RM}^{\sf fin}$.
The resulting table for $\ov{g}$ gives a right-ideal morphism in 
${\sf tfl}{\cal RM}^{\sf fin}$, with fixed input length $m\!+\!1$ and 
fixed output length $n\!+\!1$.

In the case of ${\sf pfl}M_{2,1}$, if $g$ is given by a table with
${\rm domC}(g) \subseteq \{0,1\}^m$, then the complementary prefix code of
${\rm domC}(g)$ is easily found.

\smallskip

The choices in the construction of $\,\ov{g}$ from $g$ are of course 
nondeterministic, but the inverse process $\,\varrho: \ov{g} \mapsto g\,$ 
is a function, and a homomorphism $F_0 \twoheadrightarrow M_{2,1}$,
respectively $S_0 \twoheadrightarrow {\sf tfl}M_{2,1}$. 
 \ \ \  \ \ \  $\Box$

\medskip

\noindent {\bf Remark.} In part (1) of the algorithm we do not know whether
the time complexity is polynomial in $\|g\|$ because we do not know an 
efficient algorithm for finding a finite complementary prefix code of a 
finite prefix code.

\bigskip

Besides the above nondeterministic completion procedure, there are many 
other ones; there are many monoid homomorphisms from submonoids of 
${\sf tot}M_{2,1}$ onto $M_{2,1}$, respectively submonoids of 
${\sf pfl}M_{2,1}$ onto ${\sf tfl}M_{2,1}$.

\smallskip

The above algorithm does not run in polynomial time (except in the {\sf pfl}
case). In the next two subsections we give algorithms that run in log-space, 
but use a different form of input than the above algorithm: $g$ will either 
be given by a generator word over $\Gamma$ or $\Gamma \cup \tau$, or (in the 
{\sf pfl} case) by an acyclic circuit.

\subsection{Generator-based completion}

\noindent {\bf Generator-based completion for $M_{2,1}$ }

\smallskip

Let $\,\varrho: F_0 \twoheadrightarrow M_{2,1}$, with $F_0 \le$
${\sf tot}M_{2,1}$, be as in Theorem \ref{THMnondetCompletion}.  Let 
$\Gamma_{\rm\!M}$ be any finite set such that $\Gamma_{\rm\!M} \cup \tau$
generates $M_{2,1}$; similarly, let $\Gamma_{\rm\!tot}$ be any finite set 
such that $\Gamma_{\rm\!tot} \cup \tau$ generates ${\sf tot}M_{2,1}$.

For each generator $\gamma$ in the finite set $\Gamma_{\rm\!M}\,$ we 
choose a word $\ov{\gamma} \in$ $(\Gamma_{\rm\!tot} \cup \tau)^*$ such that 
$E_{\ov{\gamma}}$  $\in$ $\varrho^{-1}(\gamma) \subseteq F_0$.

For a transposition $\tau_{i,j} \in \tau$ with $1 \le i < j$ we choose
$\ov{\tau}_{i,j}$ as follows: For all $x \in \{0,1\}^{\ge j}$, 

\smallskip

 \ \ \  \ \ \  $\ov{\tau}_{i,j}(0x) \,=\, 0x$,  

\smallskip

 \ \ \  \ \ \  $\ov{\tau}_{i,j}(1x) \,=\, \tau_{i+1,j+1}(1x)$,

\smallskip

\noindent so $\,{\rm domC}(\ov{\tau}_{i,j}) = \{0,1\}^{j+1}$ $=$
${\rm imC}(\ov{\tau}_{i,j})$.
Then $\,\ov{\tau}_{i,j} \in \varrho^{-1}(\tau_{i,j})$ $\subseteq F_0$, so 
$\ov{\tau}_{i,j}$ is indeed a completion of $\tau_{i,j}$
according to Theorem \ref{THMnondetCompletion}(1).
To represent $\ov{\tau}_{i,j}$ be a word in 
$(\Gamma_{\rm\!tot} \cup \tau)^*$ we use the following formula:
For every string $bx \in \{0,1\}^{j+1}$ with $b \in \{0,1\}$,

\smallskip

 \ \ \  \ \ \ $\ov{\tau}_{i,j}(bx) \ = \ $
$\big( (b')^{j+1} \ \ {\sf and}_{j+1} \ \ bx\big)$
$ \ \ {\sf or}_{j+1} \ \ $
$\big( b^{j+1} \ \ {\sf and}_{j+1} \ \ \tau_{i+1,j+1}(bx) \big)\,$,

\smallskip

\noindent where $b' = {\sf not}(b)$, $z^n$ denotes the concatenation of
$n$ copies of $z$, ${\sf and}_n$ (or ${\sf or}_n$) is the bitwise {\sf and}
(respectively {\sf or}) of two bitstrings of length $n$.

Based on this formula, one can easily build an acyclic boolean circuit for 
$\ov{\tau}_{i,j}$, and this circuit then yields a word in 
$(\Gamma_{\rm\!tot} \cup \tau)^*\,$ (by Theorem \ref{circ2lepM}(1)). Thus
for $\tau_{i,j}$ we can construct a word for $\ov{\tau}_{i,j}$ in space
$O(\log j)$.  

This proves the following:

\begin{cor} \label{CORcompletGener}
 \ For any $\,w = w_n \,\ldots\, w_1$ $\in (\Gamma_{\rm\!M} \cup \tau)^*$
with $w_i \in \Gamma_{\rm\!M} \cup \tau$ (for $i = 1, \ldots, n$), let 
$\,E_w = w_n \circ \ldots \circ w_1(.)$ $\,\in\,$  $M_{k,1}$ be the function
generated.
 \ By completing each generator and concatenating we obtain $\,W \,=\,$ 
$\ov{w}_n \,\ldots\, \ov{w}_1$ $\in$  $(\Gamma_{\rm\!tot} \cup \tau)^*$, 
which generates the total function $\,E_W \, =\,$  
$\ov{w}_n \circ \ldots \circ \ov{w}_1(.)$ 
$\,\in\,$  $F_0\, \subseteq {\sf tot}M_{k,1}$.

\smallskip

Then $\,E_W\,$ is a completion of $\,E_w\,$, \ i.e., 
  \ $\varrho(E_W) = E_w$.  

\smallskip

The word $W$ can be computed from $w$ in log-space.
 \ \ \  \ \ \  \ \ \  $\Box$
\end{cor}
In other words, a completion of a function $f \in M_{k,1}$, given by a
sequence of generators, can be efficiently obtained by completing the 
generators and composing them. 
 \ The same approach works for ${\cal RM}^{\sf fin}$.

\bigskip

\noindent {\bf Generator-based completion for ${\sf pfl}{\cal RM}^{\sf fin}$}

\smallskip 

For ${\sf pfl}{\cal RM}^{\sf fin}$ we obtain more results, based on the 
generating sets $\Gamma_{\rm \!pfl} \cup \tau$, respectively 
$\Gamma_{\rm \!tfl} \cup \tau$, with $\Gamma_{\rm \!pfl}$ and 
$\Gamma_{\rm \!tfl}$ finite.

For the finitely many generators $\gamma \in \Gamma_{\rm \!pfl}$ we choose
$\ov{\gamma} \in \varrho^{-1}(\gamma) \subseteq S_0$ nondeterministically as in 
Prop.\ \ref{NondetComplAlgo}. 
For a position transposition $\tau_{i,j} \in \tau$ with $1 \le i < j$ we
choose

\smallskip

 \ \ \  \ \ \  $\ov{\tau}_{i,j} = \tau_{i+1,j+1}$.

\smallskip

\noindent For this choice, $\,\tau_{i+1,j+1}(0z) \in$ $0 \{0,1\}^*\,$ for 
all $z \in \{0,1\}^{\ge j+1}$ and for all $z \in$ $\{0,1\}^{\omega}$; so 
$\tau_{i+1,j+1}$  $\in$  $S_0$ ($\,\subseteq$ 
${\sf tfl}{\cal RM}^{\sf fin}$), and $\,\varrho(\,\tau_{i+1,j+1}\,)$  $=$
$\tau_{i,j}$. 

\medskip

\noindent Remarks: (1) Many completions of $\tau_{i,j}$ are possible,
but the above one has the advantage of being itself in $\tau$.
(2) Note that $\,\tau_{i+1,j+1} \,\not\in\, F_0$, so the above completion
of $\tau_{i,j}$ only works for ${\sf plep}M_{2,1}$, not for $M_{2,1}$.

\smallskip

\begin{pro} \label{PROPcompletionWord} {\bf (completion in terms of
generators).}

\smallskip

\noindent For every $w \in (\Gamma_{\sf\!pfl} \cup \tau)^*$ there exists a 
word $\,\ov{w} \in (\Gamma_{\sf\!tfl} \cup \tau)^*$ such that

\medskip

\noindent {\small \rm (1)} \ $E_{\ov{w}} \,\in\, S_0\,$ $(\,\subseteq$ 
${\sf tfl}{\cal RM}^{\sf fin})$ \ is a completion of $E_w$, and $\ov{w}$ can
be found from $w$ by a deterministic algorithm that runs in linear time and 
log-space.

\medskip

\noindent {\small \rm (2)} \ $\,|w|_{\Gamma_{\sf\!pfl} \cup \tau}$ 
 $\,\le\,$
 $|\ov{w}|_{\Gamma_{\sf\!tfl} \cup \tau}$  $\,\le\,$
 $c \, |w|_{\Gamma_{\sf\!pfl} \cup \tau}\,$, 
 \ for some constant $c \ge 1$; 

\smallskip

 \ $\,{\rm maxindex}_{\tau}(w)$  $\,\le\,$
${\rm maxindex}_{\tau}(\ov{w}) \,\le\, {\rm maxindex}_{\tau}(w) + 1\,$;

\smallskip

 \ $\,\ell_{\rm in}(\ov{w}) \,=\, \ell_{\rm in}(w) + 1$,
 \ \ and \ \ $\ell_{\rm out}(\ov{w}) + 1 \,=\, \ell_{\rm out}(w) + 1$;
 \ hence,

\smallskip

 \ $\,|w| \le |\ov{w}| \le c\,|w|\,$, \ for some constant $c \ge 1$.

\medskip

\noindent {\small \rm (3)} \ $\,\ov{w} \in (\Gamma_{\sf\!tfl} \cup \tau)^*$
can be chosen so that in addition to the above,

\smallskip

 \ $\,{\sf depth}(\ov{w}) \,\le\,$
  $c \,+\, \max\{{\sf depth}(w),\, \log_2 |w|_{\Gamma_{\sf\!pfl}} \}\,$, 
  \ \ for some constant $c \ge 1$.
\end{pro}
{\sc Proof.} (1) This is proved in the same way as Corollary 
\ref{CORcompletGener}, except that a simpler form for $\ov{\tau}_{i,j}$ is
now used.
The construction of $\ov{w}$ implies that it can be computed from $w$ in
linear times and log-space (if every transposition $\tau_{i,j}$ is encoded
in the form $0^i 1^j$ over the alphabet $\{0,1\} \cup \Gamma_{\sf\!pfl}$, 
assuming $\,\{0,1\}^* \cap \Gamma_{\sf\!pfl} = \0$).

\medskip

\noindent (2) The claims about $|\ov{w}|_{\Gamma_{\sf\!tfl} \cup \tau}$
and ${\rm maxindex}_{\tau}(\ov{w})$ follow immediately from the construction
of $\ov{w}$.

The claims about $\ell_{\rm in}(\ov{w})$ and $\ell_{\rm out}(\ov{w})$ follow
from the definition. Let $x \in$ $\{0,1\}^m$. By induction on $i$ $\in$
$\{1,\ldots,k\}$, one proves that
$\,\ov{w}_i \circ \ldots \circ \ov{w}_1(1x)$ is defined iff
$w_i \circ \ldots \circ w_1(x)$ ($\in \{0,1\}^n$) is defined.
And for $x \not\in {\rm Dom}(w_i \circ \ldots \circ w_1)$,
$\,\ov{w}_i \circ \ldots \circ \ov{w}_1(1x) \in 0 \{0,1\}^n$.
Moreover, $\,\ov{w}_i \circ \ldots \circ \ov{w}_1(0 \{0,1\}^m)$
$\subseteq$  $0 \{0,1\}^n$.
Hence, $\ell_{\rm in}(\ov{w}) = \ell_{\rm in}(w) + 1$ and
$\ell_{\rm out}(\ov{w}) = \ell_{\rm out}(w) + 1$.

\smallskip

\noindent (3) The claim about ${\sf depth}(\ov{w})$ follows from Prop.\
\ref{PROPcompletCircuit}(2) below, and and the definition of depth in terms
of circuits (Def. \ref{DEFdepthRM}).
 \ \ \  \ \ \ $\Box$

\subsection{Circuit-based completion}

\noindent We now construct circuit completions that do not use generators 
nor Theorems \ref{circ2lepM} and \ref{circ2PartLepM} (about the conversion 
between words and circuits or partial circuits). 
By a circuit we always mean an acyclic circuit.

\begin{pro} \label{PROPcompletCircuit} {\bf (boolean circuit for the 
completion).}

\smallskip

\noindent Let $\varrho$ be as in Theorem {\rm \ref{THMnondetCompletion}}.
Suppose that $f \in {\sf pfl}{\cal RM}^{\sf fin}$ has a
{\em partial circuit} $C_{\!f}$ with input-length $m$, output-length $n$, 
and size $s$.  \ Then we have:

\medskip

\noindent {\small \rm (1)} \ $f$ has a completion $\ov{f} \,\in\,$ 
$\varrho^{-1}(f) \,\subseteq\, {\sf tfl}{\cal RM}^{\sf fin}\,$ that has a
boolean circuit with input-length $m\!+\!1$,
output-length $n\!+\!1$, and size $\,\le c\,s\,$ (for some constant 
$c \ge 1$).

\medskip

\noindent {\small \rm (2)} \ If, in addition, $C_{\!f}$ has depth $d$, then 
$f$ has a completion $\,\ov{f} \,\in\, \varrho^{-1}(f)\,$ that has a boolean
circuit with depth $\,\le\, O(1)$  $\,+\,$ 
$\max\{d, \ \lceil \log_2 s \rceil$ $+$ 
$2\,\lceil \log_2 (n\!+\!1) \rceil \}\,$, 
input-length $m\!+\!1$, output-length $n\!+\!1$, and size $\,O(s)$.

The boolean circuit for $\ov{f}$ can be computed from $C_{\!f}$ in log-space.

\medskip

\noindent {\small \rm (3)} \ Under the assumptions of {\rm (2)}, the
{\em classical completion} $\,\tilde{f} \in {\sf tfl}{\cal RM}^{\sf fin}\,$
(defined in Subsection {\rm 9.1}) has a boolean circuit with $\,{\rm depth}$
 $\,\le\, O(1)$ $+$  $\max\{d, \ $ 
$\lceil \log_2 s \rceil + 2\,\lceil \log_2 (n\!+\!1) \rceil\}\,$, 
input-length $m$, output-length $n$, and {\em size} $\,\le O(s)$.

The boolean circuit for $\tilde{f}$ can be computed from the
partial circuit $C_{\!f}$ in log-space.
\end{pro}
{\sc Proof.} (1) \ This follows from Prop.\ \ref{PROPcompletionWord}(2) 
(about the word-lengths of $f = E_w$ and $\ov{f} = E_{\ov{w}}$), combined 
with Theorems \ref{circ2lepM} and \ref{circ2PartLepM} (about the connection 
between word-length and circuit-size).
 \ (1) also follows from the construction in part (2) below.

\medskip

\noindent (2) We pick the completion
\[ \ov{f} \ = \  \left[ \begin{array}{ccccccccc}
0u^{(1)}&\ldots&0u^{(r)}&1x^{(1)}&\ldots&1x^{(k)}&1z^{(1)}&\ldots&1z^{(h)}\\
0^{n+1} &\ldots&0^{n+1} &1y^{(1)}&\ldots&1y^{(k)}&0^{n+1} &\ldots & 0^{n+1}
\end{array}        \right]\, ,
\]
\noindent with $\varrho(\ov{f}) = f$. 
Here $U = \{u^{(1)}, \,\ldots\,, u^{(r)}\} = \{0,1\}^m$, with
$\,r = 2^m$; $ \ {\rm domC}(f) = P = \{x^{(1)}, \,\ldots\,, x^{(k)}\}$ 
$\subseteq \{0,1\}^m$, and $\,y^{(i)} = f(x^{(i)}) \in \{0,1\}^n\,$ for 
$i \in [1, k]$; $ \ P' = \{z^{(1)}, \,\ldots\,, z^{(h)}\} = $
$\{0,1\}^m \minus P\,$ (complementary prefix code), with $h = 2^m - k$.
In other words, for all inputs $x_0\,x$ with $x_0 \in \{0,1\}$ and
$x \in \{0,1\}^m:$

\medskip

 \ \ \  \ \ \ $\ov{f}(x_0 x) \,=\, 1\,f(x)$ \ \ if \ \ $x_0 = 1$ \ and
 \ $x \in P$;

\smallskip

 \ \ \ \ \ \ $\ov{f}(x_0 x) \,=\, 0^{n+1}$ \ \ \ if \ \ $x_0 = 0$ \ {\bf or}
 \ $x \in \{0,1\}^m \minus P$.

\bigskip

\noindent {\sf Claim:} Suppose $f$ has a {\em partial circuit} $C_{\!f}$ 
with input-length $m$, output-length $n$, size $s$, and depth $d$.  Then 
there is a boolean circuit $C_G$ with the following properties:

\smallskip

\noindent $\circ$ \ The input-output function of $C_{\!G}$ is 
 \ $G: \ x \in \{0,1\}^m  \,\longmapsto\, G(x) \in \{0,1\}^{n+1}\,$, 
 \ where

\smallskip

 \ \ \  \ \ \ $G(x) \,=\, 1\,f(x)$ \ \ if \ $x \in P$, 

\smallskip

 \ \ \  \ \ \ $G(x) \,=\, 0^{n+1}$ \ \ \ if \ $x \in \{0,1\}^m \minus P$.

\smallskip

So, $\,G(x) \,=\, \ov{f}(1x) \ $ for all $x \in \{0,1\}^m$.

\smallskip

\noindent $\circ$ \ The boolean circuit $C_{\!G}$ has size $\,\le c_0\,s$,
and depth $\,\le c_1 + \max\{d, \ $
$\lceil \log_2 s \rceil + \lceil \log_2 n \rceil \}$ 
 \ (for some constants $c_0, c_1 \ge 1$).

\medskip

\noindent Proof of the Claim:  
 \ We modify the partial circuit $C_{\!f}$ by completing every partial gate 
$g$ as follows. To $g$ we add a new output wire that outputs $0$ if the 
gate-output is undefined on the given gate-input; and for that gate-input, 
the already existing output wires of the gate $g$ now receive the value 0.  
If on the given gate-input of $g$, the gate-output is defined, the newly 
added output wire outputs 1.
This modified gate, which is total, is the {\em completion} of the partial 
gate $g$.

The newly added output wires do not feed into any existing gates, but into a
binary tree of {\sf and}-gates, of size $\,\le s$ and depth
$\,\le 1 + \lceil \log_2 s \rceil$. Let $y_0$ be the single output bit of 
this tree of {\sf and}-gates. 
Let $\ov{y}$ be output, other than $y_0$, of the completed circuit; so 
$\ov{y}$ has length $n$.

Next, $n$ copies of $y_0$ are produced by a binary tree of $\,\le n\,$ 
{\sf fork} gates (of depth $1 + \lceil \log_2 n \rceil$). 

The output $\ov{y}$ is combined (using {\sf and}-gates) with the $n$ copies
of $y_0$. In other words, this yields the bitstring 
$ \ y_0^{ \ n} \ {\sf and}_n \ \ov{y}\,$, where $\,{\sf and}_n\,$ is the 
bitwise {\sf and} of two bitstrings of length $n$.  We also output $y_0$ 
itself.  The output of $C_{\!G}$ is thus
 \ $y_0\,(y_0^{\ n} \ {\sf and}_n \ \ov{y})$, which is 
$\,1\,f(x)\,$ if $x \in P$, and 
$\,0^{n+1}\,$ if $\,x \in$ $\{0,1\}^m \minus P$.
This completes the construction of $C_{\!G}$.

\smallskip

The claim about circuit size is straightforward.  

The depth of $C_{\!G}$ is 
$\,\le \max\{d,\, 2 + \lceil \log_2 s \rceil + \lceil \log_2 n \rceil\}$,
since the {\sf and}-tree, followed by the {\sf fork}-tree, can be put in 
parallel with the completed circuit for $G$, except for 
a delay of one gate (this is why there is ``$2+$'' in 
$\,2 + 2\,\lceil \log_2 s \rceil + \lceil \log_2 n \rceil\,$ instead of 
the previous ``$1+$'').
 \ \ \  \ \ \ [End, proof of Claim.]

\bigskip

The definition of $\,\ov{f}: \{0,1\}^{m+1} \to \{0,1\}^{n+1}$ can be
formulated as follows.  For all inputs $x_0\,x$ with $x_0 \in \{0,1\}$ and
$x \in \{0,1\}^m$,

\smallskip

 \ \ \  \ \ \  \ \ \ $\ov{f}(x_0\,x) \ = \ $
  $x_0^{ \ n+1} \ \ {\sf and}_{n+1} \ \ G(x)$,

\smallskip

\noindent where $x_0^{ \ n+1} \in \{0^{n+1}, 1^{n+1}\}\,$ is the
concatenation of $n\!+\!1$ copies of the bit $x_0$.

\smallskip

Let us now describe a log-space algorithm that constructs a boolean circuit 
$C_{\!\ov{f}}$ for $\ov{f}$, based on the partial circuit $C_{\!f}$.
First, the boolean circuit $C_{\!G}$ is constructed as the proof of the 
Claim; on input $x$, $C_{\!G}$ outputs $\ov{f}(1x)$. 
Second, $n\!+\!1$ copies of $x_0$ are made, by using a binary tree of
{\sf fork} gates, of depth $1 + \lceil \log_2 (n\!+\!1) \rceil$. This can
be done in log-space.
Finally, we use the formula
$ \ \ov{f}(x_0 x) \,=\, x_0^{ \ n+1} \ {\sf and}_{n+1} \ G(x)$; 
for this, a row of $n\!+\!1$ {\sf and}-gates combines $x_0^{ \ n+1}$ bitwise
with the output of $C_{\!G}$.

The size of $C_{\!\ov{f}}$ is $O(s)$, and the depth is $\,\le\, 3$ $+$
$\max\{d,\ 2 + \lceil \log_2 s \rceil + \lceil \log_2 n \rceil$
$+$  $\lceil \log_2 (n\!+\!1) \rceil \}$.

\medskip

\noindent (3) The boolean circuit $C_{\!\ov{f}}$ can be modified to a boolean
circuit that it computes $\tilde{f}$, as follows.
On input $x \in \{0,1\}^m$, the circuit prepends $1$, and computes 
$\ov{f}(1x)$. If the result is $1\,f(x)$ the circuit outputs $f(x)$; 
if the output is in $\,0\{0,1\}^*$, the circuit outputs $0^n$. 
 \ \ \  \ \ \ $\Box$

\section{Representation of a function in $M_{2,1}$ or ${\cal RM}^{\sf fin}$
by a \\ 
union of partial circuits with disjoint domains}

\subsection{Preliminaries}

Recall the definitions of $\ell_{\rm in}(C)$, $\ell_{\rm out}(C)\,$ (Def.\ 
\ref{DEF_IOlenCircuits}), circuit size (Def.\ \ref{DEFsizedepthCircuits}),
$\ell_{\rm in}(w)$, $\ell_{\rm out}(w)\,$ (Def.\ \ref{DEF_IOlengthRMfin}), 
and various word-lengths and maxlen($\cal F$) (Def.\ \ref{word_length_size}).

\begin{pro} \label{PROPdomCimCLen}\!\!.

\smallskip

\noindent {\small \rm (1)} For any partial acyclic circuit $C:$

\smallskip

 \ \ \  \ \ \  $\ell_{\rm in}(C), \ \ell_{\rm out}(C) \ \le \ |C|$. 

\medskip

\noindent {\small \rm (2)} Let $\Gamma_{\rm\!RM}$ be a finite set such that 
$\Gamma_{\rm\!RM} \cup \tau$ generates ${\cal RM}^{\sf fin}$. For every word 
$w \in$ $(\Gamma_{\rm\!RM} \cup \tau)^*$ let $ \ |w| \,=\, $
$|w|_{\Gamma_{\rm\!RM} \cup \tau} \,+\, {\rm maxindex}_{\tau}(w)\,$. 
Then we have:

\medskip

 \ \ \  \ \ \  $\ell_{\rm in}(w) \ \le \ c\, |w|_{\Gamma_{\rm\!RM}}$ 
$\,+\,$  ${\rm maxindex}_{\tau}(w)$ \ $\le \  c\,|w|$,  

\smallskip

 \ \ \  \ \ \  $\ell_{\rm out}(w) \ \le \ 2\,c\,|w|_{\Gamma_{\rm\!RM}}$  
$\,+\,$ ${\rm maxindex}_{\tau}(w)$   \ $\le \  2\,c\,|w|$,

\medskip

\noindent where $c = {\rm maxlen}(\Gamma_{\rm\!RM})$.
\end{pro}
{\sc Proof.} (1) The inequality is immediate from the definition of
circuit size (Def.\ \ref{DEFsizedepthCircuits}). 

\smallskip

\noindent (2) Let $\,\ell =$ 
$c\,|w|_{\Gamma_{\rm\!RM}} + {\rm maxindex}_{\tau}(w)$.
For the first inequality it suffices to show that $E_w(x)$ is defined for
every $x \in \{0,1\}^{\ge \ell}$ such that $x\,\{0,1\}^{\omega}$ $\subseteq$
${\rm Dom}(E_w)\,\{0,1\}^{\omega}$.

When $w$ is applied to $x$, the generators in $\tau$ do not cause any 
change in length, and can only have an effect within length 
$\,\le {\rm maxindex}_{\tau}(w)$; indeed, the generators of $w$ in 
$\Gamma_{\rm\!RM}$ cannot cause a decrease of more than 
$\,c\,|w|_{\Gamma_{\rm\!RM}}$, so the string that results at any moment 
has length $\,\ge {\rm maxindex}_{\tau}(w)$.
And the $|w|_{\Gamma_{\rm\!RM}}$ generators in $\Gamma_{\rm\!RM}$ cause a
length-increase of at most $c \cdot |w|_{\Gamma_{\rm\!RM}}$.

Hence, if $|x| \ge \ell$ then $x$ is not too short for any of the
generators to be applied (if $x \{0,1\}^{\omega} \in {\rm Dom}(E_w)$).
Hence $m \le \ell$.
And in the output of the action of $w$ on $x \in \{0,1\}^{\ell}$ is a most
$c \cdot |w|_{\Gamma_{\rm\!RM}}$ longer than $x$.
 \ \ \  \ \ \ $\Box$

\begin{defn} \label{DEFcommInOutLen} {\bf (common input-length and
common output-length).}

\smallskip

\noindent Let $\Gamma_{\sf\!pfl}$ be a finite set such that
$\Gamma_{\sf\!pfl} \cup \tau$ generates ${\sf pfl}{\cal RM}^{\sf fin}$.

\smallskip

\noindent {\small \rm (1)} For a set of words $\,\{u_1, \,\ldots\,, u_k\}$
$\subseteq$ $(\Gamma_{\sf\!pfl} \cup \tau)^*$, the {\em common input-length}
is

\smallskip

\hspace{1,8cm} $\ell_{\rm in}(u_1, \,\ldots\,, u_k)$ $ \ = \ $
$\max\{\ell_{\rm in}(u_1), \,\ldots\,, \ell_{\rm in}(u_k)\}\,$,

\smallskip

where $\ell_{\rm in}(u_i)$ is as in Def.\ {\rm
\ref{DEF_IOlengthRMfin}}.

\medskip

\noindent {\small \rm (2)} If $\,\delta(E_{u_i})\,$ is the same for all
$i \in [1,k]$, then the {\em common output-length} of
$\,\{u_1, \,\ldots\,, u_k\}$ is

\smallskip

\hspace{1,8cm} $\ell_{\rm in}(u_1, \,\ldots\,, u_k) \,+\, \delta(E_{u_i})\,$,

\smallskip

which is the same for all $i$.

\smallskip

If $\delta(E_{u_i})$ is not the same for all $i \in [1,k]$ then
$\,\{u_1, \,\ldots\,, u_k\}$ has no common output-length.
\end{defn}

\begin{lem} \label{LEMcommInOutLen}
 \ Under the conditions of Def.\ {\rm \ref{DEFcommInOutLen}}, the common
input-length $m$ and the common output-length $n$ satisfy

\smallskip

 \ \ \  \ \ \  $m, \ n \ \le \ $
$\max\big\{ \,{\rm maxlen}(\Gamma_{\sf\!pfl})$ $\!\cdot\!$
$|u_i|_{\Gamma_{\sf\!pfl}}$   $\,+\,$
${\rm maxindex}_{\tau}(u_i) \,:\, 1 \le i \le k \,\big\}$.
\end{lem}
{\sc Proof.} This is proved in a similar way as Prop.\ \ref{PROPdomCimCLen}.
 \ \ \ $\Box$

\bigskip

\noindent {\bf Depth of words in
{\boldmath ${\cal RM}^{\sf fin}$ and $M_{2,1}$}}

\smallskip

We saw that words in $(\Gamma_{\sf\!tfl} \cup \tau)^*$ or
$(\Gamma_{\sf\!pfl} \cup \tau)^*$ are closely related to boolean circuits,
respectively partial circuits.
Words in $(\Gamma_{\rm\!RM} \cup \tau)^*$, where $\Gamma_{\rm\!RM}$ is finite 
set such that $\Gamma_{\rm\!RM} \cup \tau$ generates ${\cal RM}^{\sf fin}$ 
(and hence $M_{2,1}$) are a generalization of partial circuits.
However, the previous definition of depth of an element of
$(\Gamma_{\sf\!tfl} \cup \tau)^*$ or $(\Gamma_{\sf\!pfl} \cup \tau)^*$ are
based on circuits that are directly obtained from words of generators.
This approach does not directly apply to $(\Gamma_{\rm\!RM} \cup \tau)^*$ 
since a generator in $\Gamma_{\rm\!RM}$ need not have a fixed input-length 
and a fixed output-length, so no circuit digraph with fixed in- and 
out-degrees can be formed.

Nevertheless, Def.\ \ref{DEFdepthM21wx} gives the depth of a word $w \in$ 
$(\Gamma_{\rm\!RM} \cup \tau)^*\,$ for a given input $x \in \{0,1\}^*\,$.
Let $w = w_k \,\ldots\, w_1$ with $w_i \in \Gamma_{\rm\!RM} \cup \tau$ (for
$1 \le i \le k$). For a fixed input $x \in \{0,1\}^*\,$ of length 
$\,|x| \ge \ell_{\rm in}(w) = m$, each $w_i(.)$ has a fixed input-length, 
and a fixed output-length. So for a given $w \in$ 
$(\Gamma_{\rm\!RM} \cup \tau)^*\,$ and $x \in \{0,1\}^{\ge m}$, there is a 
partial circuit $C_{w,x}$ such that $C_{w,x}(x) = E_w(x)\,$
(as constructed for partial circuits in {\rm Theorem \ref{circ2PartLepM}}).
If $E_w(x)$ is undefined then $C_{w,x}$ and its depth are undefined.

\begin{defn} \label{DEFdepthM21wx} {\bf (depth of a word in 
${\cal RM}^{\sf fin}$).}
 
\smallskip

\noindent Let $\Gamma_{\rm\!RM}$ be a finite set such that 
$\Gamma_{\rm\!RM} \cup \tau$ generates ${\cal RM}^{\sf fin}$ and hence 
$M_{2,1}$. Let $\,w \in$ $(\Gamma_{\rm\!RM} \cup \tau)^*$, and let
$\,\ell_{\rm in}(w) = m\,$ be the input-length of $w$ ({\rm Def.\
\ref{DEF_IOlengthRMfin}}). For any $x \in {\rm Dom}(E_w) \cap \{0,1\}^m$, let
$\,C_{w,x}\,$ be the partial circuit just constructed for $w$ on input $x$.
Then the {\em depth} of $w$ is defined by

\medskip

 \ \ \  \ \ \ ${\sf depth}(w) \,=\,$
$\max \big\{{\sf depth}(C_{w,x}) \,:\,$
$x \in {\rm Dom}(E_w) \cap \{0,1\}^m \big\}$.
\end{defn}
If $E_w = \theta\,$ (the empty function) then ${\sf depth}(w)$ is not
defined.

\subsection{Unambiguous union}

\begin{defn} \label{DEFunamunion} {\bf (unambiguous union).}

\smallskip

\noindent {\small \bf (1) Sets:} For an indexed family of sets
$\,(D_i : i \in I)\,$ the {\em unambiguous union}, denoted by
 \ {\large $!$}$\bigcup_{i \in I} D_i\,$,
is the set of elements that belong to {\em exactly one} of the sets $D_i:$

\medskip

 \ \ \  \ \ \  {\large $!$}$\bigcup_{i \in I} D_i$  $ \ = \ $
$\{d \in \bigcup_{i \in I} D_i \,:\, (\exists!\,j \in I)[\,d \in D_j\,]\,\}$.

\medskip

\noindent In particular, if $\,D_i \cap D_j = \0\,$ for all $i, j \in I$ 
with $i \ne j$, then
 \ {\large $!$}$\bigcup_{i \in I} D_i \,=\, \bigcup_{i \in I} D_i$.

\bigskip

\noindent {\small \bf (2) Functions:} Consider a family of functions
$\,(f_i : i \in I)$, where $f_i = (X, r_i, X)\,$ (see Def.\ {\rm
\ref{DEFfunction}}); so all $f_i$ have the same source and target set $X$.
The {\em unambiguous union} of that family of functions is denoted by

\hspace{3,3cm} {\large $!$}$\bigcup_{i \in I} f_i$

\medskip

\noindent and defined to be $\,(X, r, X)$ where $r$ is the restriction of 
the relation $\,\bigcup_{i \in I} r_i\,$ to the unambiguous union
$\,${\large $\,!$}$\bigcup_{i \in I} {\rm Dom}(r_i)$.

\smallskip

In particular, if all the functions $f_i$ have two-by-two disjoint domains, 
then {\large $\,!$}$\bigcup_{i \in I} f_i$ $\,=\,$ $\bigcup_{i \in I} f_i$.
\end{defn}

\noindent The empty set is a neutral element for unambiguous union:
If $D_{j_0} = \0$ then \ {\large $!$}$\bigcup_{i \in I} D_i$  $\,=\,$
{\large $!$}$\bigcup_{i \in I \minus \{j_0\}} D_i$.

\smallskip

\noindent For unambiguous union the order of the sets or functions in a 
family does not matter, i.e., the unambiguous union operation satisfies 
generalized commutativity. 

\bigskip

\noindent {\bf Associativity, and relation between unambiguous union and 
symmetric difference:}

\smallskip

\noindent For any {\em two} sets $D_1, D_2:$ $ \ D_1$$\,!$$\cup\,$$D_2$  
$=$ $D_1\!\!\vartriangle\!\!D_2$ \ (symmetric difference).  
The symmetric difference is associative.
Since $\vartriangle$ and $\,!$$\cup$ are the same when applied to two sets,
the unambiguous union, as a binary operation, is {\em associative}.

\smallskip

\noindent The symmetric difference of a finite family of sets
$(D_1, \,\ldots\,, D_N)$ satisfies

\smallskip

 \ \ \  \ \ \ {\Large $\Delta$}$_{i=1}^{^N} D_i$ $ \ = \ $
$\{x \in \bigcup_{i=1}^{^N} D_i \,:\, x$ occurs in $D_i$ for an {\em odd}
number of indices $i\}$

\smallskip

 \ \ \  \ \ \ $=\,$
$( \ldots ((D_1 \vartriangle D_2) \vartriangle D_3) \vartriangle $
$\ldots$ $\vartriangle D_{n-1}) \vartriangle D_n$.

\medskip

\noindent However, for unambiguous union the situation is different. As a
variable-arity operation,

\smallskip

 \ \ \  \ \ \ {\large $!$}$\bigcup_{i \in I} D_i$
$ \ = \ \{x \in \bigcup_{i=1}^{^N} D_i \,:\, x$ occurs in $D_i$ for exactly
{\em one} index $i\}$.

\smallskip

\noindent This is usually different from
$ \ (\ldots ((D_1 \,!$$\cup\, D_2) \,!$$\cup\, D_3) \,!$$\cup\,$
$\ldots$ $\,!$$\cup\,D_{n-1}) \,!$$\cup\, D_n\,$,
since the latter is the same as the symmetric difference. So,
non-binary unambiguous union (of a family) is {\em not associative}. 

\smallskip

For an infinite family the unambiguous union
 \ {\large $!$}$\bigcup_{i \in I}D_i$ \ is well-defined but the symmetric
difference is not defined.
 \ \ \  \ \ \ [End, Remark.]

\smallskip

\begin{pro} \label{PROPunambM}\!\!\!.

\smallskip

\noindent {\small \rm (1)} \ For any family $(f_1,\ldots, f_k)$ with
$f_i \in M_{2,1}$ (or $\in {\cal RM}^{\sf fin}$) for all $i \in [1,k]$: 

\medskip

 \ \ \  \ \ \ {\large $!$}$\bigcup_{i=1}^k f_i \ \in \, M_{2,1}$ 
 \ \ (respectively
$\, \in {\cal RM}^{\sf fin}$).

\bigskip

\noindent {\small \rm (2)} \ Consider any  \ $f_i \in$ 
${\rm plep}M_{2,1} \minus \{\theta\}$ \ (or $\in$ 
${\sf pfl}{\cal RM}^{\sf fin} \minus \{\theta\}$), \ for $i \in [1,n]$.

 \ If $\,\delta(f_i)$ is the same for every $i$, then

\medskip

 \ \ \  \ \ \ {\large $ \ !$}$\bigcup_{i=1}^k  f_i$
 \ $\in \ {\sf plep}M_{2,1}$ 
 \ \ \ (respectively $\,\in$ ${\sf pfl}{\cal RM}^{\sf fin}$).
\end{pro}
{\sc Proof.} (1) follows straightforwardly from the definition of unambiguous
union.
(2) follows easily from the definitions of $\,!$$\cup$, $\,\delta(.)$,
$\,{\sf plep}M_{2,1}$, and ${\sf pfl}{\cal RM}^{\sf fin}$.
 \ \ \ $\Box$

\bigskip

\noindent The input-output length difference function $\delta(.)\,$ (see 
Def.\ \ref{DEFdelta} for total functions, and  Prop.\ 
\ref{PROPiolDiffpartial} for partial functions) can be generalized to 
$\,{\cal RM}^{\sf fin}$, $\,M_{2,1}$, and $G_{2,1}$:

\begin{defn} \label{DEFdeltaM} {\bf (general input-output length difference
{\boldmath $\delta_M(.)$}).} 

\smallskip

\noindent For $f \in {\cal RM}^{\sf fin}\,$ we define

\medskip

\hspace{2,5cm} $\delta_M(f)$  $\,=\,$
$\{\,|f(x)| - |x| \,: \ x \in {\rm domC}(f)\}$.
\end{defn}
Then $\,\delta_M$ is a total function from ${\cal RM}^{\sf fin}$ onto
${\cal P}_{\rm fin}({\mathbb Z})\,$
(the set of all finite subsets of ${\mathbb Z}$).

For $f \in {\sf pfl}{\cal RM}^{\sf fin}$), we have:
 \ $\delta_M(f) = \{\delta(f)\}$.
In particular, for the empty function $\theta$ we have
$\,\delta_M(\theta) = \0$; \ so $\delta_M(\theta)$ is defined,
while $\delta(\theta)$ is undefined.

\begin{pro} \label{PROPdeltaMWelldef}\!\!\!.

\smallskip

\noindent {\small \rm (1)} \ If $\,g, f \in {\cal RM}^{\sf fin}$ satisfy
$\,g \equiv_{\rm fin} f$, then $ \ \delta_M(g) = \delta_M(f)$.
 \ Hence $\delta_M$ is well-defined on $M_{2,1}$.

\medskip

\noindent {\small \rm (2)} \ For all $g, f \in M_{2,1}$:

\medskip

\hspace{2,5cm} $\delta_M(g \circ f) \ \subseteq \ \delta_M(g) + \delta_M(f)$

\medskip

where $+$ denotes elementwise addition in
$\,{\cal P}_{\rm fin}({\mathbb Z})$.

\medskip

\noindent {\small \rm (3)} \ Let $\Gamma$ be a finite set such that
$\Gamma \cup \tau$ generates ${\cal RM}^{\sf fin}$. Then for every word
$\,w \in$ $(\Gamma \cup \tau)^*$, 

\medskip

\hspace{2,5cm} $\delta_M(E_w) \ \subseteq \ $
$[-c_{_{(-)}}\,|w|_{\Gamma}, \ c^{^{(+)}}\,|w|_{\Gamma}]\,$, 

\medskip

where $\,c^{^{(+)}} \,=\,$
$\max \{\,j :\, j \ge 0,\,j \in \delta_M(\gamma),\,\gamma \in \Gamma\}$, 
and $\,c_{_{(-)}} \,=\,$
$\max \{\,|j| :\, j \le 0,\, j \in \delta_M(\gamma),\,\gamma \in \Gamma\}$.
\end{pro}
{\sc Proof.} Both (1) and (2) are straightforward.
Elementwise addition of $X, Y \in {\cal P}_{\rm fin}({\mathbb Z})$
is defined by $ \ X + Y = \{x+y : x \in X, y \in Y\}$; hence, if
$X = \0$ or $Y = \0$ then $X + Y = \0$.

\noindent (3) This follows from (2), using the fact that 
$\,\delta_M(\tau_{i,j}) = \{0\}\,$ for all $\tau_{i,j} \in \tau$. 
 \ \ \  $\Box$

\bigskip

\noindent The inclusion ``$\subseteq$'' in (2) can be strict, e.g.\ if
$g \circ f = \theta$, but $g \ne \theta \ne f$.
So the function $\delta_M(.)$ is not a semigroup homomorphism.

The invariance of $\delta(.)$ under $\,\equiv_{\rm fin}\,$ (Prop.\
\ref{PROPdeltaMWelldef}(1)), and the fact that $f$ is a right-ideal morphism,
imply that 

\smallskip

\hspace{1,75cm} $\delta_M(f)$ $\,=\,$
$\{|f(x)| - |x| : x \in {\rm domC}(f)\}$

\smallskip

\hspace{3,0cm} $=\,$     $\{|f(x)| - |x| : x \in {\rm Dom}(f)\}\,$.
 
\smallskip

\noindent This is a finite set (although ${\rm Dom}(f)$ is infinite), and
part (3) gives bounds on that set.

\smallskip

Since the monoids $M_{2,1}$, $\,{\rm tot}M_{2,1}$, and ${\rm plep}M_{2,1}$
are congruence-simple (see Subsection 7.1), $\delta_M(.)$ is not a semigroup
morphism on these monoids.

\bigskip 

Theorem \ref{THMwordforDisjU} will use the functions ${\sf C}(.)$ and 
${\sf W}(.)$ from Theorems \ref{circ2lepM} and \ref{circ2PartLepM}.
 \ Let $\Gamma_{\rm\!RM}$ be a finite generating set of 
${\cal RM}^{\sf fin}$ (and hence $M_{2,1}$), and let $\Gamma_{\sf\!pfl}$ be 
a finite set such that $\,\Gamma_{\sf\!pfl} \cup \tau\,$ generates 
${\sf pfl}{\cal RM}^{\sf fin}$ (and hence ${\sf plep}M_{2,1}$).

\medskip

\begin{thm} \label{THMwordforDisjU} {\bf (I. From an unambiguous union of
partial circuits to a $(\Gamma_{\rm\!RM} \cup \tau)$-word).}

\medskip

\noindent {\small \bf (1)} {\bf Unambiguous union of partial circuits with
the same {\boldmath $\delta$}-value}

\medskip

\noindent An unambiguous union of acyclic partial circuits with the same 
$\delta$-value is equivalent an acyclic partial circuit, or equivalently, to
a word in $\,(\Gamma_{\rm\!pfl} \cup \tau)^*$.
 \ In detail:

\medskip

\noindent Let $\,\{u_1, \,\ldots\,, u_k\}$  $\subseteq$
$(\Gamma_{\sf\!pfl} \cup \tau)^*\,$ be such that for all $i \ne j$ in
$[1,k]:$

\smallskip

 \ \ \  \ \ \ $\delta(E_{u_i}) = \delta(E_{u_j})$ \ \ and
 \ \ ${\rm Dom}(E_{u_i})$  $\cap$  ${\rm Dom}(E_{u_j}) = \0$.

\smallskip

\noindent Let $m$ and $n$ be the common input-length, respectively 
output-length, of $\,(E_{u_1},\,\ldots\,, E_{u_k})$.

\medskip

\noindent Then there exists a partial circuit $C$, and hence a word 
$\,w = {\sf W}(C) \in$ $(\Gamma_{\sf\!pfl} \cup \tau)^*\,$ with 
$E_w = E_C$, such that:

\medskip

\noindent $\,\bullet$ \ \ \  \ $!$$\bigcup (E_{u_1}, \,\ldots\,, E_{u_k})$
$ \ = \ $  $\bigcup_{i=1}^k E_{u_i}$
$\,=\, E_C \,\in\, {\sf pfl}{\cal RM}^{\sf fin}$.

\medskip

\noindent $\,\bullet$ \ \ \  \ $|C|$ $\,\le\,$    
$O\big(mk + nk + \sum_{i=1}^k |u_i|_{\Gamma_{\sf \!pfl} \cup \tau} \big)$.

\medskip

\noindent $\,\bullet$ \ \ \  \ $\ell_{\rm im}(C) = m\,$, 
 \ \ $\ell_{\rm out}(C) = n$. 

\medskip

\noindent $\,\bullet$ \ \ \  \ ${\sf depth}(C)$  $\,\le\,$
$O\big(\log k + \log n + \log m$ $\,+\,$
$\max\{ \max\{{\sf depth}(u_i), $
 $\,\log_2 |u_i|_{\Gamma_{\sf\!tfl} \cup \tau}\} \,:\, i \in [1,k] \} \big)$.

\medskip

\noindent $\,\bullet$ \ \ \  \ There is a deterministic algorithm
that computes $\,C\,$ on input $\,(u_1, \,\ldots\,, u_k)$, in log-space

\hspace{0.09in} (in terms of the input-length $\,\sum_{i=1}^k |u_i|$, where
any $\tau_{i,j} \in \tau$ is given the length $j$).

\bigskip

\noindent {\small \bf (2)} {\bf Unambiguous union of partial circuits with
different {\boldmath $\delta$}-values}

\medskip

\noindent An unambiguous union of partial circuits with different 
$\delta$-values is equivalent to a word in 
$\,(\Gamma_{\rm\!RM} \cup \tau)^*$.
 \ In detail:

\medskip

\noindent 
Let $\,\{u_1, \,\ldots\,, u_k\} \subseteq (\Gamma_{\sf\!pfl} \cup \tau)^*\,$
be such that for all $i \ne j$ in $[1,k]$:

\smallskip

 \ \ \  \ \ \  $\delta(E_{u_i}) \ne \delta(E_{u_j})$ \ \ \ (hence
$\,{\rm Dom}(E_{u_i}) \,\cap\, {\rm Dom}(E_{u_j}) = \0$).

\smallskip

\noindent Let $m$ be the common input-length of 
$\,(E_{u_1},\,\ldots\,, E_{u_k})$,
let $\,n_i = m + \delta(E_{u_i})\,$ (for  $i \in [1,k]$), and let
$ \ n_{\rm max} = \max\{n_i : i \in [1,k] \}$.

\medskip

\noindent Then there exists a word $\,w \in (\Gamma_{\rm\!RM} \cup \tau)^*$
such that

\medskip

\noindent $\,\bullet$ \ \ \  \ $!$$\bigcup (E_{u_1}, \,\ldots\,, E_{u_k})$
$ \ = \ $  $\bigcup_{i=1}^k E_{u_i}$
$\,=\, E_w \,\in\, {\cal RM}^{\sf fin}$.

\medskip

\noindent $\,\bullet$ \ \ \  \ $|w|_{\Gamma_{\rm\!RM} \cup \tau}$  $\,\le\,$
$\,O\big(\sum_{i=1}^k |u_i|_{\Gamma_{\sf\!pfl} \cup \tau} \big)\,$.

\medskip

\noindent $\,\bullet$ \ \ \  \ $\ell_{\rm in}(w) = m$.

\medskip

\noindent $\,\bullet$ \ \ \  \ ${\sf depth}(w)$  $\,\le\,$
$O(1) \,+\,  \max\big\{\max\{{\sf depth}(u_i), \,$
$\lceil \log_2 |u_i|_{\Gamma_{\sf\!tfl} \cup \tau} \rceil \}$ 
$ \,:\, i \in [1,k] \big\}$

\medskip

\noindent $\,\bullet$ \ \ \  \ There is a deterministic algorithm that
computes $\,w\,$ on input $\,(u_1, \,\ldots\,, u_k)$, in log-space

\hspace{0.09in} (in terms of the input-length $\,\sum_{i=1}^k |u_i|$, where
any $\tau_{i,j} \in \tau$ is given the length $j$).
\end{thm}
{\sc Proof.} (1) \ By Prop.\ \ref{PROPunambM}(2),
$ \ !$$\bigcup (E_{u_1}, \,\ldots\,, E_{u_k})$ $\,\in\,$
${\sf pfl}{\cal RM}^{\sf fin}$, so there exists $w \in$
$(\Gamma_{\sf\!pfl} \cup \tau)^*$ such that $\,E_w = $
$ \ !$$\bigcup (E_{u_1}, \,\ldots\,, E_{u_k})$.
We will construct a partial circuit for $w$; then Theorem \ref{circ2PartLepM}
immediately yields $w$ as a word over $\Gamma_{\sf\!pfl} \cup \tau$.
The partial circuit to be constructed has an input-output function
$\,x \in \{0,1\}^m$  $\,\longmapsto\,$  $E_w(x) \in \{0,1\}^n$, where $m$ is
the common input-length of $\,(E_{u_i} : i \in [1,k])$, and $n$ is the
common output-length.

The partial circuit is constructed so as to have four levels.

\medskip

\noindent {\sf Level 1:} 
The input is $x \in \{0,1\}^m$. The output is $x^k\,$ (i.e., $k$ copies of 
$x$).  It is produced by $\,m$ binary trees, each of depth 
$\,\lceil \log_2 k \rceil$, and each consisting of $\,k\!-\!1\,$ 
{\sf fork}-gates.
It is straightforward to construct the boolean circuit of level 1 in space 
$O(\log(km))$. 

\medskip

\noindent {\sf Level 2:} A preliminary remark: If we would apply
$\,(E_{u_1},\,E_{u_2}, \,\ldots\,, E_{u_k})\,$ to $\,(x, \,\ldots\,, x)$
$\,(= x^k)$, the result would be the concatenation
$\,E_{u_1}(x)\, E_{u_2}(x)\,\ldots\, E_{u_k}(x)$, which is undefined for all
$x\,$ (by disjointness of the domains).
So instead, we first use the completion $\ov{E}_{u_i}$ of each $E_{u_i}$.

Level 2 takes input $x^k$, and consists of the
$k$ circuits $\,\ov{C}_{u_1}, \,\ldots\,, \ov{C}_{u_k}$ in parallel, where
$\ov{C}_{u_i}$ is a circuit for the completion of ${\sf C}(u_i)\,$
(constructed in Prop.\ \ref{PROPcompletCircuit}). Each $\ov{C}_{u_i}$ takes
one copy of $1x$ as input (so, for this, a 1 is prepended to each of the $k$ 
copies of $x$, i.e., $(1x)^k$ is produced).
Each $\ov{C}_{u_i}$ has the same output-length, namely $n\!+\!1$.
Recall that $\,\ov{E}_{u_i}(1x)$  $=$  $1\,E_{u_i}(x)\,$ if $\,E_{u_i}(x)\,$
is defined, and $\,= 0^{n+1}\,$ otherwise.
The output of level 2 is the concatenation $ \ \ov{E}_{u_1}(1x)$
$\,\ldots\,$ $\ov{E}_{u_k}(1x)\,$ of length $\,k\,(n\!+\!1)$.
So, the boolean circuit for level 2 has input-output function

\smallskip

 \ \ \ $x^k \longmapsto$ $\, \ov{E}_{u_1}(1x) \ \ov{E}_{u_2}(1x)$ 
$ \ \ldots \ $   $\ov{E}_{u_{k-1}}(1x) \ \ov{E}_{u_k}(1x)\,$.

\smallskip 

\noindent The size of the boolean circuit for level 2 is

\smallskip

 \ \ \ $\le \ k \,+\, \sum_{i=1}^k |\ov{u}_i|_{\Gamma_{\sf\!pfl} \cup \tau}$
$ \ < \ c\, \sum_{i=1}^k |u_i|_{\Gamma_{\sf\!pfl} \cup \tau}\,$,

\smallskip

\noindent for some constant $c > 1$. The $+k$ is for producing the $k$ 1s
for $(1x)^k$, and the size of the boolean circuit $\ov{C}_{u_i}$ is
$\,O(|u_i|_{\Gamma_{\sf\!tfl} \cup \tau})$, by Prop.\ 
\ref{PROPcompletCircuit}.

\smallskip

\noindent The depth of this boolean circuit is 

\smallskip

 \ \ \ $\, \le\,$ $O(1) \,+\, $
$\max\{{\sf depth}(u_i), \,$
$\log_2 |u_i|_{\Gamma_{\sf\!tfl} \cup \tau} + 2\,\log_2 (n\!+\!1) \,:\,$
$i \in [1,k] \}$,

\smallskip

\noindent by Prop.\ \ref{PROPcompletCircuit}. And $\ov{C}_{u_i}$ and 
$\ov{u}_i$ can be computed from $u_i$ in log-space, by Prop.\ 
\ref{PROPcompletCircuit} and Theorem \ref{circ2lepM}(1).  

\medskip

\noindent {\sf Level 3:} This level has input $ \ \ov{E}_{u_1}(1x)$ 
$\,\ldots\,$ $\ov{E}_{u_k}(1x)$, and its output is the completion of 
$\,$!$\cup (E_{u_1},$  $\,\ldots$ $\,,$  $E_{u_k})(x)$. \ Equivalently, the 
output is $\,\ov{E}_{u_j}(1x)$, where $j \in [1,k]$ is the unique index for 
which $\,E_{u_j}(x)\,$ is defined; \ otherwise the output is $\,0^{n+1}$, if
$\,E_{u_i}(x)\,$ is undefined for all $i \in [1,k]$.

To obtain this output, we form the bitwise {\sf or} of the $k$ bitstrings 
$\,\ov{E}_{u_i}(1x)$, for $i \in [1,k]$. 
Let $y_j^{(i)}$ be the $j$th output bit of $\ov{E}_{u_i}(1x)$, for
$j \in [0,n]$, $\,i \in [1,k]$ \ (using the fact that all $\,E_{u_i}(x)\,$ 
are $\,0^{n+1}$, except perhaps one).  Then the $j$th output bit is

\smallskip

 \ \ \  \ \ \ $y_j = \ ${\large \sf OR}$_{i=1}^k\,y_j^{(i)}$
 \ \ \ for every $j \in [0,n]$.

\smallskip

\noindent Thus, level 3 consists of $n\!+\!1$ {\sf or}-trees, each of
size $k\!-\!1$ and depth $\lceil \log_2 k \rceil$. The $j$th {\sf or}-tree 
(for $j \in [0,n]$) has $k$ input bits and one output bit, namely $y_j$.  
The total output of level 3 is the concatenation $\,y_0 y_1 \,\ldots\, y_n$.

The boolean circuit for level 3 has size $(n\!+\!1)\,(k\!-\!1)$ and depth 
$\lceil \log_2 k \rceil$.

\smallskip

\noindent Levels 1, 2, and 3 together form an acyclic boolean (total) 
circuit.

\medskip

\noindent {\sf Level 4:} This level turns $\,y_0 y_1\,\ldots\,y_n$, which
is the output of completion $\,!$$\cup (E_{u_1},\,\ldots\,,$  $E_{u_k})(x)$,
 \ into the final partial output $\,!$$\cup (E_{u_1},\,\ldots\,,$  
$E_{u_k})(x)\,)\,$. This is equal to $\,y_1 \,\ldots\, y_n\,$ if $y_0 = 1$, 
and is undefined if $y_0 = 0$.

We build a partial circuit for level 4 as follows.  We use the partial 
one-bit gate $h_0$ such that $h_0(1) = 1$, and $h_0(0)$ is undefined. 
First, a partial $h_0$-gate is connected to the leftmost wire (that
carries the bit $y_0$).
Second, $n$ copies of the output-wire of the $h_0$-gate are made (by a binary 
tree of depth $\lceil \log_2 n \rceil$, made of $n\!-\!1$ {\sf fork}-gates); 
then we form 

\smallskip

 \ \ \  \ \ \ $(h_0(y_0))^n \ {\sf and}_n \ y_1\,\ldots\,y_n$ 

\smallskip

\noindent (the bitwise {\sf and} of two bitstrings of length $n$).
Thus the result is undefined if $y_0 = 0$; and it is
$\,y_1 \,\ldots \, y_j \,\ldots \, y_n\,$ if $y_0 = 1$.

The size the partial circuit of level 4 is $\,n\!-\!1 + n$, and the 
depth is $\lceil \log_2 n \rceil + 1$.  
 \ \ \  \ \ \  \ \ \ {\small [End, Level 4]}

\medskip

This completes the construction of a partial circuit for
$\,!$$\cup (E_{u_1}, \,\ldots\,, E_{u_k})$.  
As mentioned at the beginning of the proof, an equivalent word 
$w = {\sf W}(C) \in (\Gamma_{\sf\!pfl} \cup \tau)^*$ can
now be obtained by Theorem \ref{circ2PartLepM}.

By adding up the sizes and the depths of the four levels we obtain that the
total size of the partial circuit is 

\smallskip

 \ \ \  \ \ \ $\le\,$
$O\big(mk + nk + \sum_{i=1}^k |u_i|_{\Gamma_{\sf \!pfl} \cup \tau} \big)$;

\smallskip

\noindent and the depth is 

 \ \ \  \ \ \ $\le\,  O\big(\log k + \log n + \log m + \,$
$\max\{ \max\{{\sf depth}(u_i), $
 $\,\log_2 |u_i|_{\Gamma_{\sf\!tfl} \cup \tau}\} \,:\, i \in [1,k] \} \big)$.

\smallskip

\noindent Complexity of computing $w$: For levels 1, 3, and 4, the log-space 
computability is straightforward. For level 2, Prop.\ 
\ref{PROPcompletCircuit} yields log-space computability of $\ov{u_i}$.

\bigskip

\noindent (2) \ Since $\,!$$\bigcup (E_{u_1}, \,\ldots\,, E_{u_k}) \in$
${\cal RM}^{\sf fin}$, by Prop.\ \ref{PROPunambM}(1), there exists a word
$w \in (\Gamma_{\rm\!RM} \cup \tau)^*$ such that $\,E_w =$
$\,!$$\bigcup (E_{u_1}, \,\ldots\,, E_{u_k})$.  We construct the word $w$ in 
four levels similar to the ones in part (1).
Let $m$ be the common input-length of $(E_{u_i} : i = 1, \ldots, k)$.
There is no common output-length. 

We first construct a (slightly generalized) partial circuit for 
$\,!$$\bigcup (E_{u_1}, \,\ldots\,, E_{u_k})$, and then a word in 
$(\Gamma_{\rm\!RM} \cup \tau)^*$.

\medskip

\noindent {\sf Level 1:} The input is $x \in \{0,1\}^m$, the output is 
$x^k\,$ (i.e., $k$ copies of $x$), and a boolean circuit for this is built 
as in part (1). 

\medskip

\noindent {\sf Level 2:} The input is $x^k$. Level 2 consists of the $k$
boolean circuits $\,\ov{C}_{u_1}, \,\ldots\,, \ov{C}_{u_k}$ in parallel, 
just as in part (1), except that now $\ov{C}_{u_i}$ has output-length 
$n_i+1$, which is different for every $i \in [1,k]$.
The output of level 2 is the concatenation
$\,\ov{E}_{u_1}(1x) \ \ldots \ \ov{E}_{u_k}(1x)$, of length 
$\,k + \,\sum_{i=1}^k n_i$. 
By the same reasoning as in part (1) we find that the size of the boolean 
circuit $C^{(2)}$ of level 2 is 

\smallskip

 \ \ \  \ \ \ $|C^{(2)}| \,\le\,$
$k \,+\, \sum_{i=1}^k |\ov{u}_i|_{\Gamma_{\sf\!pfl} \cup \tau}$
$\,<\, O\big(\sum_{i=1}^k |u_i|_{\Gamma_{\sf\!pfl} \cup \tau}\big)\,$,

\smallskip

\noindent and the depth of level 2 is 

\smallskip

 \ \ \  \ \ \ ${\sf depth}(C^{(2)}) \,\le\,$ $O(1) \,+\, $
$\max\{{\sf depth}(u_i), \,$
$\log_2 |u_i|_{\Gamma_{\sf\!tfl} \cup \tau} + 2\,\log_2 (n_i+1) \,:\,$
$i \in [1,k] \}$,

\medskip

\noindent {\sf Level 3:} This level has input 
$\,\ov{E}_{u_1}(1x) \ \ldots \ \ov{E}_{u_k}(1x)$, of length
$\,k + \,\sum_{i=1}^k n_i\,$. Its output is $\,\ov{E}_{u_j}(1x)\,$
$(\,= 1\,E_{u_j}(x))$, where $j \in [1,k]$ is the unique index for which
$\,E_{u_j}(x)\,$ is defined; otherwise the output is $0\,$  (if 
$\,E_{u_i}(x)\,$ is undefined for all $i \in [1,k]$).

As in level 3 of part (1), let $y_0$ be the {\sf or} of the leading bits
of the $k$ strings $\ov{E}_{u_i}(1x)\,$ (for $i \in [1,k]$). Hence,
$\,y_0 = 1\,$ if $ \ !$$\bigcup (E_{u_1}, \,\ldots\,, E_{u_k})(x)\,$
is defined, and $y_0 = 0$ otherwise.
As in level 3 of (1), we can use an {\sf or}-tree (of size $k\!-\!1$ and
depth $\lceil \log_2 k \rceil$) to compute $y_0$.

Since the strings $\ov{E}_{u_i}(1x)\,$ (for $i \in [1,k]$) all have 
different lengths, the output cannot be obtained by an {\sf or}-trees 
as in level 3 of part (1).
Instead, we first apply the transformation $\xi_i$ to each
$\ov{E}_{u_i}(1x)$, so the output of level 3 is 
 \ $y_0 \ \xi_1(\ov{C}_{u_1}(1x))$ $ \ \ldots \ $
 $\xi_k(\ov{C}_{u_k}(1x))$.
Each $\,\xi_i(.)$ (for $i \in[1,k]$) is the right-ideal morphism defined by  

\smallskip

 \ \ \  \ \ \ $\xi_i: \ z_0\,y \in \{0,1\}\,\{0,1\}^{n_i}$ $\,\longmapsto\,$
$\xi_i(z_0 y) \in \{\e\} \cup \{0,1\}^{n_i}\,$, \ where

\smallskip

 \ \ \  \ \ \ $\xi_i(0y) = \e\,$, \ \ and 

\smallskip

 \ \ \  \ \ \ $\xi_i(1y) = 1y\,$;

\smallskip

\noindent so $\,{\rm domC}(\xi_i) = \{0,1\}^{n_i + 1}$.
To express $\xi_i$ by a word in $(\Gamma_{\rm\!RM} \cup \tau)^*$
we introduce the right-ideal morphism

\smallskip

 \ \ \  \ \ \ $\zeta: \ z_0 z_1 \in \{0,1\}^2 \,\longmapsto\, \zeta(z_0 z_1)$
$\,\in\,$ $\{\e\} \cup \{0,1\}\,$, \ where

\smallskip

 \ \ \  \ \ \ $\zeta(0 z_1) = \e\,$, \ and 

\smallskip

 \ \ \  \ \ \ $\zeta(1 z_1) = z_1\,$;

\smallskip

\noindent so $\,{\rm domC}(\zeta) = \{0,1\}^2$. 

We can either include $\zeta$ into the generating set 
$\Gamma_{\rm\!RM}$, or express it by a fixed word in 
$(\Gamma_{\rm\!RM} \cup \tau)^*$.

\noindent Now for all $\,y = y_1 \,\ldots\,y_j\,\ldots\ y_{n_i}$, with 
$y_j \in$ $\{0,1\}\,$ for $j \in [1,n_i]$, we have   

\smallskip

 \ \ \  \ \ \ $\xi_i(z_0\,y_1\,\ldots\,y_j\,\ldots\,y_{n_i})$  $\,=\,$  
$\zeta(z_0 y_1) \ \ldots \ \zeta(z_0 y_j) \ \ldots \ \zeta(z_0 y_{n_i})$.

\smallskip

\noindent Indeed, if $z_0 = 0$ this expression becomes $\e$, and if 
$z_0 = 1$ the expression becomes $y$.  

Thus, $\xi_i(z_0\,y_1\,\ldots\,y_{n_i})$ can be computed by a tree of 
{\sf fork}-gates (of size $n_i - 1$ and depth $\lceil \log_2 n_i \rceil$),
that makes $n_i$ copies of $z_0$; this tree can be converted into a 
$(\Gamma_{\sf\!tfl} \cup \tau)$-word of size $O(n_i)$.  
The $n_i$ copies of $z_0$ are then combined with 
$y_1\,\ldots\,y_j\,\ldots\,y_{n_i}$ in a row of $n_i$ instances of 
$\zeta(.)$. This yields a word for $\xi_i$ over 
$\{\zeta,\,{\sf fork}\}$ of size $O(n_i)$. 

Hence the output $\,y_0 \ \xi_1(\ov{C}_{u_1}(1x))$ $ \ \ldots \ $
$\xi_k(\ov{C}_{u_k}(1x))\,$ of level 3 can be computed by a word $w^{(3)}$
$\in (\Gamma_{\rm\!RM} \cup \tau)^*$ of size 

\smallskip

 \ \ \  \ \ \ $|w^{(3)}|_{\Gamma_{\rm\!RM} \cup \tau}$
$\,\le\,$ $k-1$ $+$ $\sum_{i=1}^k O(n_i)$ 
$ \ \le \ O(\sum_{i=1}^k n_i)$,

\smallskip

\noindent where, $k\!-\!1$ is for computing $y_0$, and $O(n_i)$ is for 
$\xi_i$.  \ And the depth is  

\smallskip

 \ \ \  \ \ \ ${\sf depth}(w^{(3)}) \,\le\,$ 
$\lceil \log_2 k \rceil$ $+$ 
$\max\{\lceil \log_2 n_i \rceil : i \in [1,k]\}$ $\,+\,$ $O(1)$.

\medskip

\noindent {\sf Level 4:} This level has input $0$ if $\,E_{u_i}(x)\,$ is
undefined for all $i \in [1,k]$; or it is $\,\ov{E}_{u_j}(1x)$
$(= 1\,E_{u_j}(x))$, where $j \in [1,k]$ is the unique index for which
$\,E_{u_j}(x)\,$ is defined.
The circuit for level 4 consists of the single right-ideal morphism $e$, 
defined by ${\rm domC}(e) =$ $\{1\}$, $\,{\rm imC}(e) = \{\e\}$, and 

\smallskip

 \ \ \  \ \ \ $e(1) = \e\,$ 

\smallskip

\noindent with $e(0)$ undefined.

Then the output of level 4 is either $\,E_{u_j}(x)$, if $y_0 = 1$ and $j$ is
as given for the input of level 4; or the output is undefined (if $y_0 = 0$).

Hence the size and depth of level 4 are constants. 

\medskip

\noindent The total size of the word $w$ of part (2) is $ \ \le\,$
$O\big(\sum_{i=1}^k (|u_i|_{\Gamma_{\sf\!pfl} \cup \tau} + n_i) \big)$.

\noindent The total depth of $w$ is  

${\sf depth}(w) \,\le\,$
$O(1) \,+\,  \max\big\{\max\{{\sf depth}(u_i), \,$
$\lceil \log_2 |u_i|_{\Gamma_{\sf\!pfl} \cup \tau} \rceil$  $+$ 
$2\,\lceil \log_2 (n_i+1) \rceil\} \ : \ i \in [1,k] \big\}$. 

\noindent Finally, since $\,n_i \le $
$c \,|u_i|_{\Gamma_{\sf\!pfl} \cup \tau}\,$ by Prop.\ \ref{PROPdomCimCLen}, 
the term involving $n_i$ can be removed from the formulas for size and depth.

\smallskip

All the levels are easily computed in log-space. For levels 1 and
2, this is similar to part (1); level 3 has a very simple structure, and 
level 4 has just one gate.  
 \ \ \  \ \ \  \ \ \ $\Box$

\smallskip

\begin{defn} \label{DEFexpressCircUnion}
{\bf \ ($\,\{\cdot,\, !$$\cup\}$-expressions, their length, depth and 
evaluation).}

\medskip

\noindent Let $\Gamma$ be one of the finite sets $\Gamma_{\rm\!RM}$, 
$\Gamma_{\sf\!pfl}$, or $\Gamma_{\sf\!tfl}$, seen above.

\medskip

\noindent {\small \bf (1)} \ A {\em $\,\{\cdot, !$$\cup\}$-expression} over
$\Gamma \cup \tau$ is defined recursively as follows.

\smallskip

\noindent $ \ \bullet$ \ {\rm (Words)} 
Any word $w \in (\Gamma \cup \tau)^*$ is a $\{\cdot, !$$\cup\}$-expression.

\smallskip

\noindent $ \ \bullet$ \ {\rm (Concatenation)} If $\,W_1$ and $\,W_2$ are
$\{\cdot, !$$\cup\}$-expressions, then so is $\,W_2 \cdot W_1 \ $ (also 
denoted 

by $\,W_2\,W_1$).

\smallskip

\noindent $ \ \bullet$ \ {\rm (Unambiguous union)} If $(W_k,\,\ldots\,,W_1)$
is a sequence of $\{\cdot, !$$\cup\}$-expressions with $k \ge 2$, then

$\,!$$\cup\,(W_k,\,\ldots\,,W_1)$ is a $\{\cdot, !$$\cup\}$-expression.

\medskip

\noindent The {\em length} $\,|W|_{\Gamma \cup \tau}$ and the {\em depth} of a
$\{\cdot, !$$\cup\}$-expression $W$ are defined recursively:

\smallskip

\noindent $ \ \bullet$ \ For any word $w \in (\Gamma \cup \tau)^*$, the 
length $\,|w|_{\Gamma \cup \tau}$ and $\,{\sf depth}(w)$ are as in 
Definitions {\rm \ref{word_length_size}} and {\rm \ref{DEFdepthRM}}.

\smallskip

\noindent $ \ \bullet$ \ If $W_1$ and $W_2$ are
$\{\cdot, !$$\cup\}$-expressions, then

\smallskip
 
 \ \ \  \ \ \ $|W_2 \cdot W_1|_{\Gamma \cup \tau}$  $\,=\,$
$|W_1|_{\Gamma \cup \tau} \,+\, |W_2|_{\Gamma \cup \tau}\,$, 

\smallskip

 \ \ \  \ \ \ ${\sf depth}(W_2 \cdot W_1) \,=\,$ 
${\sf depth}(W_2) + {\sf depth}(W_1)\,$.

\smallskip

\noindent $ \ \bullet$ \ If $(W_k,\,\ldots\,,W_1)$ is a sequence of
$\{\cdot, !$$\cup\}$-expressions with $k \ge 2$, then

\smallskip

 \ \ \  \ \ \ $|!$$\cup\,(W_k,\,\ldots\,,W_1)|_{\Gamma \cup \tau}$
$\,=\,$  $\sum_{i=1}^k |W_i|_{\Gamma \cup \tau}\,$, 

\smallskip

 \ \ \  \ \ \ ${\sf depth}(!$$\cup\,(W_k,\,\ldots\,,W_1))$
$\,=\,$  $\max\{{\sf depth}(W_i) : i \in [1,k]\}\,$.

\medskip

\noindent {\small \bf (2)} \ The {\em evaluation function} $E_W$ of a
$\{\cdot, !$$\cup\}$-expression $W$ is defined recursively:

\smallskip

\noindent $ \ \bullet$ \ A word $w \in (\Gamma \cup \tau)^*$ has an
evaluation function $\,E_w \in M_{2,1}$ (or $\in {\cal RM}^{\sf fin}$),
obtained by composing 

the generators from right to left, as they appear in $w$. 

\smallskip

\noindent $ \ \bullet$ \ The concatenation $\,W_2 \cdot W_1\,$ of
$\{\cdot, !$$\cup\}$-expressions has the evaluation function
$\,E_{W_2} \circ E_{W_1}(.)\,$

(composition of functions).

\smallskip

\noindent $ \ \bullet$ \ The unambiguous union
$\,!$$\cup\,(W_k,\,\ldots\,,W_1)$ of $\,\{\cdot, !$$\cup\}$-expressions has
the evaluation function

$\,!$$\cup\,(E_{W_k},\,\ldots\,,E_{W_1})$ \ \ (unambiguous union of
functions).
\end{defn}

\smallskip

\begin{lem} \label{LEMemptyintersectDom} {\bf (distributivity of composition
over union with disjoint domains).}

\smallskip

\noindent {\small \rm (1)} Let $f, g, \alpha, \beta$ be functions
$X \to X$ for some set $X$, such that
$ \ {\rm Dom}(\alpha) \,\cap\, {\rm Dom}(\beta) = \0$. \ Then

\smallskip

\noindent {\small \rm (1.a)} \hspace{0.25cm}
${\rm Dom}(\alpha \circ f(.))$  $\,\cap\,$
${\rm Dom}(\beta \circ f(.))$ $=$  $\0$, 

\smallskip

\noindent {\small \rm (1.b)} \hspace{0.25cm}
${\rm Dom}(f \circ \alpha(.))$  $\,\cap\,$
${\rm Dom}(g \circ \beta (.))$ $=$  $\0$.

\medskip

\noindent {\small \rm (2)} Let $(f_i : i \in I)$ and
$(g_j : j \in J)$ be families of functions on $X$ with
$\,{\rm Dom}(f_{i_1}) \cap {\rm Dom}(f_{i_2})$  $=$ $\0\,$ for all
$i_1 \ne i_2$ in $I$, and
$\,{\rm Dom}(g_{j_1}) \cap {\rm Dom}(g_{j_2})$  $=$  $\0\,$ for all
$j_1 \ne j_2$ in $J$. \ Then:

\medskip

\noindent {\small \rm (2.a)} \ For all $i_1, i_2, j_1, j_2$ such that
$\,(j_1, i_1)$ $\ne$ $(j_2, i_2):$
 \ \ \  ${\rm Dom}(g_{j_1} \circ f_{i_1})$ $\cap$
${\rm Dom}(g_{j_2} \circ f_{i_2})$ $=$  $\0$.

\medskip

\noindent {\small \rm (2.b)} \ (Distributivity:)
 \ \ \ \big({\large $!$}$\bigcup_{j \in J} g_j$\big)
$\circ$   \big({\large $!$}$\bigcup_{i \in I} f_i$\big)
$ \ = \ $  {\large $!$}$\bigcup_{(i,j) \in I \x J} \,g_j \circ f_i$.
\end{lem}
{\sc Proof.} (1.a)  For any function $X \stackrel{F}{\to} X$,
$ \ {\rm Dom}(F) = F^{-1}(X)$.  Hence,
$\,{\rm Dom}(\alpha \circ f(.))$  $\,\cap\,$  ${\rm Dom}(\beta \circ f(.))$
$\,=\,$  $f^{-1} \alpha^{-1}(X)$  $\,\cap\,$  $f^{-1} \beta^{-1}(X)$
$\,=\,$  $f^{-1} \big(\alpha^{-1}(X)$  $\,\cap\,$  $\beta^{-1}(X) \big)$
$\,=\,$  $f^{-1}(\0) = \0$.

\smallskip

\noindent (1.b)  We have $\,{\rm Dom}(f \circ \alpha(.))$
$\subseteq {\rm Dom}(\alpha)$; indeed, if $f(\alpha(x))$ is defined, then
$\alpha(x)$ is defined too.  Hence,
$\,{\rm Dom}(f \circ \alpha(.))$  $\,\cap\,$ ${\rm Dom}(g \circ \beta (.))$
$\,\subseteq\,$  ${\rm Dom}(\alpha) \,\cap\, {\rm Dom}(\beta) = \0$.

\smallskip

\noindent (2) If $i_1 \ne i_2$ then (1.b) implies
$\,{\rm Dom}(g_{j_1} \circ f_{i_1})$ $\cap$
${\rm Dom}(g_{j_2} \circ f_{i_2})$ $=$  $\0\,$ for all $j_1, j_2$
(including if $j_1 = j_2$).
If $i_1 = i_2$ and $j_1 \ne j_2$ then (1.a) implies
$\,{\rm Dom}(g_{j_1} \circ f_{i_1})$ $\cap$
${\rm Dom}(g_{j_2} \circ f_{i_1})$ $=$ $\0$.

By disjointness of the domains, the distributivity of composition over
unambiguous union now follows from the distributivity of composition over
union.
 \ \ \  \ \ \ $\Box$

\smallskip

\begin{pro} \label{LEMcircCompUnambig} {\bf (partial circuits for the 
composition of unambiguous unions).}

\smallskip

\noindent Consider two finite sets

\smallskip

 \ \ \ $\{(u^{(1)},\,d_u^{(1)}), \ \ldots \ , (u^{(h)},\,d_u^{(h)})\}$,
$ \ \{(v^{(1)},\,d_v^{(1)}), \ \ldots \ , (v^{(k)},\,d_v^{(k)})\}$
$ \ \subseteq \ (\Gamma_{\sf\!pfl} \cup \tau)^* \x {\mathbb Z}$

\smallskip

\noindent satisfying the following assumptions:

\smallskip

for all $i \in [1,h]:$
 \ if $E_{u^{(i)}} \ne \theta$ then $\,d_u^{(i)} = \delta(u^{(i)})$;

\smallskip

for all $j \in [1,k]:$
 \ if $E_{v^{(j)}} \ne \theta$ then $\,d_v^{(j)} = \delta(v^{(j)})$;

\smallskip

for all $i_1 \ne i_2$ in $[1,h]:$
$ \ {\rm Dom}(E_{u_{i_1}}) \cap {\rm Dom}(E_{u_{i_2}}) = \0$;

\smallskip

for all $j_1 \ne j_2$ in $[1,k]:$
$ \ {\rm Dom}(E_{v_{j_1}}) \cap {\rm Dom}(E_{v_{j_2}}) = \0$.

\medskip

\noindent Let $m$ be the common input-length of
$ \ \{ v_j u_i \,:\, i \in [1,h], \, j \in [1,k]\}$;

\smallskip

\noindent let \ $n_b \,=\,$
$m + \max\{d_u^{(i)} + d_v^{(j)} :\,i \in [1,h],\, j \in [1,k]\}$;

\smallskip

\noindent let \ $n_a \,=\,$
$m + \min\{d_u^{(i)} + d_v^{(j)} :\,i \in [1,h],\, j \in [1,k]\}$.

\medskip

\noindent We restrict every $E_{v_j u_i}(.)\,$ to $\{0,1\}^{\ge m} \ $
(for $\,i \in [1,h],\, j \in [1,k]$). \ Then we have:

\bigskip

\noindent {\small \bf (1)} \textbf{\textit{(Decomposition according to 
    $\delta$-values).}}

\medskip

\noindent For every $r \in [n_a,n_b]$ there exists 
$\,w^{(r)} \in$ $(\Gamma_{\sf\!pfl} \cup \tau)^*$ such that

\medskip

\noindent {\small \rm (1.1)} \ \ $E_{w^{(r)}} \ = \ $  {\large $!$}$\bigcup$
$\big(E_{v_j} \circ E_{u_i}(.) \,:\, (i,j) \in [1,h] \x [1,k],$
$ \ \delta(E_{v_j u_i}) = r\!-\!m \big)$;

\medskip

 \  \ \ \ moreover, by {\rm Def.\ \ref{DEFexpressCircUnion}}, this is the
evaluation function of the expression

\smallskip

 \ \  \ \ \  \ \ \ {\large $!$}$\bigcup \big(v_j u_i \,:\,$ 
    $(i,j) \in [1,h] \x [1,k], \ \delta(E_{v_j u_i}) = r\!-\!m \big)$;

\medskip

\noindent {\small \rm (1.2)} \ \ $E_{w^{(r)}} \in$
${\sf pfl}{\cal RM}^{\sf fin}$;

\medskip

\noindent {\small \rm (1.3)} \ \ if $E_{w^{(r)}} \ne \theta\,$ then 
$E_{w^{(r)}}$ satisfies:

\medskip

 \ \  \ \ \  \ \ \ ${\rm imC}(E_{w^{(r)}}) \subseteq \{0,1\}^r$,
 \ \ ${\rm domC}(E_{w^{(r)}}) \subseteq \{0,1\}^m$, \ \ so 

\smallskip
 
 \ \  \ \ \  \ \ \ $\delta(E_{w^{(r)}}) = r\!-\!m$;

\medskip

\noindent {\small \rm (1.4)} \ \
$|w^{(r)}|_{\Gamma_{\sf pfl} \cup \tau}$  $ \ \le \ $
$O\big(|I| \, |J| \,m \, n$  $\,+\,$
$\sum (\,|v_j u_i|_{\Gamma_{\sf pfl} \cup \tau} \,:\,$
$(i,j) \in [1,h] \x [1,k], \ \delta(E_{v_j u_i}) = r\!-\!m)\,\big)$;

\medskip

\noindent {\small \rm (1.5)} \ \ ${\rm maxindex}_{\tau}(w)$  $\,\le\,$
$O\big(\max(\, \{|I|\,|J|\,m, \ |I|\,|J|\,n\}$

\smallskip

\hspace{5cm} $\cup$  
$\{{\rm maxindex}_{\tau}(v_j u_i) \,:\, (i,j) \in [1,h] \x [1,k],$
$ \ \delta(E_{v_j u_i}) = r\!-\!m \} )\, \big).$

\medskip

\noindent {\small \rm (1.6)} \ \ ${\sf depth}(w^{(r)})$   $\,\le\,$
$O\big(\log_2(hk)$      $\,+\,$

\smallskip

\hspace{1,5cm} $\max\{{\sf depth}(v_j) + {\sf depth}(u_i)$  $+$
  $\log_2 |v_ju_i|_{\Gamma_{\sf pfl}} \,:\,$
  $i \in [1,h], \, j \in [1,k], \ \delta(E_{v_j u_i}) = r\!-\!m \} \big)$;

\medskip

\noindent {\small \rm (1.7)} \ \ $w^{(r)}$ and $\,m\!-\!r\,$
can be computed from $\,\{u_1,\,\ldots\,,u_h\}$, $\{v_1,\,\ldots\,,v_k\}$

 \  \ \ \  in log-space.

\bigskip

\noindent {\small \bf (2)} \textbf{\textit{(Recombination into an 
unambiguous union).}}

\medskip

\noindent {\small \bf (2.1)} \ \  ${\rm Dom}(E_{w^{(r)}})$  $\cap$
${\rm Dom}(E_{w^{(s)}}) = \0$ \ for all $r \ne s$ in $[n_a,n_b]$.

\medskip

\noindent {\small \bf (2.2)} 
 \ \ {\large $!$}$\bigcup \big(E_{w^{(r)}} : r \in [n_a, n_b]\big) \ $
$(\in M_{2,1})\,$ is the evaluation function of the expression

\medskip

 \ \  \ \ \  \ \ \ {\large $!$}$\bigcup (v_j u_i \,:\,$ 
$(i,j) \in [1,h] \x [1,k])$
$ \ = \ $
\big({\large $!$}$\bigcup_{j=1}^k v_j\big)$
$\cdot$  \big({\large $!$}$\bigcup_{i=1}^h u_i \big)$;

\medskip

 \ \  \ \ \  \ \ \ and $\delta_M$ of this function is
$ \ \{ r\!-\!m \,:\, r \in [n_a, n_b], \  E_{w^{(r)}} \ne \theta\}$.

\medskip

 \ \ \  \ \ The set $ \ \{(w^{(r)},\,r\!-\!m) : r \in [n_a,n_b]\}\,$ can be
computed from $\,\{u_1,\,\ldots\,,u_h\}, \, \{v_1,\,\ldots\,,v_k\}$

 \ \ \  \ \ in log-space. 
\end{pro}
Moreover, a word $w \in (\Gamma_{\rm\!RM} \cup \tau)^*$ for the disjoint 
union of functions in (2.2) can be found in log-space, by Theorem \ref{THMwordforDisjU}(2).) 

\medskip

\noindent {\sc Proof.} (1.1), (1.2) and (1.3) follow from Theorem
\ref{THMwordforDisjU}(1), since for $E_{w^{(r)}}$ we only consider words
$v_ju_i$ with the same $\delta$-value (namely $r\!-\!m$).
(1.4) follows from the fact that for $x \in {\rm Dom}(E_{w^{(r)}})$ we
have $|E_{w^{(r)}}(x)| - |x| = r$, while for $x \in {\rm Dom}(E_{w^{(s)}}$
we have $|E_{w^{(s)}}(x)| - |x| = s \ne r$.
(1.5), (1.6), (1.7), and (1.8) follow from Theorem \ref{THMwordforDisjU}(1),
applied to the unambiguous union
 \ {\large $!$}$\bigcup \big(v_j u_i \,:\, (i,j) \in [1,h] \x [1,k],$
$\, m + \delta(E_{v_j u_i}) = r \big)$.
 \ (2) follows immediately from the properties of the functions $E_{w^{(r)}}$  
proved in (1).
 \ \ \  \ \ \ $\Box$

\bigskip

\noindent 
So far we have proved the following: (1) The unambiguous union of a finite
family of functions in $M_{2,1}$ (or in ${\cal RM}^{\sf fin}$) is in
$M_{2,1}$ (respectively in ${\cal RM}^{\sf fin}$), by Prop.\
\ref{PROPunambM}. 
 \ (2) The unambiguous union of a finite family of functions in
${\rm plep}M_{2,1}$ (or ${\sf pfl}{\cal RM}^{\sf fin}$) with the same
$\delta$-value is in ${\rm plep}M_{2,1}$ (respectively
${\sf pfl}{\cal RM}^{\sf fin}$), and a circuit for it can be constructed
(Theorem \ref{THMwordforDisjU}).
If the functions have different $\delta$-values then their unambiguous union
is just in $M_{2,1}$ (or in ${\cal RM}^{\sf fin}$), and a word in 
$(\Gamma_{\rm\!RM} \cup \tau)^*$ can be constructed by Theorem 
\ref{THMwordforDisjU}.
 \ (3) The composition of two finite unambiguous unions of functions in
$M_{2,1}$ (or in ${\cal RM}^{\sf fin}$) is a function in $M_{2,1}$
(respectively in ${\cal RM}^{\sf fin}$); a finite set of ${\rm plep}M_{2,1}$
(or ${\sf pfl}{\cal RM}^{\sf fin}$) circuits for it can be constructed, one
for each $\delta$-value (Prop.\ \ref{LEMcircCompUnambig}).

\subsection{Decomposition into an unambiguous union of partial circuits}

The next Theorem shows that, conversely, any ``generalized circuit''
$w \in $ $(\Gamma \cup \tau)^*\,$ (where $\Gamma$ is finite and 
$\Gamma \cup \tau$ generates ${\cal RM}^{\sf fin}$, hence $M_{2,1}$) can be 
decomposed into an equivalent unambiguous union of partial circuits, one for 
each number in $\delta_{\!M}(w)$.

The proof of the Theorem uses the following finite generating set of 
$M_{2,1}\,$ (see \cite[{\small Section 3}]{BiThompsMon}):

\bigskip

\noindent $a_0 = \left( \hspace{-.08in} \begin{array}{l|l}
0 & 1  \\
1 & 0
\end{array} \hspace{-.08in} \right)$,
 \ \ \ $a_1 = \left( \hspace{-.08in} \begin{array}{r|r|r}
00 & 01 & 1  \\
00 &  1 & 01
\end{array} \hspace{-.08in} \right)$,
 \ \ \ $a_2 = \left( \hspace{-.08in} \begin{array}{r|r|r}
0 & 10 & 11  \\
10 & 0 & 11
\end{array} \hspace{-.08in} \right)$,
 \ \ \ $\tau_{1,2} = $
$\left( \hspace{-.08in} \begin{array}{r|r|r|r}
00 & 01 & 10 & 11  \\
00 & 10 & 01 & 11
\end{array} \hspace{-.08in} \right)$,

\smallskip

\noindent
$a_3 = (\varepsilon \to 0)$,
 \ \ \ $a_4 = (0 \to \varepsilon)$,
 \ \ \ $a_5 = (0 \to 00)$,
 \ \ \ $a_6 = (00 \to 0)$,
 \ \ \ $a_7 = \left( \hspace{-.08in} \begin{array}{l|l}
0 & 1  \\
0 & 0
\end{array} \hspace{-.08in} \right)$,
 \ and  \\
$a_8 = \left( \hspace{-.08in} \begin{array}{r|r}
0 & 1  \\
0 & 01
\end{array} \hspace{-.08in} \right)$ .

\bigskip

\noindent Let $\Gamma_{\!M}$ $=$ $\{a_i : i \in [0,8]\}$; then 
$\Gamma_{\!M} \cup \tau$ generates $M_{2,1}$.
Since we will use a generating set of the form $\Gamma \cup \tau$, the 
generator $\tau_{1,2}$ above will be considered to be only in $\tau$, 
and not in $\Gamma_{\!M}$.

The elements of $\Gamma_{\!M}$ are {\sf pfl} (plep with fixed input-length) 
except for $a_1$, $a_2$, and $a_8$.  The non-{\sf pfl} generators can be 
written as unambiguous unions of {\sf pfl} functions as follows:

\bigskip

\noindent
$a_1 = \left( \hspace{-.08in} \begin{array}{r|r|r}
00 & 01 & 1 \\
00 &  1 & 01 
\end{array} \hspace{-.08in} \right)$
$ \ = \ $
{\Large $!$$\bigcup$}$\, \big\{ (01 \to 1), \ (00 \to 00), $
$ \ (1 \to 01) \big\}$
$ \ = \ $ {\large $!$$\bigcup$}$\,\{ a_{1,-1}, \ a_{1,0}, \ a_{1,1} \}$\,;

\bigskip

\noindent the $a_{1,k}\,$ (for $k \in \{-1,0,1\} = \delta_M(a_1)\,$) have 
2-by-2 disjoint domains;

\bigskip

\noindent $a_2 = $
$\left( \hspace{-.08in} \begin{array}{r|r|r}
 0 & 10 & 11  \\
10 &  0 & 11
\end{array} \hspace{-.08in} \right)$
$ \ = \ $
{\Large $!$$\bigcup$}$\,\big\{\! (10 \to 0), \ (11 \to 11),
 \ (0 \to 10) \big\}$
$ \ = \ $ {\large $!$$\bigcup$}$\,\{ a_{2,-1}, \ a_{2,0}, \ a_{2,1} \}\,$;

\bigskip

\noindent the $a_{2,k}\,$ (for $k \in \{-1,0,1\} = \delta_M(a_2)\,$) have 
2-by-2 disjoint domains.

\bigskip

\noindent $a_8 = \left( \hspace{-.08in} \begin{array}{r|r}
0 & 1  \\
0 & 01
\end{array} \hspace{-.08in} \right)$
$ \,=\, $
{\large $!$}$\bigcup\, \{ (0 \to 0), \ (1 \to 01) \}$
$ \,=\, $
{\large $!$}$\bigcup\, \{ a_{8,0}, \ a_{8,1} \}\,$;

\medskip

\noindent moreover, the $a_{8,k}\,$ (for $k \in \{0,1\} = \delta_M(a_8)\,$)
have disjoint domains;

\bigskip

\noindent From $\Gamma_{\!M}$ we thus obtain the set

\medskip

 \ \ \ $\Gamma_{\!M,\,{\sf pfl}}$ $=$   
$\{a_0,\,a_{1,-1},\, a_{1,0},\, a_{1,1},\,$
$a_{2,-1},\, a_{2,0},\, a_{2,1},\,$
$a_3,\, a_4,\, a_5,\, a_6,\, a_7,\, a_{8,0},\, a_{8,1}\}\,$.

\medskip

\noindent So, $\Gamma_{\!M,\,{\sf pfl}} \subseteq $
${\sf pfl}{\cal RM}^{\sf fin}$; and $\Gamma_{\!M,\,{\sf pfl}}$ consists of 
the {\sf pfl} elements of $\Gamma_{\!M}$, together with the eight newly 
introduced {\sf pfl} fragments of the non-plep elements of $\Gamma_{\!M}$.  

\smallskip

The monoid generated by $\Gamma_{\!M,\,{\sf pfl}} \cup \tau$ in
${\cal RM}^{\sf fin}$ is a submonoid of ${\sf pfl}{\cal RM}^{\sf fin}$ that,
under the map of the congruence $\equiv_{\rm fin}$, maps onto 
${\rm plep}M_{2,1}$.  If we want to generate all of
${\sf pfl}{\cal RM}^{\sf fin}$ we can add the gates $\,\{ {\sf and}, 
{\sf or}, {\sf fork}\}\,$ to $\Gamma_{\!M,\,{\sf pfl}}\,$ (the gate 
$\,{\sf not}\,$ is already included).

\medskip

In Theorem \ref{THM_Mtoplep}(2), $\,{\sf fLog\,unif\,AC}_2$ denotes the 
parallel complexity class ${\sf AC}_2$ with log-space uniformity. 

\smallskip 

\begin{thm} \label{THM_Mtoplep} {\bf (II. Decomposing functions in
${\cal RM}^{\sf fin}$ or $M_{2,1}$ into an 

\hspace{3,2cm} unambiguous union of partial circuits).}

\medskip

\noindent {\small \bf (1)} \textbf{\textit{Function form of the Theorem:}}

\smallskip

For any $f \in M_{2,1}$, let $\,m = {\rm maxlen}\big({\rm domC}(f)\big)$. 
And let
$\,n_a = {\rm minlen}\big(f(\{0,1\}^m)\big)\,$ and $\,n_b$ $=$  
${\rm maxlen}\big(f(\{0,1\}^m)\big)$. Then there exist functions
$\,f^{(n_a)}, \,\ldots\,, f^{(n_b)}$  $\in$  ${\rm plep}M_{2,1}$ such that:

\medskip

\noindent $\,\bullet$ \hspace{0.25cm} $f \,=$
$ \ \bigcup_{i=n_a}^{ \ \ n_b} f^{(i)}$
$ \ = \ $  {\large $!$}$\bigcup_{i=n_a}^{ \ \ n_b} f^{(i)}$.

\medskip

\noindent $\,\bullet$ \hspace{0.25cm} For all $i \ne j:$
 \ \ ${\rm Dom}(f^{(i)}) \cap {\rm Dom}(f^{(j)}) = \0$.

\medskip

\noindent $\,\bullet$ \hspace{0.25cm} For all $i \in [n_a, n_b]$: 
$ \ f^{(i)}$ is the restriction of $f$ such that 

\smallskip

\hspace{0.75cm} ${\rm domC}(f^{(i)})$  $\,=\,$ 
$\{0,1\}^m \,\cap\, f^{-1}(\{0,1\}^i)$;   \ \ hence

\smallskip

\hspace{0.75cm} $f^{(i)}(\{0,1\}^m) \,=\, {\rm im}C(f^{(i)})$  
$ \ \subseteq \ $  $\{0,1\}^i$, 

\smallskip

\hspace{0.75cm} $\delta(f^{(i)}) = i\!-\!m$ \ \ if
  \ \ $f^{(i)} \ne \theta$.

\bigskip

\noindent {\small \bf (2)} \textbf{\textit{Circuit form of the Theorem:}}

\smallskip

Let $\,\Gamma_{\!M} = \{a_i : i \in [0,8]\}$ and $\,\Gamma_{\!M, {\sf pfl}}$
$\subseteq \Gamma_{\sf\!pfl}\,$ be as described before the Theorem.

\smallskip

For any $\,w \in (\Gamma_{\!M} \cup \tau)^*$, let $\,m = \ell_{\rm in}(w)$.
Then there exist $\,n_b \ge n_a \ge 0$, and there exist partial circuits
$\,W^{(n_a)},$  $\,\ldots\,,$  $W^{(n_b)}$  $\,\in\,$
$(\Gamma_{\sf\!pfl} \cup \tau)^*\,$ such that the following holds:

\medskip

\noindent {\small \bf (2.1)} \ Decomposition of $w$ into partial circuits:

\medskip

\noindent $\,\bullet$ \hspace{0.25cm}
$E_w \,= \ \bigcup_{i=n_a}^{ \ n_b} E_{W^{(i)}}$
$ \ = \ $  {\large $!$}$\bigcup_{i=n_a}^{ \ n_b} E_{W^{(i)}}$;
 \ so $E_w$ is the evaluation function of
$\,$ {\large $!$}$\bigcup_{i=n_a}^{ \ n_b} \!W^{(i)}$.

\medskip

\noindent $\,\bullet$ \hspace{0.25cm} For all $i \ne j:$
 \ \ ${\rm Dom}(E_{W^{(i)}}) \,\cap\, {\rm Dom}(E_{W^{(j)}}) = \0$.

\medskip

\noindent $\,\bullet$ \hspace{0.25cm} $\,m \,\le\,$
${\rm maxlen}(\Gamma_{\!M}) \cdot |w|_{\Gamma_{\!M}}$  $\,+\,$
${\rm maxindex}_{\tau}(w)\,$;

\smallskip

\hspace{0,13cm} $n_a \,\le\, n_b \,\le\, $
$2 \ {\rm maxlen}(\Gamma_{\!M}) \cdot |w|_{\Gamma_{\!M} \cup \tau}$
$\,+\,$  ${\rm maxindex}_{\tau}(w)$.

\bigskip

\noindent {\small \bf (2.2)} \ For all $i \in [n_a, n_b]:$

\smallskip

\noindent $\,\bullet$ \hspace{0.25cm}
$E_{W^{(i)}}$ is the restriction of $E_w$ such that

\smallskip
 
 \ \ \  \ \ \  \ \ \ ${\rm domC}(E_{W^{(i)}})$  $\,=\,$
$\,\{0,1\}^m \,\cap\, E_w^{ \ -1}(\{0,1\}^i) \ \subseteq \ \{0,1\}^m$.

\smallskip

\hspace{0,11cm} Hence $E_{W^{(i)}}$ (but not the word $W^{(i)}$) is uniquely
determined by $E_w$; 

\smallskip

\hspace{0,11cm} $E_{W^{(i)}}(\{0,1\}^m) \,=\, {\rm im}C(E_{W^{(i)}})$
$ \ \subseteq \ $  $\{0,1\}^i$,

\smallskip

\hspace{0,11cm} $\delta(E_{W^{(i)}}) = i\!-\!m$ \ \ if
 \ \ $E_{W^{(i)}} \ne \theta$.

\medskip

\noindent $\,\bullet$ \hspace{0.25cm}
$|W^{(i)}|_{\Gamma_{\!M,\,{\sf pfl}} \cup \tau}$  $\,\le\,$
${\rm maxlen}(\Gamma_{\!M}) \, |w|_{\Gamma_{\!M} \cup \tau}$  $\,+\,$
${\rm maxindex}_{\tau}(w)\,$;

\smallskip

\hspace{0,11cm}
${\rm maxindex}_{\tau}(W^{(i)}) \,\le\, {\rm maxindex}_{\tau}(w)$.

\medskip

\noindent $\,\bullet$ \hspace{0.25cm} ${\sf depth}(w)$  $\,=\,$
$\max\{{\sf depth}(W^{(i)}) : i \in [n_a,n_b]\}$
 \ \ \ (where ${\sf depth}(w)$ was defined in {\rm Def.\ 
\ref{DEFdepthM21wx}}).

\medskip

\noindent $\,\bullet$ \hspace{0.25cm} On input $w$, the sequence
$\,(W^{(n_a)}, \,\ldots\,, W^{(n_b)})\,$ can be computed in
$\,{\sf fLog\,unif\,AC}_2$

\hspace{0.11cm} (where the length of the input $w$ is taken to be
$ \ {\rm maxlen}(\Gamma_{\!M}) \ |w|_{\Gamma_{\!M}}$
$+$ ${\rm maxindex}_{\tau}(w)$).

\bigskip

\noindent {\small \bf (3)} \  The results of {\rm (2)} also hold with
$M_{2,1}$ replaced by $\,{\cal RM}^{\sf fin}$, and $\,{\rm plep}M_{2,1}$
replaced by $\,{\sf pfl}{\cal RM}^{\sf fin}$.
\end{thm}
{\sc Proof.} (1) We represent $f \in M_{2,1}$ by a maximally extended
right-ideal morphism in ${\cal RM}^{\sf fin}$, that we also call $f$.
Let $m = {\rm maxlen}({\rm domC}(f))$; we assume $m \ge 1$ (otherwise we
choose $m=1$).
Let $n_b = {\rm maxlen}(f(\{0,1\}^m))$, and let
$n_a = {\rm minlen}(f(\{0,1\}^m))$.

For every $i \in [n_a,n_b]$, let $f^{(i)}$ be the restriction of $f$ with
$\,{\rm domC}(f^{(i)})$  $\,=\,$
$\{0,1\}^m$  $\cap$  $f^{-1}(\{0,1\}^i)$.
Then $f^{(i)} \in {\rm plep}M_{2,1}$, and the other properties in part (1)
of the Theorem are easily verified.

Note that some of the functions $f^{(i)}$ could be the empty function
$\theta$.

\medskip

\noindent (2) Let $\,\Gamma_{\!M}$ and $\Gamma_{\!M,\,{\sf pfl}} (\subseteq$
$\Gamma_{\sf\!pfl})$ be as above. For any $w \in (\Gamma_{\!M} \cup \tau)^*$,
there exists $m \ge 1$ such that by Prop.\ \ref{PROPdomCimCLen}:

\smallskip

 \ \ \  $m \,\le\,$
${\rm maxlen}(\Gamma_{\!M}) \cdot |w|_{\Gamma_{\!M}}$  $\,+\,$
${\rm maxindex}_{\tau}(w)$,

\smallskip

 \ \ \ ${\rm domC}(E_w)$ is $\equiv_{\rm fin}$-equivalent to a
subset of $\{0,1\}^m$, \ \ and

\smallskip

 \ \ \ $E_w(\{0,1\}^m) \subseteq \{0,1\}^{\le n}$ \ \ for some $n \ge 1$ with

\smallskip

 \ \ \ $n \,\le\,$
$2 \ {\rm maxlen}(\Gamma_{\!M}) \cdot |w|_{\Gamma_{\!M}}$
$\,+\,$     ${\rm maxindex}_{\tau}(w)$.

\smallskip

\noindent We represent $E_w$ by a right-ideal morphism restricted to
$\,\{0,1\}^{\ge m}$, i.e., ${\rm domC}(E_w) \subseteq \{0,1\}^m$.
In ${\rm Dom}(E_w) \subseteq \{0,1\}^{\ge m}\,$ we have:
$\,E_w \equiv_{\rm fin}\,w_{\ell} \circ \,\ldots\, \circ w_1(.)$, where 
$w = w_{\ell} \,\ldots\, w_1$ with $w_t \in \Gamma_{\!M} \cup \tau$ for 
all $t \in [1, \ell]$. So, when $w_{\ell} \circ \ldots \circ w_1(x)$ is 
calculated for $|x| \ge m$, the functions $w_t$ are composed and they are 
defined without any need of a further restriction of right-ideal morphisms 
(see also  \cite[{\small Cor.\ 3.7}]{BiTh},
\cite[{\small Thm.\ 4.5}]{BiThompsMonV3},
\cite[{\small Prop.\ 6.5}]{BinMk1}, \cite[{\small Def.\ 2.4}]{BiEvalProb}).

Among the generators that appear in $\{w_{\ell}, \,\ldots\,, w_1\}$, some
are in ${\rm plep}M_{2,1}$,  namely those in
$\{a_0,$ $a_3,$ $a_4,$ $a_5,$ $a_6,$ $a_7\}$  $\cup$  $\tau$.
For the other generators, namely those in $\{a_1, a_2, a_8\}$, we use the
decompositions into unambiguous unions with disjoint domains:
$\,a_1 =$  $!$$\cup \{a_{1,-1},\, a_{1,0},\, a_{1,1}\}$ $\,( = a_1^!)$,
$\,a_2 =$  $!$$\cup \{a_{2,-1},\, a_{2,0},\, a_{2,1}\}$ $\,( = a_2^!)$, and
$\,a_8 =$  $!$$\cup \{a_{8,0},\, a_{8,1}\}$ $\,( = a_8^!)$.
For $w_t \in \{a_1, a_2, a_8\}$, we denote the relevant decomposition of 
$w_t$ into an unambiguous union by $w_t^!$; for $w_t \in {\rm plep}M_{2,1}$
we have  $\,w_t^! = w_t$.

\smallskip

Main idea of the proof: 
In the composition $\,E_w = w_{\ell} \circ \ldots \circ w_1(.)\,$ we
will replace each $w_t$ by $w_t^!$, for all $t \in [1,\ell]$. Then, by using
distributivity of composition over union in 
$\,w_{\ell}^! \circ \ldots \circ w_1^!(.)\,$ (Lemma 
\ref{LEMemptyintersectDom}) , we obtain an expression for $E_w$ in the form 
of an unambiguous union of {\sf pfl}-functions with disjoint domains.
If, moreover, we collect together the terms of equal $\delta$-value we obtain
$\,E_w \,= \ !$$\cup (E_{W^{(n_a)}}, \,\ldots\,, E_{W^{(n_b)}})$, where for
every $i \in [n_a,n_b]$: either $E_{W^{(i)}} = \theta\,$ or
$\delta(E_{W^{(i)}}) = i\!-\!m$.

The claimed properties of the words $W^{(i)}$ in the Theorem do not depend
on the order in which the distributivity steps are carried out in
$\,w_{\ell}^! \circ \ldots \circ w_1^!(.)$. 
However, the complexity of the algorithm for finding the words $W^{(i)}$
depends on this order.

\smallskip

We implement this idea in two algorithms.

In the first algorithm, called the \textbf{\textit{sequential algorithm}},
the generators and distributivity are applied from right to left (starting
with $w_1^!$, then $w_2^!\, w_1^!$, then $w_3^!\, w_2^!\, w_1^!$, etc.).
At each step, Theorem \ref{THMwordforDisjU} and Prop.\ 
\ref{LEMcircCompUnambig} are used in order to collect the terms with the
same $\delta$-value. This algorithm runs in polynomial time (but it is not
an {\sf NC}-algorithm).

\smallskip

In the second algorithm, called the \textbf{\textit{parallel algorithm}},
a binary tree of pairwise compositions is used.  So we compute the sequence

\smallskip

\hspace{2,3cm} $\ldots \ , \ \ \ldots \ ,$
$ \ \ (w_{2h}^! \circ w_{2h-1}^!), \ \ \ldots \ , \ $
$(w_4^! \circ w_3^!)$,  $(w_2^! \circ w_1^!)\,$;

\smallskip

\noindent then we again compose the results in pairs, and again in pairs, 
etc. In each step, distributivity and Prop.\ \ref{LEMcircCompUnambig} are 
applied. There are $\lceil \log_2 \ell \rceil$ parallel steps.
This can be carried out as a parallel algorithm, since all the pairs used in
the same step can be composed independently at the same time. We will see 
that this is an ${\sf AC}_2$-algorithm.

\smallskip

We next describe the two algorithms in detail. In the algorithms,
``$:=$'' denotes assignment, ``$=$'' is equality, and in a loop \ ``for
$k \in [a,b]: \ldots$'' \ the variable $\,k$ ranges over the integers in
$[a,b]$ in increasing order.

\bigskip

\noindent {\sf Sequential algorithm}:

\smallskip

\noindent {\footnotesize \sf
/* Input: \  $w = w_{\ell} \,\ldots\, w_1\,$ (with $w_t \in$
$\Gamma_{\!M} \cup \tau$ for $t \in [1,\ell]$);

\noindent \ * \ Output: \ $\{(W^{(i)},\,i\!-\!m) : i \in [n_a,\, n_b]\}$
$\,\subseteq\,$  $(\Gamma_{\!M,\,{\sf pfl}} \cup \tau)^* \x {\mathbb Z}$,
satisfying the conditions of the Theorem.
 \ \ */ }

\medskip

\noindent $m \,:=\,$
${\rm maxlen}(\Gamma_{\!M}) \cdot |w|_{\Gamma_{\!M} \cup \tau}$
$\,+\,$  ${\rm maxindex}_{\tau}(w)\,$;
\hspace{0.75cm} {\footnotesize \sf /* will not change */}

\smallskip

\noindent $n_b := m$; \ \ \ $n_a := m$; \ \ \ $W^{(n_a)} := \varepsilon$;
\hspace{1.5cm}  {\footnotesize \sf /* initialization, will change */}

\medskip

\noindent for $\,t \in [1,\, \ell]:$

\smallskip

{\footnotesize \sf
/* {\sl Loop assertion:} \ At the beginning of the loop (i.e., at value
   $t\!-\!1$), numbers $n_b \ge n_a \ge 0$, and words

 \ * \ $\,W^{(n_a)},\,\ldots\,, W^{(n_b)} \in $
    $(\Gamma_{\!M,\,{\sf pfl}} \cup \tau)^*\,$
    have been found such that for all $i, j \in [n_a, n_b]:$

\smallskip

 \ * \hspace{1.1cm} $E_{w_{t-1} \,\ldots\, w_1}$  $=$
    {\large $!$}$\bigcup_{r=n_a}^{ \ n_b} E_{W^{(r)}}\,$,

\smallskip

 \ * \hspace{1.1cm} ${\rm Dom}(E_{W^{(i)}}) \,\cap\, {\rm Dom}(E_{W^{(j)}})$
    $=$ $\0\,$ if $i \ne j\,$,

 \ * \hspace{1.1cm} $E_{W^{(i)}}$ is {\sf pfl} with
 $\,{\rm domC}(E_{W^{(i)}}) \subseteq \{0,1\}^m$ and
  \ $E_{W^{(i)}}(\{0,1\}^m) \subseteq \{0,1\}^i$.
 \ \ */  }

\medskip

{\footnotesize \sf /* Case 1: \ $w_t$ is {\sf pfl} \ */ }

if $w_t \in$ $\{a_0, a_3, a_4, a_5, a_6, a_7\} \cup \tau\,$ then$:$

\smallskip

\hspace{1,0cm} for $i \in [n_a,\,n_b]:$ 
 \ \ $W^{(i+\delta(w_t))}_{_{\rm\!new}}$ $:=$  $\,w_t \cdot W^{(i)}\,$;

\smallskip

\hspace{1,0cm} $n_b := n_b + \delta(w_t)\,$;
                  \ \ \ $n_a := n_a + \delta(w_t)$;

\smallskip

\hspace{1,0cm} for $j \in [n_a,\,n_b]:$
  \ \ $W^{(j)} := W^{(j)}_{_{\rm\!new}}\,$;
    
\medskip

{\footnotesize \sf /* Case 2: \ $w_t$ is {\sl not} plep; so, 
$w_t \in \{a_1, a_2, a_8\}$.
 \ Note that $\delta(a_{i,k}) = k$ \ if $\,i \in \{1,2,8\}$,
$\,k \in \{-1,0,1\}$ \ */
}
 
\smallskip

if $w_t = a_1$ \ \ \ {\footnotesize \sf /* Case 2.$a_1$: */ }
  
then:
 
\smallskip

\hspace{1,0cm} for $j \,\in\, [n_a\!-\!1, \ n_b\!+\!1]:$

\smallskip

\hspace{2,0cm} Let $W^{(j)}_{_{\rm\!new}} \in$
$(\Gamma_{\!M,\,{\sf pfl}} \cup \tau)^*\,$ be the word constructed in
Theorem \ref{THMwordforDisjU} for

\smallskip

\hspace{2,0cm} {\large $!$}$\bigcup$
$\big\{a_{1,k} \cdot W^{(i)} \ : \ i \in [n_a,\, n_b]\,$ and
$\,k \in \{-1,0,1\}\,$ such that $\, i + \delta(a_{1,k}) = j \big\}\,$;

\smallskip

\hspace{1,0cm} $n_b := n_b\!+\!1$; \ \ \ $n_a := n_a\!-\!1$;

\smallskip

\hspace{1,0cm} for $j \in [n_a,\,n_b]:$
            \ \ $W^{(j)} := W^{(j)}_{_{\rm\!new}}\,$;

\medskip

if $w_t = a_2$ \ \ \ {\footnotesize \sf /* Case 2.$a_2$: */ }

then:

\smallskip

\hspace{1,0cm} for $j \,\in\, [n_a\!-\!1, \ n_b\! +\! 1]:$

\smallskip

\hspace{2,0cm} Let $W^{(j)}_{_{\rm\!new}} \in$
$(\Gamma_{\!M,\,{\sf pfl}} \cup \tau)^*\,$ be the word constructed in
Theorem \ref{THMwordforDisjU} for

\smallskip

\hspace{2,0cm}  {\large $!$}$\bigcup$
$\big\{a_{2,k} \cdot W^{(i)} \ : \ i \in [n_a,\, n_b],$ and
$\, k \in \{-1,0,1\}\,$ such that $\,i + \delta(a_{2,k}) = j \big\}$;

\smallskip

\hspace{1,0cm} $n_b := n_b\!+\!1$; \ \ \ $n_a := n_a\!-\!1$;

\smallskip

\hspace{1,0cm} for $j \in [n_a,\,n_b]:$
      \ \ $W^{(j)} := W^{(j)}_{_{\rm\!new}}\,$;

\medskip

if $w_t = a_8$ \ \ \ {\footnotesize \sf /* Case 2.$a_8$: */ }

then:

\smallskip

\hspace{1,0cm} for $j \,\in\, [n_a, \ n_b\!+\! 1]:$

\smallskip

\hspace{2,0cm} Let $W^{(j)}_{_{\rm\!new}} \in$
$(\Gamma_{\!M,\,{\sf pfl}} \cup \tau)^*\,$ be the word constructed in
Theorem \ref{THMwordforDisjU} for

\smallskip

\hspace{2,0cm} {\large $!$}$\bigcup$
$\big \{a_{8,k} \cdot W^{(i)} \ : \ i \in [n_a,\, n_b]\,$ and
$\, k \in \{0,1\}\,$ such that $\,i + \delta(a_{8,k}) = j \big\}\,$;

\smallskip

\hspace{1,0cm} $n_b := n_b\!+\! 1$; \ \ \ {\footnotesize \sf
    /* $n_a$ is not changed */ }

\smallskip

\hspace{1,0cm} for $j \in [n_a,\,n_b]:$
      \ \ $W^{(j)} := W^{(j)}_{_{\rm\!new}}\,$;

\medskip

\noindent {\footnotesize \sf /* End of the ``for $t$'' loop; all the 
generators in $w$ have been applied. */ }

\medskip

\noindent output $\,\{ (W^{(i)}, \,i\!-\!m) \,:\, i \in [n_a, n_b] \}$.

\smallskip

\noindent {\footnotesize \sf [End, Sequential algorithm.]  }

\bigskip

\noindent Correctness of the Sequential Algorithm follows from the following
claim. 

\medskip

\noindent {\sf Claim.} If the {\em loop assertion} holds at the beginning of
an execution of the for-loop then it holds at the end of the loop for the new
$n_a$, $n_b$, and $(W^{(n_a)}, \,\ldots\,, W^{(n_b)})$.

\smallskip

\noindent {\sf Proof of Claim.} \ In case 1, 
$W^{(i)}$ and $w_t$ are {\sf pfl}; moreover, ${\rm domC}(E_{W^{(i)}})$ 
$\subseteq$  $\{0,1\}^m$ and $ \ E_{W^{(i)}}(\{0,1\}^m)$  $\subseteq$ 
$\{0,1\}^i$. Therefore the new $W^{(i+\delta(w_t))}$ is also {\sf pfl}, and 
$\,E_{W^{(i)}}(\{0,1\}^m)$ $\subseteq$ $\{0,1\}^{i+\delta(w_t)}$.

In case 2, $\,W^{(i)}$ and each fragment $w_{t,k}$ are {\sf pfl}, and
$W^{(i)}$ satisfies $\,{\rm domC}(E_{W^{(i)}}) \subseteq \{0,1\}^m\,$ and
  \ $E_{W^{(i)}}(\{0,1\}^m) \subseteq \{0,1\}^i$. Therefore,
$W^{(i)} \cdot w_{t,k}\,$ is {\sf pfl}; and the construction of
$W^{(j)}_{_{\rm\!new}}$ makes it {\sf pfl} (with input-length $m$ and
output-length $j$).

In both cases 1 and 2, at the end of the loop we have
$\,{\rm Dom}(E_{W^{(i_1)}}) \cap {\rm Dom}(E_{W^{(i_2)}})$ $=$ $\0\,$
(if $i_1 \ne i_2$); this follows from Lemma \ref{LEMemptyintersectDom}.

At the end of the loop with value $t$, $ \ E_{w_t w_{t-1}\,\ldots\, w_1}$ 
$\,=\,$ {\large $!$}$\bigcup_{s=n_a}^{ \ n_b} E_{W^{(s)}}$. This holds 
because the relation is true at the beginning of the loop (for
$w_{t-1} \,\ldots\, w_1$). Moreover, in case 1, $w_t$ is applied in the 
next loop, so the relation then holds for $\,w_t w_{t-1} \,\ldots\, w_1$. 
In case 2, the fragments of $w_t^!$ are applied (and the disjoint union of 
the fragments is equal to $w_t$); the operation 
``$W^{(j)}_{\rm\!new} := \ldots$'' permutes the functions and sorts them
according to output-length, but this does not change the disjoint union 
(which is commutative).
 \ \ \  \ \ \ {\small \sf [This proves the Claim.]}

\medskip

For the length we have $\,|W^{(i)}| \le$
${\rm maxlen}(\Gamma_{\!M}) \, |w|_{\Gamma_{\!M}}$ $+$
${\rm maxindex}_{\tau}(w)$, since in every execution of the loop,
$|W^{(i)}|$ grows by at most by length $c\,$ ($\,=\,$ 
${\rm maxlen}(\Gamma_{\!M})$), or
increases to some length $\,\le {\rm maxindex}_{\tau}(w)$.

\smallskip

The computation of final expressions $(W^{(n_a)}, \,\ldots\,, W^{(n_b)})$ 
for $w$ takes polynomial time. Indeed, the loop is iterated
$\,\le \ell\,$ ($= |w|_{\Gamma_{\!M} \cup \tau}$) times.
Moreover, each execution of the body of the loop takes polynomial time:
For case 1, it is straightforward that the execution takes linear time.
In case 2, the set 

 \ \ \ $W^{(j)}_{_{\rm new}}$ $\,=\,$ $\{a_{r,k} \cdot W^{(i)} \ :$
$ \ i \in [n_a,\, n_b], \ k \in \{-1,0,1\}\,$ (or $\in \{0,1\}$), 
such that $j = i + \delta(a_{r,k})\}\,$

\smallskip

\noindent has 3 elements (if $a_r = a_1$ or $a_r = a_2$), or 2 elements 
(if $a_r = a_8$); e.g.\ for $a_1$, $\,\delta(a_{1,k})$ ranges over 
$\{-1,0,1\}$, and $i$ is determined by $\delta(a_{1,k})\,$ (since 
$i = j - \delta(a_{1,k})$); so

\smallskip

 \ \ \ $W^{(j)}_{_{\rm new}}$  $\,=\,$ $\{ a_{1,-1} \cdot W^{(j+1)},$
 \ $a_{1,0} \cdot W^{(j)},$ \ $a_{1,1} \cdot W^{(j-1)}\}$.

\smallskip

\noindent In the formation of the set $W^{(j)}_{_{\rm new}}$, no new copy of
any existing $W^{(i)}$ is made; the exiting $W^{(i)}s$ are simply either 
left unchanged, or rearranged and concatenated with one letter $a_{r,k}$.
Thus, $W^{(j)}_{_{\rm new}}$ is constructed in linear time.

\bigskip

The next algorithm is an ${\sf AC}_2$ algorithm, with parallel for-loops. 
In $\,$``parallelfor $i \in [a,b]$ ... '',$\,$ the $\,b\!-\!a\!+\!1\,$ 
iterations of the body of the loop are independent and are executed 
simultaneously. See again the earlier descriptions of the sequential 
algorithm and the parallel algorithm. 

\bigskip

\noindent {\sf Parallel algorithm}

\smallskip

\noindent {\footnotesize \sf
/* Input: \  $w = w_{\ell} \,\ldots\, w_1\,$ (with $w_t \in$
$\Gamma_{\!M} \cup \tau$ for $t \in [1,\ell]$);

\noindent \ * \ Output: \ $\{(W^{(j)},\,j\!-\!m) : j \in [n_a,\, n_b]\}$
$ \ \subseteq \ $  $(\Gamma_{\!M,\,{\sf pfl}} \cup \tau)^* \x {\mathbb Z}$,
satisfying the conditions of the Theorem.
 \ \ */ }

\smallskip

\noindent parallelfor $t \in [1,\ell]:$ \ $x_t := w_t^!\,$;

\smallskip

\noindent $D = \lceil \log_2 \ell \rceil\,$;

\smallskip

\noindent if $\ell$ is not a power of 2 then:

parallelfor $t \in [\ell\!+\!1,\, 2^D]:$ \ $x_t := {\mathbb 1}\,$; 

\hspace{2,5cm} {\footnotesize \sf /* attach a prefix of $\,2^D\!-\!\ell\,$
                copies of ${\mathbb 1}$ to $x$, acting as the 
                identity function in compositions \ */}

$\ell := 2^D\,;$ \hspace{1,0cm}  {\footnotesize \sf /* Now, 
                               $\,\ell = |x| = 2^D$ */}

\smallskip

\noindent parallelfor $t \in [1,\ell]:$

 $n_{a,t} := m + \min(\delta_M(x_t))\,$;

 $n_{b,t} := m + \max(\delta_M(x_t))\,$;

\smallskip

\noindent for $d \in [0, D\!-\!1]:$

\medskip

{\footnotesize \sf
/* \ {\sl Loop assertion:}

 \ * \ For all $i \in [1,\, 2^{D-d}]:$ \ Assume that $x_i$ has been
    decomposed into an unambiguous union

 \ * \ $\,x_i =$ {\large $!$}$\bigcup$
$\big(x_i^{(j)} : j \in [n_{a,i}, n_{b,i}] \big)$, \ with
$x_i^{(j)}$ $\in$  $(\Gamma_{\!M,\,{\sf pfl}} \cup \tau)^*\,$,
and $\,{\rm imC}(E_{x_i^{(j)}}) \subseteq \{0,1\}^j$.

 \ * \ Assume that $x_i$ is represented by a set
 \ $\{(x_i^{(j)}, \ d_i^{(j)}) : j \in [n_{a,i}, n_{b,i}]\}$
$ \ \subseteq \ $  $(\Gamma_{\!M,\,{\sf pfl}} \cup \tau)^* \x {\mathbb Z}$,

 \ * \ with $d_i^{(j)} = j\!-\!m$; and $d_i^{(j)} = \delta(E_{x_i^{(j)}})\,$
 if $\,E_{x_i^{(j)}} \ne \theta$.
 \ \  */ }

\medskip

parallelfor $h \in [1,\, 2^{D-d-1}]:$

\smallskip

\hspace{0.7cm} $N_{a,h} := m + \min\{d_{2h}^{(j)} + d_{2h-1}^{(j')} : $
    $j \in [n_{a,2h},\,n_{b,2h}]$,
    $j' \in [n_{a,2h-1},\,n_{b,2h-1}] \}\,$;

\smallskip

\hspace{0.7cm} $N_{b,h} := m + \max\{d_{2h}^{(j)} + d_{2h-1}^{(j')} : $
    $j \in [n_{a,2h},\,n_{b,2h}]$,
    $j' \in [n_{a,2h-1},\,n_{b,2h-1}] \}\,$;

\smallskip

\hspace{0.7cm} Let $\,X_h = $
$\{(X_h^{(j)},\,D_h^{(j)}) : j \in [N_{a,h},\, N_{b,h}]\}$
$ \ \subseteq \ (\Gamma_{\!M,\,{\sf pfl}} \cup \tau)^* \x {\mathbb Z}\,$

\hspace{0.7cm} be the set constructed in Prop.\ \ref{LEMcircCompUnambig} for

\hspace{0.7cm} $x_{2h} \circ x_{2h-1}$ $\,=\,$
$\{(x_{2h}^{(j)},\,d_{2h}^{(j)}) : j \in [n_{a,2h},\,n_{b,2h}]\}$
$\circ$
$\{(x_{2h-1}^{(j)},\,d_{2h-1}^{(j)}) : j \in [n_{a,2h-1},\,n_{b,2h-1}]\}\,$;

\medskip

parallelfor $h \in [1,\, 2^{D-d-1}]:$

\hspace{0.7cm} $n_{a,h} := N_{a,h}\,$; \ \ \ $n_{b,h} := N_{b,h}\,$;

\smallskip

\hspace{0,7cm} parallelfor $j \in [n_{a,h}, n_{b,h}]:$

\hspace{1,50cm} $(x_h^{(j)},\,d_h^{(j)}) := (X_h^{(j)},\,D_h^{(j)})\,$;
           \ \ \ {\footnotesize \sf /* \ i.e., \ $x_h := X_h$ \ */ }

\smallskip

\noindent output $x_1$ \  \ \ \ {\footnotesize \sf  /* $\,=\,$
  $\{(x_1^{(j)},\,d_1^{(j)}) : j \in [n_{a,1}, n_{b,1}]\}$ \ */} .

\medskip

\noindent {\footnotesize \sf [End, Parallel algorithm.] }

\bigskip

The parallel algorithm and the sequential algorithm produce the same output.
Indeed, the two algorithms differ only by the grouping of subexpressions,
and composition is associative and distributes over union.

\smallskip

The parallel algorithm uses $D = \lceil \log_2 \ell \rceil$ sequential steps,
where $\ell = |w|_{\Gamma_{\!M} \cup \tau}$. In step $\,d \in [0, D\!-\!1]\,$
there are $2^{D-d-1}$ simultaneous compositions of the form 
$\,x_{2h} \circ x_{2h-1}\,$ (in ``parallelfor $h \in [1, 2^{D-d-1}]$''); 
each such composition is computed in log-space (by Prop.\ 
\ref{LEMcircCompUnambig}). 
Since log-space is contained in ${\sf AC}_1$, each parallel step $d$ is in 
${\sf AC}_1$. Thus  the parallel algorithm is in ${\sf AC}_2$.

\medskip

\noindent (3) The proof for ${\cal RM}^{\sf fin}$ is similar to the proof
for $M_{2,1}$, using a generating set of the form $\Gamma \cup \tau$, where 
$\Gamma$ is finite.
Just as in the proof of (2), each generator $\gamma \in$  $\Gamma$  $\minus$
${\sf pfl}{\cal RM}^{\sf fin}$ can be decomposed as into an unambiguous union
of elements of ${\sf pfl}{\cal RM}^{\sf fin}$ with disjoint domains. From 
there on the reasoning is the same as in (2).
 \ \ \  \ \ \ $\Box$

\bigskip

\bigskip

\noindent {\bf Comments on Theorem \ref{THM_Mtoplep}}

\medskip

\noindent {\bf 1. Theorem \ref{THM_Mtoplep}(2) and the function-emptiness 
problem:} 

\smallskip

A subtle aspect of Theorem \ref{THM_Mtoplep}(2)(3) is that it produces a
correct output in polynomial time (in fact in complexity ${\sf AC}_2$) 
regardless of whether $E_w = \theta$, which is 
{\sf coNP}-complete Prop.\ \ref{deltaComplexity}(1). 

\bigskip

\noindent {\bf 2. Equivalent forms of the input-length,
and of the depth of a word in ${\cal RM}^{\sf fin}$}

\smallskip

\noindent Based on the unambiguous union decomposition of the elements of
${\cal RM}^{\sf fin}$ we obtain the following characterization of 
input-length and depth.

\begin{pro} \label{PROPioLengthMUCirc} {\bf (input-length and depth for
${\cal RM}^{\sf fin}$ and union of partial circuits).}

\noindent Let $\ell_{\rm in}(.)$ be the input-length function (Def.\ 
{\rm \ref{DEF_IOlengthRMfin}} and Lemma {\rm \ref{EquivDefInputLen}}), and 
let $\,{\sf depth}(.)$ be the depth function ({\rm Def.\ 
\ref{DEFdepthM21wx}}).
Then for any $w \in $ $(\Gamma \cup \tau)^*$ and its decomposition
$\,!$$\bigcup_{i=n_a}^{ \ \ n_b} W^{(i)}\,$, given in {\rm Theorem
\ref{THM_Mtoplep}(2)}, we have:

\medskip

 \ \ \  \ \ \ $\ell_{\rm in}(w) \ = \ \ell_{\rm in}(W^{(i)})$,
 \ \ for every $\,i \in [n_a, n_b]$;

\medskip

 \ \ \  \ \ \ ${\sf depth}(w) \ = \ $
$\max\{\,{\sf depth}(W^{(i)}) \,:\, i \in \delta_M(w)\}$.
\end{pro}
For an element $w \in (\Gamma \cup \tau)^*$ with $\,E_w \in$
${\cal RM}^{\sf fin} \minus {\sf pfl}{\cal RM}^{\sf fin}$, there is more
than one output-length; the output-lengths for the set
$\,\{ \ell_{\rm in}(w) + d \,:\,$ $d \in \delta_M(w)\}$.

\medskip

\noindent {\sc Proof.} The property of $\ell_{\rm in}(w)$ is given in the
first item of Theorem \ref{THM_Mtoplep}(2.2).
The depth property is proved in the third item of Theorem 
\ref{THM_Mtoplep}(2.2).
 \ \ \  \ \ \  \ \ \  $\Box$

\section{Appendix: The function emptiness problem, and computing 
$\delta(.)$ in ${\sf pfl}{\cal RM}^{\sf fin}$ }

We consider the following three problems about
${\sf pfl}{\cal RM}^{\sf fin}$ with generating set $\Gamma_{\sf\!pfl}$ 
$\cup$  $\tau$, with $\Gamma_{\sf\!pfl}$ finite. 
The function emptiness problem also applies to ${\cal RM}^{\sf fin}$ with
generating set $\Gamma_{\sf\!fin} \cup \tau$, with $\Gamma_{\sf\!fin}$ 
finite.

\bigskip

\noindent {\small \rm (1)} {\bf The function emptiness problem} is the following decision problem. 

\smallskip

\noindent Input: \ $w \in (\Gamma_{\sf\!pfl} \cup \tau)^*$, or
$\,w \in (\Gamma_{\sf\!fin} \cup \tau)^*$.

\smallskip

\noindent Question:
 \ Is $\,E_w = \theta$ \ (the empty function, as an element of
${\rm pfl}{\cal RM}^{\sf fin}$, respectively of ${\cal RM}^{\sf fin}$)?

\smallskip

In ${\sf pfl}{\cal RM}^{\sf fin}$, this is equivalent to the question:
 \ Is $\,\delta(w)\,$ defined?

In ${\cal RM}^{\sf fin}$ it is equivalent to the question:
 \ Is $\,\delta_{\!M}(w) = \0\,$?

\bigskip

\noindent {\small \rm (2)} {\bf The computation problem for $\delta$ in the
 non-empty case} is the following premiss problem.

\smallskip

\noindent Input: \ $w \in (\Gamma_{\sf\!pfl} \cup \tau)^*$.

\smallskip

\noindent Premiss: \ $E_w \ne \theta$.

\smallskip

\noindent Output: \ $\delta(E_w) \ $ (a number in ${\mathbb Z}$, written in
sign-magnitude binary notation).

\bigskip

\noindent {\small \rm (3)} {\bf The computation problem for $\delta$}
is the following input-output problem.

\smallskip

\noindent Input: \ $w \in (\Gamma_{\sf\!pfl} \cup \tau)^*$.

\smallskip

\noindent Output: \ $\delta(E_w) \ $ (a number in ${\mathbb Z}$, written 
in sign-magnitude binary notation).

\medskip

The {\sf coNP}-completeness of function emptiness problem for $M_{2,1}$ over 
$\Gamma \cup \tau$ was already proved in \cite[{\small Prop.\ 6.2}]{BiLR}.
In Prop. \ref{deltaComplexity}(1) this is now proved for 
${\cal RM}^{\sf fin}$ and ${\sf pfl}{\cal RM}^{\sf fin}$.

\smallskip

\begin{pro} \label{deltaComplexity} {\bf (complexity of $\delta(.))$}.

\smallskip

\noindent {\small \rm (1)} The function emptiness problem for
${\sf pfl}{\cal RM}^{\sf fin}$, or for ${\cal RM}^{\sf fin}$, is
{\sf coNP}-complete (with respect to log-space many-one reduction).

\smallskip

\noindent {\small \rm (2)} The premiss problem for computing $\delta(.)$ in
${\sf pfl}{\cal RM}^{\sf fin}$ belongs deterministic log-space.
(However, checking the premiss is the negation of the function emptiness 
problem, hence it is {\sf NP}-complete.)

\smallskip

\noindent {\small \rm (3)} The computation problem for $\delta(.)$ in 
${\sf pfl}{\cal RM}^{\sf fin}$ is {\sf coNP}-hard.
\end{pro}
{\sc Proof.} (1) The function emptiness problem is in {\sf coNP}. Indeed for
any $w$ in $(\Gamma_{\sf\!pfl} \cup \tau)^*$ or 
$(\Gamma_{\sf\!fin} \cup \tau)^*$ we have: \ $E_w = \theta$ \ iff
 \ $(\forall x \in \{0,1\}^m)[ E_w(x)$ is undefined]. 
Here we can choose $m$ to be $\,m = $
${\rm maxlen}(\Gamma_{\sf\!pfl}) \,|w|_{\Gamma_{\sf\!pfl} \cup \tau}$
$+$  ${\rm maxindex}_{\tau}(w)\,$; \ hence $m$ is easily computed.
And for a given $w$ and $x$, it can be checked in polynomial time (in terms
of $|w| + |x|$) whether $E_w(x)$ is undefined.

\smallskip

To show {\sf coNP}-hardness we first consider the following problem:
The input is an acyclic circuit $C$ with a single output wire and $m$ input
wires, for any $m \in {\mathbb N}$; and the question is whether for all 
inputs in $\{0,1\}^m$ the output of $C$ is 0. Let us call this problem 
NonSat; it is the negation of the satisfiabilty problem for acyclic circuits,
and it is well-known to be {\sf coNP}-complete (see e.g.\ 
\cite{Papadim, JESavage}).

Next, NonSat is reduced to the function emptiness problem. Given a circuit
$C$, we find ${\sf W}(C) \in$ $(\Gamma_{\sf\!tfl} \cup \tau)^*$ such that 
$E_{{\sf W}(C)}$ is the input-output function of $C$; by Theorem 
\ref{circ2lepM}, ${\sf W}(C)$ can be computed in log-space.  Thus $C$ has 
output 0 on all inputs in $\{0,1\}^m\,$ iff $\,[{\sf W}(C)](x) = 0\,$ for 
all $x \in \{0,1\}^m$. Next, we define $\zeta \in$ 
${\rm pfl}{\cal RM}^{\sf fin}$ as follows: 
$ \ \zeta(1z) = 1z$, and $\,\zeta(0z)\,$ is undefined (for all
$z \in \{0,1\}^*$).
Let $w_0 \in$ $(\Gamma_{\sf\!pfl} \cup \tau)^*$ be such that 
$\,E_{w_0} = \zeta$.
Then $C$ is a yes-case of NonSat iff $\,E_{w_0\,{\sf W}(C)} = \theta$.

\smallskip

(2) Let $\,w = w_k\,\ldots\,w_1$ with $w_t \in \Gamma_{\sf\!pfl} \cup \tau$
for $t \in [1,k]$.  By Prop.\ \ref{PROPiolDiffpartial}(3), $\delta(w) =$
$\sum_{t=1}^k \delta(w_t)\,$ if $E_w \ne \theta$. And for $\gamma \in$
$\Gamma_{\sf\!pfl} \cup \tau$, $\,\delta(\gamma)$ is bounded by constants 
from above and from below, since $ \Gamma_{\sf\!pfl}$ is finite, and
$\delta(\tau_{i,j}) = 0$ for all $\tau_{i,j} \in \tau$.
We now use the well-known and easy fact that the sum of $k$ bounded integers 
can be computed in space $O(\log k)$.

\smallskip

(3) This follows immediately from (1) and (2).
 \ \ \  \ \ \ $\Box$

\bigskip

\noindent {\bf Remark:} \ We have for all $w \in$ 
$(\Gamma_{\sf\!fin} \cup \tau)^*$: $\,E_w = \theta\,$ in 
${\cal RM}^{\sf fin}$ \ iff \ $E_w = \theta\,$ in $M_{2,1}$. 
Indeed, $E_w = \theta$ iff ${\rm domC}(E_w) = \0$, and this holds whether
$E_w$ is taken in ${\cal RM}^{\sf fin}$ or in $M_{2,1}$.
Hence the above Proposition also applies to $M_{2,1}$ and 
${\sf plep}M_{2,1}$.


\bigskip



\medskip

\noindent {\small birget@camden.rutgers.edu }

\end{document}